\documentclass[3p]{elsarticle}

\usepackage{array}
\usepackage{color}
\usepackage{xcolor,colortbl}
\usepackage{natbib}
\usepackage{tabularx}
\usepackage{graphicx}
\usepackage{amsmath}
\usepackage{amssymb}
\usepackage{amsthm}
\usepackage{empheq}
\usepackage{multirow}
\usepackage{subfig}
\usepackage{moreverb}
\usepackage{nicefrac} 
\usepackage{dsfont}
\usepackage{bm}
\usepackage{multirow}
\usepackage{soul}
\usepackage{lmodern}

\newcommand\BibTeX{{\rmfamily B\kern-.05em \textsc{i\kern-.025em b}\kern-.08em
T\kern-.1667em\lower.7ex\hbox{E}\kern-.125emX}}

\usepackage{hyperref} 

\newcommand{\x}{\mbf{x}}

\newcommand{\mbf}[1]{\mathbf{#1}}           %

\newcommand{\Q}{\mathbf{Q}}
\renewcommand{\S}{\mathbf{S}}
\renewcommand{\u}{\mathbf{u}}
\newcommand{\w}{\mathbf{w}}
\newcommand{\q}{\mathbf{q}}
\newcommand{\F}{\mathbf{F}}
\newcommand{\f}{\mathbf{f}}
\newcommand{\g}{\mathbf{g}}
\newcommand{\h}{\mathbf{h}}

\newcommand{\halb}{\frac{1}{2}}

\newcommand{\be}{\begin{equation}}
\newcommand{\ee}{\end{equation}}
\newcommand{\bdm}{\begin{displaymath}}
\newcommand{\edm}{\end{displaymath}}

\newcommand{\bea}{\begin{eqnarray} }
\newcommand{\eea}{\end{eqnarray} }

\newcommand{\aposteriori}{\textit{a posteriori} }

\newcommand{\dev}{\textnormal{dev}} 

\newcommand{\AAA}{{\boldsymbol{A}}}

\newcommand{\vv}{{\mathbf{v}}}

\newcommand{\JJ}{{\mathbf{J}}}
\newcommand{\QQ}{{\mathbf{Q}}}
\newcommand{\QV}{{\mathbf{V}}}

\newcommand{\tr}{\textnormal{tr}}

\DeclareMathOperator{\de}{d\!}
\renewcommand{\vec}[1]{\bm{\mathrm{#1}}}
\newcommand{\up}[1]{\ensuremath{\mathrm{#1}}}
\newcommand{\od}[2]{\frac{\de{#1}}{\vphantom{l^l}\de{#2}}}
\newcommand{\abs}[1]{\left|#1\right|}

\newcommand{\Qe}{\vec{Q}_\up{e}}
\newcommand{\Qs}{{\vec{Q}^\ast}}
\newcommand{\Js}{{\vec{J}^\ast}}
\newcommand{\Ss}{{\vec{S}^\ast}}
\newcommand{\Bs}{{\vec{B}^\ast}}
\newcommand{\Cs}{{\vec{C}^\ast}}
\newcommand{\ts}{{t^\ast}}

\newcommand\trasp[1]{{#1}^{\mathsf T}}

\newfont{\numerikEleven}{ecrm1000}
\newfont{\numerikTen}{cmss10}
\newfont{\numerikNine}{cmss9}
\newfont{\numerikEight}{cmss8}

\hypersetup{
hidelinks,
    colorlinks=true,                          
    linkcolor=black, % Couleur des liens internes
    citecolor=black, % Couleur des num�ros de la biblio dans le corps
    urlcolor=black,
    allcolors=black} % Couleur des url

\journal{Journal of Computational Physics}

\begin{document} 

\begin{frontmatter} 
%-------------------------------------------------------
% TITLE
\title{\textbf{Space-time adaptive ADER discontinuous Galerkin schemes for nonlinear hyperelasticity with material failure}} 
%-------------------------------------------------------
%-------------------------------------------------------
% AUTHORS

\author[UniTN]{\small Maurizio Tavelli}
\ead{m.tavelli@unitn.it}

\author[UniTN]{\small Simone Chiocchetti}
\ead{simone.chiocchetti@unitn.it}

\author[UniTN,NSC]{\small Evgeniy Romenski}
\ead{evrom@math.nsc.ru}

\author[LMU]{\small Alice-Agnes Gabriel}
\ead{gabriel@geophysik.uni-muenchen.de}

\author[UniTN]{\small Michael Dumbser$^{*}$}
\ead{michael.dumbser@unitn.it}
\cortext[cor1]{Corresponding author}

%-------------------------------------------------------
% INSTITUTIONS
\address[UniTN]{Laboratory of Applied Mathematics, University of Trento, Via Mesiano 77, 38123 Trento, Italy}
\address[NSC]{{Sobolev Institute of Mathematics, 4 Acad. Koptyug Avenue, 630090 Novosibirsk, Russia}}
\address[LMU]{Ludwig-Maximilians Universit\"at M\"unchen (LMU), Theresienstr. 41, 80333 M\"unchen, Germany}
%-------------------------------------------------------

%-------------------------------------------------------
% ABSTRACT
\begin{abstract}
We are concerned with the numerical solution of a \textit{unified} first order hyperbolic formulation of continuum 
mechanics that goes back to the work of Godunov, Peshkov and Romenski \cite{GodRom1972,GodRom2003,PeshRom2014} (GPR model) 
and which is an extension of nonlinear hyperelasticity that is able to describe simultaneously nonlinear elasto-plastic 
solids at large strain, as well as viscous and ideal fluids. The proposed governing PDE system also contains the effect 
of heat conduction and can be shown to be symmetric and thermodynamically compatible, as it obeys the first and second 
law of thermodynamics. 

In this paper we extend the GPR model to the simulation of nonlinear dynamic rupture processes, which can be achieved by adding  
an additional scalar to the governing PDE system. This extra parameter describes the material damage and is governed 
by an advection-reaction equation, where the stiff and highly nonlinear reaction mechanisms depend on the ratio  
of the local equivalent stress to the yield stress of the material. The stiff reaction mechanisms are integrated in time 
via an efficient exponential time integrator. Due to the multiple spatial and temporal scales involved in the problem of crack 
generation and  propagation, the model is solved on space--time adaptive Cartesian meshes using  
high order accurate discontinuous Galerkin finite element schemes endowed with an a posteriori subcell finite volume limiter. 

A \textit{key feature} of our new model is the use of a twofold \textit{diffuse interface approach} that allows the 
cracks to form anywhere and at any time, independently of the chosen computational grid, which is simply adaptive 
Cartesian (AMR). This is substantially different from many fracture modeling approaches that need to resolve discontinuities explicitly, 
such as for example dynamic shear rupture models used in computational seismology, 
where the geometry of the rupture fault needs to prescribed \textit{a priori}. We furthermore make use of a scalar 
volume fraction function $\alpha$ that indicates whether a given spatial point is inside the solid ($\alpha=1$) or 
outside ($\alpha=0$), thus allowing the description of solids of arbitrarily complex shape. 

We show extensive numerical comparisons with experimental results for stress-strain diagrams of different real materials 
and for the generation and propagation of fracture in rocks and pyrex glass at low and high velocities. Overall, a very 
good agreement between numerical simulations and experiments is obtained. The proposed model is also naturally able to describe \textit{material fatigue}. 
\end{abstract}
%-------------------------------------------------------

%-------------------------------------------------------
% KEY WORDS
%-------------------------------------------------------
\end{frontmatter}

%-----------------------------------
% CONTENTS
%  This will deseapear in the submitted version
%\tableofcontents
%------------------------------------

%=========================================================================
%==========         I N T R O D U C T I O N
% 
\section{Introduction} \label{sec:introduction}

Nonlinear solid mechanics with large deformations and material failure is of fundamental importance in  mechanical and civil engineering for the computation and design of structures under dynamic loads. 
Particularly relevant for engineering applications are crack formation in intact materials and the phenomenon
of material fatigue. Also in geophysics nonlinear large deformation solid mechanics is relevant for the 
description of dynamic rupture processes in earthquakes.
One major difficulty in numerical modeling of crack propagation is the introduction of strong discontinuities in the displacement field in the vicinity of the crack. Numerical methods which require the discontinuities to be exactly resolved by the mesh generally require the geometry of the fractures 
to be predefined and thus typically 
only permit small deformations. Alternative methods which allow representing discontinuities at the sub-element level, such as the eXtended Finite
Element Method (XFEM) \cite{Moes1999}, introduce singularities when an interface intersects a cell vertex, and can become difficult to implement in an efficient manner in 3D. 

 While solid mechanics is usually formulated more  naturally in Lagrangian coordinates, it was shown in \cite{GodunovRomenski72} that a mathematical formulation of solid mechanics is also possible in Eulerian coordinates. This formulation was put in the more general  mathematical framework of symmetric hyperbolic and thermodynamically compatible (SHTC) systems in Godunov-type form \cite{God1961} in a series of papers \cite{God1972MHD,Rom1998,GodRom1995,GodRom2003}. Similar hyperelastic models were also used later in \cite{Barton2009, BartonRom2010, Barton2012, PeshRom2014, DumbserPeshkov2016, Iollo2017, AbateIolloPuppo, Jackson2019, Jackson2019a}, including strain hardening and plastic deformations. 
In order to model complex geometries, diffuse interface methods can be employed, where a scalar volume fraction represents the occupation of each control volume by the solid material. The most important contributions concerning  
diffuse interface models for the description of compressible hyperelastic multi-material flows have been made in a series of papers by Favrie and Gavrilyuk and collaborators, see 
\cite{Gavrilyuk2008,FavrieGavrilyukSaurel,FavrGavr2012,Ndanou2014,NdanouFavrieGavrilyuk,Barton2019}.

The main objectives of this paper are to show that nonlinear elasto-plastic material behaviour under large strains with material damage can be described via nonlinear hyperelasticity \cite{GodunovRomenski72} and the framework of SHTC systems in Godunov form \cite{God1961,Rom1998,SHTC-GENERIC-CMAT} and that the resulting governing PDE system can be robustly solved with high order accurate discontinuous Galerkin (DG) finite element schemes with \textit{a posteriori} subcell finite volume limiter, 
see \cite{Dumbser2014}. 

\textcolor{black}{Currently the generally accepted model of continual failure and rupture is the phase-field model, which was proposed by Francfort and Marigo in \cite{Francfort1998} for brittle materials and served as the basis for further theoretical and numerical developments. The phase-field model consists in representing the crack as a diffuse interface between intact zones of the medium using the order parameter field and is usually formulated as a minimization of the energy functional and is solved by the finite element method. Numerous generalizations of the model \cite{Francfort1998} for the brittle and ductile fracture and the fracture of porous material have been developed and their description can be found, for example, in \cite{Ambati2015, Ambati2015Rev, Miehe2016, Carollo2018, Spetz2020} and references therein. Note that all phase-field models are well suited for \textit{quasi-static} problems, and the corresponding partial differential equations contain second-order (parabolic) terms, which is different from the first order hyperbolic PDE system proposed in this paper that allows to construct a \textit{rate-dependent} and fully dynamic model of material failure, including also \textit{material fatigue}.} 

\textcolor{black}{
Numerical modeling of secondary cracks during dynamic fracture propagation is also an important problem to consider in earthquake physics \cite{Kikuchi1975}: high stress concentrations at the earthquake rupture tip feedback with off-fault anelastic processes such as continuum damage \cite{Lyakhovsky2005,Bhat2012,Thomas2017,Thomas2018}, explicit tensile and shear fracturing \cite{Yamashita2000,Dalguer2003,Okubo2019,Okubo2020} or off-fault plasticity \cite{Andrews2005,Gabriel2013}.
}  

In this paper we first present an extended diffuse interface version of the Godunov-Peshkov-Romenski (GPR) 
{\it hyperbolic} model of continuum mechanics \cite{PeshRom2014,DumbserPeshkov2016,BoscheriDumbser2016,HypoHyper2} that includes an additional evolution equation for a scalar material \textcolor{black}{damage variable} $\xi$, where $\xi=0$ indicates the fully intact material, while $\xi=1$ represents a fully damaged material. The dynamics of the material damage is driven by a reaction-type source term (similar to reaction source terms in chemistry), that depends on the ratio of equivalent stress to yield stress. This idea goes back to \cite{Resnyansky2003,Romenskii2007,Resnyansky2009}, where the dynamics of failure waves was discussed in the context of hyperelasticity. 
To account for nonlinear elasto-plastic materials, the GPR model includes a strain relaxation mechanism, see \cite{PeshRom2014,DumbserPeshkov2016}. We also suggest some simple mixture rules for the computation of the Lam\'e parameters and the strain relaxation time scale of the material as a function of the \textcolor{black}{damage variable} $\xi$.  
\textcolor{black}{All above gives us a highly flexible model that allows us to simulate large deformations with rate- and temperature- dependent failure of brittle and ductile materials.}
Note that an anisotropic thermodynamically consistent damage model based on the hyperelastic Godunov-Romenski approach has already been proposed in \cite{Barton2016}. In \cite{Barton2016} a damage  deformation gradient was introduced in addition to an elastic and a plastic deformation gradient.  

The numerical methods employed in this paper essentially rely on the framework of Godunov-type finite volume schemes, see e.g. \cite{godunov,hll,munz91,Toro:1994,roe,woodwardcol84,eno,toro3,toro4,toro-book}, which are then used as \textit{a posteriori}  subcell limiter \cite{Dumbser2014,Sonntag,Sonntag2} of a high order discontinuous Galerkin finite element method \cite{cbs1,cbs2,cbs3,cbs4} on space-time adaptive Cartesian AMR meshes \cite{Berger-Oliger1984,berger85,Berger-Colella1989,Bell1994,Berger-Leveque1998,Zanotti2015a}. The concept of \textit{a posteriori} limiting was introduced in a series of papers by Clain and Loub\`ere and co-workers in \cite{CDL1,CDL2,CDL3,ADER_MOOD_14} and relies on the \textit{a posteriori} detection of troubled cells and subsequent local order reduction of the scheme, in order to ensure essentially non-oscillatory behaviour and positivity of the numerical solution. For more details, the reader is referred to the above references. To the best knowledge of the authors, this is the first time that high order discontinuous Galerkin finite element schemes are applied to a nonlinear hyperelastic model of elasto-plastic media with material failure and diffuse interface representation of the geometry of the medium. 

The rest of this paper is organized as follows. In Section \ref{sec:model} we introduce and discuss the underlying mathematical model. In Section \ref{sec:ader} we present the numerical method employed in this paper to solve the governing partial differential equations. Our approach is based on high order discontinuous Galerkin (DG) finite element schemes with \textit{a posteriori} subcell finite volume limiting on space-time adaptive meshes (AMR). In order to deal with the very stiff reaction source terms in the evolution equation of the \textcolor{black}{damage variable}, an accurate and efficient ODE solver based on exponential time integration is presented. The time integration method is validated for different material behaviours (brittle, ductile and fatigued materials) against standard solvers like LSODA. In Section \ref{sec:results} we present computational results for a set of different test cases containing crack generation, including also comparisons with experimental data for some of the test problems. The paper closes with Section \ref{sec:conclusion}, where concluding remarks and an outlook to future research are given.

\section{Mathematical model} 
\label{sec:model} 

\subsection{Governing partial differential equations} 
\label{sec:pde}

A diffuse interface formulation for moving nonlinear solids of arbitrary geometry and at large strain is given by the following PDE system in Eulerian coordinates: 
\begin{subequations}\label{eqn.GPR}
    \begin{align}
& \frac{\partial \alpha}{\partial t}+v_k \cdot \frac{\partial \alpha}{\partial x_k}=0,  
\label{eqn.alpha} \\ 
&\frac{\partial \rho}{\partial t}  + \frac{\partial (\rho v_k )}{\partial x_k} =0, 
\label{eqn.mass} \\
&\frac{\partial (\alpha \rho v_i)}{\partial t}  +\frac{\partial \left(  \alpha \rho v_i v_k + \alpha p \delta_{ik} - \alpha \sigma_{ik}  \right) }{\partial x_k}   =0, 
\label{eqn.momentum} \\ 
&\frac{\partial A_{ik}}{\partial t} +\frac{\partial A_{im}v_m}{\partial x_k} + 
v_m \left( \frac{\partial A_{ik}}{\partial x_m} -\frac{\partial A_{im}}{\partial x_k}  \right) = 
-\frac{1}{\theta_1(\tau_1)} E_{A_{ik}}, 
\label{eqn.A} \\  
& \frac{\partial J_k}{\partial t} + \frac{\partial \left( v_m J_m + T \right) }{\partial x_k} + 
v_m \left( \frac{\partial J_{k}}{\partial x_m} -\frac{\partial J_{m}}{\partial x_k} \right) = - \frac{1}{\theta_2(\tau_2)} E_{J_k},  
\label{eqn.J} \\ 
& \frac{\partial \xi}{\partial t} + v_k \cdot \frac{\partial \xi}{\partial x_k} = - \theta E_\xi, 
\label{eqn.xi} \\ 
& \frac{\partial \rho S }{\partial t} + \frac{\partial \left( \rho S v_k + \rho E_{J_k} \right)}{\partial x_k} = 
\frac{\rho}{T} \left( \frac{1}{\theta_1} E_{A_{ik}} E_{A_{ik}} + \frac{1}{\theta_2} E_{J_k} E_{J_k}  + \theta E_\xi E_\xi \right) \geq 0,  
\label{eqn.entropy} \\ 
& \frac{\partial \rho {E}}{\partial t}  + \frac{ \partial \left( v_k \rho {E} +v_i(p\delta_{ik} - \sigma_{ik}) \right) }{\partial x_k}=0.  
\label{eqn.energy} 
\end{align}
\end{subequations} 

Throughout this paper we make use of the Einstein summation convention over repeated indices. 
Here \eqref{eqn.alpha} is the evolution equation for the colour function $\alpha$ that is needed in the diffuse interface approach (DIM) as introduced in \cite{Tavelli2019,FrontierADERGPR} for the description of solids of arbitrary geometry; \eqref{eqn.mass} is the mass conservation law and $\rho$ is the material density; \eqref{eqn.momentum} is the momentum conservation law and $v_i$ is the velocity field; \eqref{eqn.A} is the evolution equation for the distortion field $\mathbf{A}=A_{ik}$ (basis triad \footnote{\textcolor{black}{Note that $A_{ik}$ transforms as a set of three vectors and not as a rank two tensor under coordinate transformations. In the absence of strain relaxation source terms ($\tau_1 \to \infty$) $A_{ik}$ is the deformation gradient of the material.}}); \eqref{eqn.J} is the evolution equation for the specific thermal impulse $J_k$ constituting the heat conduction in the medium via a hyperbolic (non Fourier--type) model; \eqref{eqn.xi} is the evolution equation for the material \textcolor{black}{damage variable} $\xi\in[0,1]$, where $\xi=0$ indicates fully intact material and $\xi=1$ fully damaged material. Finally, \eqref{eqn.entropy} is the entropy inequality and \eqref{eqn.energy} is the energy conservation law.
Other thermodynamic parameters are defined via the total energy potential $E=E(\rho,S,\vv,\AAA,\JJ,\xi)$: 
$\vec{\Sigma} = {\Sigma}_{ik} = - p\delta_{ik} + \sigma_{ik}$ is the total stress tensor ($\delta_{ik}$ is the Kronecker delta); $p = \rho^2 E_\rho$ is the contribution to the stress tensor due to volume deformations;
$\sigma_{ik} =  -\rho A_{ji} E_{A_{jk}} + \rho J_i E_{J_k}$ is the contribution to the stress tensor due to shear and thermal stress, and $T = E_S$ is the temperature. 
\textcolor{black}{Note that in system \eqref{eqn.GPR} we only use a simplified diffuse interface approach, which completely neglects the dynamics of the air surrounding the solid medium. In our model the solid volume fraction $\alpha$ is only used to locate and track the shape of the moving solid, see \cite{Tavelli2019}. In comparison, in the work of Favrie and Gavrilyuk and collaborators \cite{FavrieGavrilyukSaurel,FavrGavr2012,NdanouFavrieGavrilyuk}, real multi-phase flows of compressible solids embedded in compressible fluids were considered. } 

The dissipation in the medium includes two relaxation processes: the strain relaxation (or shear stress relaxation) characterized by the scalar function $\theta_1(\tau_1) > 0$ depending on the relaxation time $\tau_1$ and the heat flux relaxation characterized by 
$\theta_2(\tau_2) > 0$, depending on the relaxation time $\tau_2$. 

Further to the above evolution equations we need to add the governing PDE for the local Lam\'e parameters $\lambda$ and $\mu$, material yield stress $Y_0$ and other material constants. For small deformations it is sufficient to assume them as constant in time, but for large deformations those parameters have to move with the local velocity field. The resulting PDEs for the  material parameters therefore read: 
\begin{equation}
\frac{\partial \lambda}{\partial t}+v_k \cdot \frac{\partial \lambda}{\partial x_k}=0, 
\qquad
\frac{\partial \mu}{\partial t}+v_k \cdot \frac{\partial \mu}{\partial x_k}=0,
\qquad
\frac{\partial Y_0}{\partial t}+v_k \cdot \frac{\partial Y_0}{\partial x_k}=0. 
\label{eqn1b} 
\end{equation}
If $\alpha=1$, then it can be immediately seen that the GPR model \eqref{eqn.GPR} is thermodynamically compatible as soon as the total energy is defined as a function of parameters of state \cite{Rom1998,GodRom2003,PeshRom2014,DumbserPeshkov2016}, since it 
satisfies the first and second principle of thermodynamics, see \eqref{eqn.energy} and \eqref{eqn.entropy}, 
which are the total energy conservation law and the entropy inequality, respectively. 
Note that system \eqref{eqn.GPR} is overdetermined, that is the number of equations is one more than the number of unknown parameters of state. In fact, as is always the case for SHTC systems, the energy equation \eqref{eqn.energy} is a consequence of all other equations and we underline that in the numerical test cases presented in this paper, we solve the \textit{energy conservation law}  
\eqref{eqn.energy} instead of the entropy inequality \eqref{eqn.entropy}.
%, but from the pure point of view of the  formulation of the mathematical model, the entropy should be considered among the vector of %unknowns and the 
%extra energy conservation law can be obtained as a \textit{consequence} of the other equations.  

Furthermore, in order to close the system one must specify the total energy potential as a function of the other state variables, i.e. $ E = E(\rho,S,\vv,\AAA,\JJ,\xi) $. This  potential then generates all terms in the fluxes and source terms by means of its partial derivatives with respect to the state variables. Therefore, as already discussed in \cite{DumbserPeshkov2016}, the energy specification is the key step 
in the model formulation. 

Here, we make the choice $E=E_1+E_2+E_3$, decomposing the energy into a contribution from the microscale $E_1$, the mesoscale $E_2$ and the macroscale $E_3$. The individual contributions read as follows:  
\begin{equation}
E_1 =  \frac{K}{2 \rho_0} \left( 1-\frac{\rho}{\rho_0}\right)^2 + c_v T_0\left(\frac{\rho}{\rho_0} \right)\left(e^{S/c_v}-1 \right), \label{E1a}
\end{equation} 
which is the equation of state (EOS) of the medium, where $\rho_0$ is the reference mass density and $K=\lambda+\frac{2}{3}\mu$ is the bulk modulus expressed in terms of the two Lam\'e parameters $\lambda$ and $\mu$, which depend on the damage variable $\xi$, $c_v$ is the heat capacity at constant volume and $T_0$ is a reference temperature. 
The mesoscale energy in our model reads as 
\begin{equation} 
E_2 =  \frac{1}{4} c_s^2 \mathring{G}_{ij} \mathring{G}_{ij} + \frac{1}{2} c_h^2 J_i J_i, \label{E2a}
\end{equation} 
where $c_s=\sqrt{\frac{\mu}{\rho_0}}$ is the shear sound speed and $c_h$ is related to the speed of heat waves in the medium (also called the second sound \cite{Peshkov1944}, or the speed of a phonon). For an alternative choice of $E_2$ with better mathematical properties, see \cite{Ndanou2014}.  
Here $\mathring{\mathbf{G}} = \mathbf{G} - \frac{1}{3} \textnormal{tr}(\mathbf{G}) \mathbf{I}$ is 
the deviator of the metric (or Finger) tensor $\mathbf{G} = \mathbf{A}^T \mathbf{A}$ that describes the deformation of the medium. In tensor notation it reads as
\begin{equation} 
G_{ik} = A_{ji} A_{jk}, \qquad    \mathring{G}_{ik} = G_{ik} - \frac{1}{3} G_{jj} \, \delta_{ik}.   
\end{equation} 
Furthermore, the macro-scale energy is the classical kinetic energy and reads as 
\begin{equation} 
E_3 =  \frac{1}{2} v_i v_i. 
\end{equation} 
The determinant of $\mathbf{A}$ has to satisfy the algebraic constraint 
\begin{eqnarray}
|\mathbf{A}|=\frac{\rho}{\rho_0}.
\label{Aconstr}
\end{eqnarray}
Note that we define the energies $E_1, E_2$ by \eqref{E1a}, \eqref{E2a} in order to be able to reduce $E_1+E_2$ in case of small deformations to the quadratic energy of small strain corresponding to Hooke's law and indeed a more complex nonlinear dependence of $E_1, E_2$ on density and strain is allowed.

\subsection{Constitutive relations for the damaged medium}
\label{sec:damage}

In order to make the model \eqref{eqn.GPR} applicable to the description of the damage processes, it is necessary to define accordingly its material parameters and constitutive relations.
Our approach is based on the mixture model proposed in \cite{Resnyansky2003} for the damage of an elastoplastic continuum in case of small deformations and generalized later for the study of the structure of failure waves in the case of finite deformations in \cite{Romenskii2007}.  
The idea of the mentioned approach is to consider a damaged material as a mixture of the intact and "fully damaged" materials.  
These two materials have their own material parameters and closing relations, such as functions characterizing the rate of shear stress relaxation. The transition from the intact material to the fully damaged material is governed by the \textcolor{black}{damage variable} $\xi \in [0,1]$ satisfying the kinetic equation \eqref{eqn.xi} with source term depending on the state parameters of the medium (pressure, stress and temperature).
Then, if to assume that in the case of small deformations the mixture parameters of state satisfy the simple mixture rules such as an additivity of small strain and continuity of stress field, one can derive the governing equations for the damaged medium considered as an elasto-plastic continuum  with material parameters and closing relations depending on the \textcolor{black}{damage variable} $\xi$. This mixture model of damaged medium allows one to describe a degradation of elastic moduli during the damage process and to fit experimental stress-strain diagrams depending on the strain rate. 

Let us assume that the elastic moduli (Lam\'e constants) of the intact material $\lambda_I$, $\mu_I$ and of the fully damaged material $\lambda_D$, $\mu_D$ are known. 
Further, assume that the material parameters of both phases, corresponding to the heat transfer processes, such as heat capacity, thermal expansion coefficient and thermal conductivity coefficients are close to each other. The latter assumptions allow us to avoid an excessive complexity of the model.   

Thus, following \cite{Resnyansky2003}, one can define the elastic moduli in the equation of state of the damaged material as
\begin{equation}
\lambda(\xi)=\frac{K_IK_D}{\tilde K}-\frac{2}{3}\frac{\mu_I\mu_D}{\tilde \mu}, \quad \mu(\xi)=\frac{\mu_I\mu_D}{\tilde \mu},
\label{Lame}
\end{equation}
where $K_I=\lambda_I+\frac{2}{3}\mu_I$, $K_D=\lambda_D+\frac{2}{3}\mu_D$, $\tilde K=\xi K_I+(1-\xi)K_D$,
$\tilde \mu=\xi \mu_I+(1-\xi)\mu_D$.
It is easy to see that for the intact material ($\xi=0$) we have the mixture moduli equal to the ones of the intact material: $\lambda=\lambda_I$, $\mu=\mu_I$, while for the fully damaged material ($\xi=1$) we obtain $\lambda=\lambda_D$, $\mu=\mu_D$. 

The mixture-model averaging method of \cite{Resnyansky2003} applied to the definition of the rate of shear stress relaxation gives us the dependence
of the shear stress relaxation time $\tau_1$ on the \textcolor{black}{damage variable} $\xi$ as follows: 
\begin{equation}
\tau_1= \left(\frac{1-\xi}{\tau_I} + \frac{\xi}{\tau_D}\right)^{-1},
\end{equation}
where $\tau_I$ and $\tau_D$ are the shear stress relaxation times for the intact and fully damaged materials respectively, which are usually highly nonlinear functions of the parameters of state. The particular choice of $\tau_I$ and $\tau_D$ that is used in this  paper reads as
\begin{equation}
\tau_I=\tau_{I0}\exp(\alpha_I-\beta_I(1-\xi)Y), \quad \tau_D=\tau_{D0}\exp(\alpha_D-\beta_D\xi Y),
\end{equation}
where $Y$ is the equivalent stress (e.g. the Von Mises stress), while $\tau_{I0},\alpha_I,\beta_I$, $\tau_{D0},\alpha_D,\beta_D$ are constants.

The parameter $\theta$ governing the rate of damage is also a nonlinear function of the parameters of state and in our numerical examples we take it in the following form: 
\begin{equation}
\theta=\theta_0 (1-\xi)(\xi+\xi_\epsilon) \left[(1-\xi)\left(\frac{Y}{Y_0}\right)^a+\xi
\left(\frac{Y}{Y_1}\right)\right],
\end{equation} 
where $\xi_\epsilon$, $Y_0$ and $Y_1,a$ are constants. 
$\xi_\epsilon$ is usually taken as $10^{-16}$ in order to provide the growth of $\xi$ with the initial data  $\xi=0$.

The derivative $E_\xi$ can be computed with the use of \eqref{E1a},\eqref{E2a},\eqref{Lame} and reads as
\begin{equation}
\frac{\partial E}{\partial \xi}=-\frac{1}{2\rho_0}\frac{(K_I-K_D)K_IK_D}{\tilde K^2}\left(1-\frac{\rho}{\rho_0}\right)^2-
\frac{1}{4\rho_0}\frac{(\mu_I-\mu_D)\mu_I\mu_D}{\tilde \mu^2}\mathring{G}_{ij} \mathring{G}_{ij}.
\end{equation}

For the sake of simplicity we assume that the material parameters characterizing thermal properties of intact and damaged material are the same, that means that $c_v$, $c_h^2$ are constant and do not depend on $\xi$. 
The function $\theta_2$, characterizing the rate of heat flux relaxation, is taken as 
$\theta_2(\tau_2) = \tau_2 \frac{c_h^2}{\rho T}$ that yields the classical Fourier heat conduction law with the thermal conductivity coefficient $\kappa=\tau_2 c_h^2 $ in the stiff relaxation limit ($\tau_2 \rightarrow 0$) \cite{DumbserPeshkov2016}.

\subsection{Discussion}  
\label{sec:discussion}

The model of Section \ref{sec:pde} generalizes the unified model of continuum mechanics presented in \cite{PeshRom2014,DumbserPeshkov2016}, taking into account also material damage processes and the possibility to simulate moving free surface problems via the use of a diffuse interface method (DIM) that simply employs a scalar colour function $\alpha$ in order to define where the solid is present ($\alpha=1$) and where it is not ($\alpha=0$), see \cite{Tavelli2019}. The additional evolution equations for the material parameters \eqref{eqn1b} are added in order to capture correctly the motion of a continuum with heterogeneous material parameters undergoing large deformations.  
Note, that the governing equations for the heat transfer in the present model differ from that considered in \cite{DumbserPeshkov2016} and are taken in the form originally proposed in \cite{Rom1998}. The latter formulation of the heat transfer processes seems to be more natural, because it can be derived by the minimization of a Lagrangian and is in agreement with the Hamiltonian GENERIC formulation \cite{SHTC-GENERIC-CMAT}, as all equations from the general SHTC class. 
      
As is noted in Section \ref{sec:pde}, if $\alpha=1$ then the PDE system \eqref{eqn.GPR} is a hyperbolic thermodynamically compatible system and advanced high-order methods can be applied to solve these equations. Nevertheless, there are very stiff algebraic source terms in the equations for distortion, thermal impulse and \textcolor{black}{damage variable} which create a significant difficulty for numerical computations. The exponential ODE integrator presented below helps to avoid problems related to the stiffness of the algebraic source terms in the governing PDE system and a series of numerical test problems has been solved successfully with the use of high order ADER-DG schemes with \textit{a posteriori} subcell finite volume limiter method in conjunction with adaptive mesh refinement (AMR). 

Another difficulty is an application of the model of damaged medium to the solution of real problems. It relates to the definition of the constitutive relations and, in particular, in the definition of the function $\theta$ characterizing the rate of change of the \textcolor{black}{damage variable}, as well as the shear and heat flux relaxation parameters $\theta_1$, and $\theta_2$. 
On the one hand, the material constants in the equation of state of the intact elastic medium can in principle be found from experimental measurements, which means that the Lam\'e constants $\lambda_I$, $\mu_I$ and the heat capacity $c_v$ can be considered as known. On the other hand, there is no way to get the sample of the fully damaged medium which appeared as a result of the deformation of the medium, 
and hence there seems to be no possibility to obtain the material constants of the damaged material by direct measurements.   

The description of the damaged material requires not only material constants in the equation of state, but also a constitutive relations for the shear stress relaxation time $\tau_1$ and for the parameter $\theta$ governing the rate of the damage. All these characteristics of the medium can in principle be found with the use of experimental stress-strain diagrams, which are usually available from standard traction, torsion and compression experiments in solid  mechanics.  
The method consists in doing a series of numerical computations and obtaining a set of stress-strain diagrams numerically, and then by variation of the material constants try to fit as much as possible the available experimental diagrams.  
Such a procedure was successfully used in the past for the closure of the nonlinear  elastoplastic Godunov-Romenski model and the idea how to do this can be found in \cite{GodRom2003}. Recently this method has also been used for the closure of complex elastoplastic media with hardening, see \cite{Barton2012}.  

The dynamic behavior of damaged materials is very complex and depends on the type of the medium. The damage followed by fracture can be brittle, or brittle-ductile, which means that the stress-strain diagrams can be completely different for different materials. In Section \ref{sec.odevalidation} two typical examples of brittle and ductile material behavior are presented. In our model, these diagrams can also depend on the  \textit{strain rate}, which is also shown in the same section.
In this paper we do not calibrate the model parameters explicitly for such a dependence, since there is no experimental data available. Nevertheless, the constitutive relations chosen in Section \ref{sec:damage} and the set of material constants presented in the Section on numerical tests gives good results in all considered test problems. 

 Last but not least, our thermodynamically compatible approach to material failure naturally includes also the phenomenon of \textit{material fatigue}, i.e. the reduced resistance of the material to stress applied over a very large number of load cycles, see Section  \ref{sec.odevalidation}.

%==================================================================================================================
\section{Space-time adaptive ADER discontinuous Galerkin finite element schemes with a posteriori subcell finite volume limiter} 
\label{sec:ader}
The equations \eqref{eqn.mass}-\eqref{eqn.energy} of the GPR model described above can be written in the following 
general form of a nonlinear system of hyperbolic PDEs with non-conservative products and stiff source terms: 
\begin{equation}    
\label{eqn.pde.nc}
    \frac{\partial \Q}{\partial t} + \nabla \cdot \bf F(\Q) + \mathbf{\mathcal{B}}(\Q) \cdot \nabla \Q = \S(\Q), 
\end{equation}
where $\Q=\Q(\x,t)$ is the state vector; $\x=(x,y,z) \in \Omega \subset \mathds{R}^d$ is the vector of spatial coordinates and 
$\Omega$ denotes the computational domain in $d$ space dimensions; ${\bf F}(\Q) = (\f, \g, \h)$ is the nonlinear flux tensor that contains 
the conservative part of the PDE system and $\mathbf{\mathcal{B}}(\Q) \cdot \nabla \Q$ is a genuinely non-conservative 
term. When written in quasilinear form, the system (\ref{eqn.pde.nc}) becomes 
\begin{equation}    \label{eq:Csyst}
\frac{\partial \Q}{\partial t}+ \mathbf{\mathcal{A}} (\Q) \cdot\nabla \Q = \S(\Q)\,,
\end{equation}
where the matrix $\mathbf{\mathcal{A}}(\Q)=\partial {\bf F}(\Q)/\partial \Q + \mathbf{\mathcal{B}}(\Q)$ includes both 
the Jacobian of the conservative flux, as well as the non-conservative product. The hyperbolicity of system \eqref{eq:Csyst} 
has been discussed in \cite{PeshRom2014}. However, for the practical implementation of the numerical schemes used in this paper, 
the eigenvectors $\mathbf{R}_n$ of the matrix $\mathbf{\mathcal{A}}_n = \mathbf{\mathcal{A}}(\Q) \cdot \mathbf{n}$ 
($\mathbf{n}$ is a unit-normal vector) will not be needed, even if they were in principle available. 

The PDE system \eqref{eqn.pde.nc} is solved by resorting to a high order one-step ADER-FV and ADER-DG method 
\cite{QiuDumbserShu,Dumbser2008,ADERNC}, which provides at the same time high order of accuracy in 
both space and time in one single step, hence completely avoiding the Runge-Kutta sub-stages that are typically
used in Runge-Kutta DG and Runge-Kutta WENO schemes. The method will be presented in the \textit{unified} framework
of $P_NP_M$ methods introduced in \cite{Dumbser2008}, which contains both, DG schemes and FV schemes as special cases
of a more general class of methods. For related work on $P_NP_M$ schemes, the reader is referred to 
\cite{luo1,luo2}. The construction of fully-discrete high order one-step schemes is typical of the ADER approach \cite{titarevtoro,toro3,toro4}.  
In the following, we only summarize the main steps, while for more details the reader is referred to 
\cite{Dumbser2008,DumbserZanotti,HidalgoDumbser,GassnerDumbserMunz,Balsara2013,Dumbser2014,Zanotti2015a,Zanotti2015b,FrontierADERGPR}.

%---------------------------------------------------------------------
\subsection{Data representation and reconstruction}

The computational domain $\Omega$ is discretized by a computational mesh that can be  structured or unstructured, composed of conforming elements denoted by $T_i$, where the index $i$ ranges from 1 to the total number of elements $N_E$. 
We will further denote the volume (area) of an individual cell by 
$|T_i| = \int_{T_i} d\x$. The discrete solution of PDE \eqref{eqn.pde.nc} is denoted by 
$\u_h(\x,t^n)$ and is represented by piecewise polynomials of maximum degree 
$N \geq 0$. Within each cell $T_i$ we have 
\begin{equation}
\label{eqn.ansatz.uh}
  \u_h(\x,t^n) = \sum_l^{\mathcal{N}} \Phi_l(\x) \hat{\u}^n_{l,i} := \Phi_l(\x) \, \hat{\u}^n_{l,i},  \quad \x \in T_i, 
\end{equation}
where we have introduced the classical Einstein summation convention over two repeated indices. 
The discrete solution $\u_h(\x,t^n)$ is defined in the space of piecewise polynomials up 
to degree $N$, spanned by a set of basis functions $\Phi_l=\Phi_l(\x)$. Throughout this paper we use a tensor-product \textit{nodal basis} for quadrilateral and hexahedral elements. The  nodal basis is given by the Lagrange interpolation polynomials passing through the Gauss-Legendre quadrature nodes on the unit element $[0,1]^d$, see \cite{stroud} for details on multidimensional quadrature. The symbol $\mathcal{N}$ denotes the number of degrees of freedom per element and is given by $\mathcal{N}=(N+1)^d$ for  quadrilateral / hexahedral Cartesian elements in $d$ space dimensions. In the framework of $P_NP_M$ methods, the 
discrete solution $\u_h$ is now \textit{reconstructed} in order to obtain for each element a piecewise polynomial $\w_h(\x,t)$ 
of degree $M\geq N$, with a total number of $\mathcal{M}$ degrees of freedom. Details on the nonlinear WENO reconstruction and on the 
$P_NP_M$ reconstruction can be found in \cite{DumbserKaeser06b,DumbserKaeser07,Dumbser2008} and are not repeated here. The number of degrees 
of freedom $\mathcal{M}$ is again $\mathcal{M}=(M+1)^d$ for tensor-product elements in $d$ space dimensions, respectively. The reconstruction step is simply abbreviated by $\w_h(\x,t) = \mathcal{R}( \u_h(\x,t) )$, and the reconstruction 
polynomial $\w_h(\x,t)$ is written as 
\begin{equation}
  \w_h(\x,t^n) = \sum_l^{\mathcal{M}} \Psi_l(\x) \hat{\w}^n_{l,i} := \Psi_l(\x) \, \hat{\w}^n_{l,i},  \quad \x \in T_i.  
 \label{eqn.recsol} 
\end{equation} 
Note that for $N=M$ the $P_NP_M$ method reduces to a classical discontinuous Galerkin finite element scheme, with the reconstruction operator
equal to the identity operator, $\mathcal{R}=\mathcal{I}$, or, equivalently, $\w_h(\x,t^n) = \u_h(\x,t^n)$, while for the case $N=0$ the method reduces to a standard high order WENO finite 
volume scheme, if a WENO reconstruction operator is adopted. 

For WENO schemes on structured meshes we have found that it is particularly convenient to adopt one-dimensional stencils, each composed by $n_e=M+1$ cells, 
which are subsequently oriented along each spatial direction. The resulting reconstruction is still multidimensional, but implemented with a dimension-by-dimension 
strategy. A complete description of this approach can be found in \cite{AMR3DCL,Zanotti2015,FrontierADERGPR}, including also the necessary details for the employed adaptive mesh refinement.  

In this paper, however, we will only use the two special limits of the general $P_NP_M$ approach, i.e. either $N=0$ (pure FV) or $N=M$ (pure DG), where the finite volume scheme is used under the form of an \textit{a posteriori} subcell finite volume limiter, see \cite{Dumbser2014,Zanotti2015a}. 

\subsection{Local space-time predictor}
\label{sec.predictor}

The discrete solution $\w_h(\x,t^n)$ is now evolved in time according to an element-local weak formulation of the governing PDE in space-time, 
see \cite{DumbserEnauxToro,Dumbser2008,HidalgoDumbser,DumbserZanotti,GassnerDumbserMunz,Balsara2013,Dumbser2014,Zanotti2015a,Zanotti2015b,FrontierADERGPR}. 
The local space-time Galerkin method is only used for the construction of an element-local predictor solution of the PDE 
\textit{in the small}, hence neglecting the influence of neighbor elements. This predictor will subsequently be inserted into the corrector 
step described in the next section, which then provides the appropriate coupling between neighbor elements via a numerical flux function 
(Riemann solver) and a path-conservative jump term for the discretization of the non-conservative product. To simplify notation, we define  
\begin{equation}
 \label{eqn.operators1}
  \left<f,g\right> =
      \int \limits_{t^n}^{t^{n+1}} \int \limits_{T_i}  f(\x, t)  g(\x, t)  \, d \x \, d t,
\qquad 
  \left[f,g\right]^{t} =
      \int \limits_{T_i}f(\x, t) g(\x, t) \,  d \x,
\end{equation}
which denote the scalar products of two functions $f$ and $g$ over the space-time element $T_i \times \left[t^n;t^{n+1}\right]$ and  
over the spatial element $T_i$ at time $t$, respectively. Within the local space-time predictor, the discrete solution of 
equation \eqref{eqn.pde.nc} is denoted by $\q_h=\q_h(\x,t)$. 
We then multiply \eqref{eqn.pde.nc} with a space-time test function $\theta_k=\theta_k(\x,t)$ and subsequently integrate over 
the space-time control volume $T_i \times \left[t^n;t^{n+1}\right]$. Inserting $\q_h$, the following weak formulation of 
the PDE is obtained: 
\begin{equation}
\label{eqn.pde.nc.weak1}
 \left< \theta_k, \frac{\partial \q_h}{\partial t}  \right>
    + \left< \theta_k, \nabla \cdot \F \left(\q_h\right) + \mathbf{\mathcal{B}}(\q_h) \cdot \nabla \q_h \right> = \left< \theta_k, \S \left( \q_h \right)  \right>.
\end{equation}
The discrete representation of $\q_h$ in element $T_i \times [t^n,t^{n+1}]$ is assumed to have the following form 
\begin{equation}
\label{eqn.st.state}
 \q_h = \q_h(\x,t) =
 \sum \limits_l \theta_l(\x,t) \hat{\q}^n_{l,i} := \theta_l \hat{\q}^n_{l,i},
\end{equation}
where $\theta_l(\x,t)$ is a space-time basis function of maximum degree $M$. 
For the basis functions $\theta_l$ we use a tensor-product of 1D nodal basis functions given by the Lagrange interpolation polynomials of the 
Gauss-Legendre quadrature points for Cartesian tensor-product elements. 
After integration by parts in time of the first term, Eqn. \eqref{eqn.pde.nc.weak1} reads 
\begin{equation}
\label{eqn.pde.nc.dg1}
 \left[ \theta_k, \q_h \right]^{t^{n+1}} - \left[ \theta_k, \w_h(\x,t^n) \right]^{t^n} - \left< \frac{\partial}{\partial t} \theta_k, \q_h \right> 
    + \left< \theta_k, \nabla \cdot \F \left(\q_h\right) + \mathbf{\mathcal{B}}(\q_h) \cdot \nabla \q_h \right> = \left< \theta_k, \S \left( \q_h \right)  \right>. 
\end{equation}
Note that the high order polynomial reconstruction of the $P_NP_M$ scheme $\w_h(\x,t^n)$ is taken into account 
in \eqref{eqn.pde.nc.dg1} in a \textit{weak sense} by the term $\left[ \theta_k, \w_h(\x,t^n) \right]^{t^n}$. 
This corresponds to the choice of a numerical flux in time direction, which is nothing else than \textit{upwinding in time}, 
according to the causality principle. 

Note further that due to the DG approximation in space-time, we may have $\q_h(\x,t^n) \neq \w_h(\x,t^n)$ in general, hence
the choice of a numerical flux in time direction is necessary. Note further that in \eqref{eqn.pde.nc.dg1} we have
\textit{not} used integration by parts in space, nor any other coupling to spatial neighbor elements. The integrals
appearing in the weak form \eqref{eqn.pde.nc.dg1}, as well as the space-time test and basis functions involved are 
conveniently written by making use of a space-time reference element $T_e \times [0;1]$. \\ 
The solution of \eqref{eqn.pde.nc.dg1} yields the unknown space-time degrees of freedom $\hat{\q}^n_{l,i}$ for each 
space-time element $T_i \times [t^n; t^{n+1}]$ and is easily achieved with a fast converging iterative scheme, see 
\cite{Dumbser2008,HidalgoDumbser,DumbserZanotti} for more details. In \cite{Jackson} it was proven that the resulting iteration matrix is nilpotent, which explains the convergence to the exact solution in a finite number of iterations for linear homogeneous systems already observed in \cite{Dumbser2008}. In \cite{FrontierADERGPR} a more general convergence proof of the space-time predictor based on fixed point arguments was given for nonlinear conservation laws. 
The above space-time Galerkin predictor has replaced the cumbersome Cauchy-Kovalewski procedure that has been initially 
employed in the original version of ADER finite volume and ADER discontinuous Galerkin schemes 
\cite{schwartzkopff,toro3,toro4,titarevtoro,dumbser_jsc,taube_jsc,DumbserKaeser07}. 
Note that very recently, a new reformulation of the ADER method has been proposed, where the reconstruction and time evolution steps  
are performed in terms of the vector of \textit{primitive variables} $\QV$ instead of using the vector $\QQ$ of
conserved quantities, see \cite{ADERPrim} and also \cite{PuppoRussoPrim}, for a similar approach in the context of Runge-Kutta WENO FV schemes.

\subsection{Fully discrete one-step finite volume and discontinuous Galerkin schemes}
\label{sec.ADERNC}

At the aid of the local space-time predictor $\q_h$, a fully discrete one-step $P_NP_M$ scheme can now be 
simply obtained by multiplication of the governing PDE system \eqref{eqn.pde.nc} by test functions  
$\Phi_k$ from the space of piecewise polynomials up to degree $N$, which are identical with the spatial basis functions of the 
original data representation 
before reconstruction, and subsequent integration over the space-time control volume $T_i \times [t^n;t^{n+1}]$. Due 
to the presence of non-conservative products, the jumps of $\q_h$ across element boundaries are taken into account in 
the framework of path-conservative schemes put forward by Castro and Par\'es in the finite volume context 
\cite{Castro2006,Pares2006} and subsequently extended to DG schemes in \cite{Rhebergen2008} and \cite{ADERNC,USFORCE2}, 
where also a generalization to the unified $P_NP_M$ framework has been provided. All these approaches are based on the 
theory of Dal Maso, Le Floch and Murat \cite{DLMtheory}, which gives a definition of weak solutions in the context of 
non-conservative hyperbolic PDE. 
For open problems concerning path-conservative schemes, the reader is referred to \cite{AbgrallKarni,NCproblems}. 

If $\mathbf{n}$ is the outward pointing unit normal vector on the surface $\partial T_i$ of element $T_i$ and  
the path-conservative jump term in normal direction is denoted by $\mathcal{D}_h^-\left(\q_h^-, \q_h^+ \right) \cdot\mathbf{n}$, 
which is a function of the left and right boundary-extrapolated data, 
$\q_h^-$ and $\q_h^+$, respectively, then we obtain the following path-conservative one-step 
$P_NP_M$ scheme, see \cite{ADERNC}: 
\begin{equation}
\label{eqn.pde.nc.gw2}
\begin{split}
\left( \int \limits_{T_i} \Phi_k \Phi_l d\x \right) \left( \hat{\u}_l^{n+1} -  \hat{\u}_l^{n} \right) +
\int \limits_{t^n}^{t^{n+1}} \int \limits_{\partial T_i} \Phi_k \, \mathcal{D}_h^-\left(\q_h^-, \q_h^+ \right)\cdot\mathbf{n} \, dS\, dt 
\\ 
+\int\limits_{t^n}^{t^{n+1}} \int \limits_{T_i \backslash \partial T_i} \Phi_k \left( \nabla \cdot \F\left(\q_h \right) + 
\mathbf{\mathcal{B}}(\q_h) \cdot \nabla \q_h \right) d\x\, dt  
= \int \limits_{t^n}^{t^{n+1}} \int \limits_{T_i} \Phi_k \S(\q_h) d\x\, dt. 
\end{split} 
\end{equation}
The element mass matrix appears in the first integral of \eqref{eqn.pde.nc.gw2}, the second term accounts for the jump in the discrete 
solution at element boundaries and the third term takes into account the smooth part of the non-conservative product. 
For general complex nonlinear hyperbolic PDE systems we use the simple Rusanov method \cite{Rusanov:1961a} (also called the local Lax Friedrichs
method), although any other kind of Riemann solver could be also used, see \cite{toro-book} for an overview of state-of-the-art Riemann solvers. 
In these regards, we would also like to point out the new general reformulation of the HLLEM Riemann solver of Einfeldt and Munz 
\cite{Einfeldt88,munz91}, within the setting of path-conservative schemes recently forwarded in \cite{NewHLLEM}, as well as the family of 
MUSTA schemes, which has been applied to the equations of nonlinear elasticity in \cite{TitarevRomenskiToro}.
\\

The extension of the Rusanov flux to the path-conservative framework requires the introduction of an additional nonconservative jump term and reads 
\begin{equation}
  \mathcal{D}_h^-\left(\q_h^-, \q_h^+ \right)\cdot\mathbf{n} = \halb \left( \F(\q_h^+) - \F(\q_h^-) \right) \cdot \mathbf{n} + 
    \halb \left( \tilde{\mathbf{\mathcal{B}}} \cdot \mathbf{n} - s_{\max} \mathbf{I} \right) \left( \q_h^+ - \q_h^- \right), 
    \label{eqn.rusanov} 
\end{equation} 
with the maximum signal speed at the element interface $s_{\max} = \max\left( \left|\boldsymbol{\Lambda}(\q_h^+) \right|, \left|\boldsymbol{\Lambda}(\q_h^-) \right| \right)$ and
the matrix $\tilde{\mathbf{\mathcal{B}}} \cdot \mathbf{n}$ given by the following path-integral along a straight line segment path $\psi$: 
\begin{equation}
 \tilde{\mathbf{\mathcal{B}}} \cdot \mathbf{n} = \int \limits_0^1 \mathbf{\mathcal{B}}\left( \psi(\q_h^-, \q_h^+, s \right) \cdot \mathbf{n} \, ds, 
\qquad 
\psi \left( \q_h^-, \q_h^+, s \right) = \q_h^- + s \left( \q_h^+ - \q_h^- \right).  
\end{equation} 
According to the suggestions made in \cite{ADERNC,USFORCE2,OsherUniversal,ApproxOsher,OsherNC}, the path-integrals can be conveniently 
evaluated numerically by the use of a classical Gauss-Legendre quadrature formula on the unit interval $[0;1]$. For an alternative choice 
of the path, see \cite{MuellerToro1,MuellerToro2}. 

This completes the brief description of the $P_NP_M$ scheme used for the discretization of the governing PDE system 
\eqref{eqn.pde.nc}. 

In the case of ADER finite volume schemes, we simply have $N=0$, $\mathcal{N}=1$, $\Phi_k=1$, and 
the limiter is directly incorporated in the \textit{nonlinear} reconstruction operator $\w_h(\x,t^n) = \mathcal{R}\left(\u_h(\x,t^n) \right)$.  
For ADER discontinuous Galerkin finite element schemes ($N=M$, $\Phi_k = \Psi_k$) a new family of \aposteriori  sub-cell finite volume limiters has been forwarded in \cite{Dumbser2014,Zanotti2015a,Zanotti2015b} and is employed throughout this work. In particular, 
in this paper the subcell finite volume limiter employs simple piecewise linear reconstruction 
based on the minmod slope limiter. 
For alternative finite volume subcell limiters in the context of DG schemes, see the work of Sonntag \& Munz \cite{Sonntag} and Meister \& Ortleb \cite{MeisterOrtleb}. 
For details on the adaptive mesh  refinement strategy used in this paper, the reader is referred to     \cite{AMR3DCL,Zanotti2015a,Zanotti2015b,Peano1,Peano2,Reinarz2020}.

\subsection{Treatment of stiff source terms via exponential time integration}
Due to the strong nonlinearities present in the reaction-like kinetics characterising
the evolution of the material \textcolor{black}{damage variable} $\xi$, a special treatment of the stiff 
source term is required. Note that the exponential time integration method described in this 
section is applied only within the subcell finite volume limiter scheme. 

Specifically, the rupture process develops at time scales that are much
quicker than what given by the CFL condition for elastic waves in the intact medium and thus
the effects of the stiff source terms might compromise the stability of the simulations. 

In order to solve the issue, we make use of a split treatment of the source, that is, 
at each timestep we solve the governing equations neglecting the source terms for strain relaxation and rupture dynamics, 
obtaining a preliminary solution, to be used in a second step as initial condition for a Cauchy problem
constructed by neglecting all the differential terms appearing in the PDE system, except for time derivatives.
Then one can evolve in time the pointwise preliminary solution by solving such initial value problem 
between the current time level and the following one. 

Unfortunately, the integration of the ordinary differential equations
arising from the split treatment of the algebraic source term 
cannot be tackled by under-resolving the sub-timestep evolution of the \textcolor{black}{damage variable} $\xi$
and the distortion matrix $\vec{A}$ with a simple implicit Euler scheme. 
Instead we decided to employ an efficient exponential integration technique 
that is not only more robust, but also more accurate.

In the following we present the procedure for applying the exponential integrator 
to a general first order system of ordinary differential equations like

% \begin{equation} \label{eq:referencesystem}
%     \od{\vec{Q}}{t} = \vec{A}(\vec{Q},\ t)\,\vec{Q} + \vec{B}(\vec{Q},\ t).
% \end{equation}
% Note that any nonlinear system of the type 
\begin{equation} \label{eq:referencesystem}
\od{\vec{Q}}{t} = \vec{S}(\vec{Q},\ t), 
\end{equation}
for which we write a linearisation about a given state $\Qs$ and time $\ts$ as
\begin{equation}
\od{\vec{Q}}{t} = \Bs + \Js(\Qs,\ \ts)\,(\vec{Q} - \Qs).
\end{equation}
In this work, the vector $\vec{Q}$ will have ten components, the first being the material damage
variable $\xi$, followed by the nine entries of the distortion matrix $\vec{A}$.
% can be cast as \eqref{eq:referencesystem} in multiple ways, for example taking $\vec{A}$ to be the Jacobian matrix
% of the nonlinear source vector $\vec{S}(\vec{Q},\ t)$ or by simply setting $\vec{B} = \vec{S}$ and neglecting $\vec{A}$.
% The specific choice of in this regard will be crucial to the success and performance of the resulting scheme. 
% In the case of relaxation systems, the equation are cast quite easily in the form \eqref{eq:referencesystem} with
% $\vec{B}(\vec{Q},\ t) = \vec{0}$,
% but it may be convenient to introduce non-homogeneous terms $\vec{B}(\vec{Q},\ t)$ in order to simplify some steps of the algorithm, 
% depending on the specific ODE at hand.
One then defines the Jacobian matrix of the source ${\Js = \vec{J}(\Qs,\ \ts)} = \partial \mathbf{S} / \partial \mathbf{Q}$ and in the same way the source  
vector evaluated at the linearisation state is ${\Bs = \vec{S}(\Qs,\ \ts)}$. 
% Moreover, for notational convenience we will define the affine transformation operator $\Sss$ such that
% \begin{equation}
% \Ss(\vec{Q};\ \Qs,\ \ts) = \Js\,\vec{Q} + \Ss = \vec{A}(\Qs,\ \ts)\,\vec{Q} + \vec{B}(\Qs,\ \ts),
% \end{equation}
Furthermore, we introduce the vector 
\begin{equation}
\Cs = \Cs(\Bs,\ \Js) = \Cs(\Qs,\ \ts)
\end{equation}
which will be used as an indicator for the adaptive timestepping scheme.
This vector may be
composed for example by the the entries of the matrix $\Js$, 
together with all the elements of the vector $\Bs$ and the state $\Qs$. A simpler general choice 
might be using only the state vector, but we recommend that some relevant 
combination of the listed variables be included, to indicate changes in the nature or in the magnitude
of the linearised source operator. For the tests presented in
this work the vector of indicator variables is comprised simply of the linearisation state $\Qs$ and of two characteristic 
relaxation times for the \textcolor{black}{damage variable} equation and for the strain relaxation equation.
% Selection of the content of $\Cs$ is carried out 
% One may further extend the contents of $\Cs$, for example it could happen that 
% some of the terms of $\Js$ are computed as a binary operator applied of two distinct 
% functions or another type of composition of more terms, and it might be useful to 
% add the individual factors to the vector of constant coefficients $\Cs$. In the following, 
% whenever we compute $\Cs(\Qs,\ \ts)$, it is intended that this is equivalent to computing both
% $\Js$ and $\Ss$ and thus constructing the affine operator $\Ss$.
% The source Jacobian matrix may be easily approximated with central finite differences, 
% or alternatively, one can also provide a simplified Jacobian matrix, constructed for example by decoupling 
% different sub-systems of equations or allowing some terms to be arbitrarily taken out of the derivative
% operator as constants. 
% 
It is then necessary to compute an exact analytical solution of the linear nonhomogeneous 
Cauchy problem 
\begin{equation} \label{eq:linearcauchy}
\left\{
\begin{aligned}
&\od{\vec{Q}}{t} = \Ss(\vec{Q};\ \Qs,\ \ts) = \Bs + \Js(\Qs,\ \ts)\,(\vec{Q} - \Qs),\\[2mm]
& \vec{Q}(t_n) = \vec{Q}_n,\\
\end{aligned}
\right.
\end{equation}
and in the general case, one can use the algorithms of Higham \cite{matexp2005} and Al-Mohy and Higham \cite{matexp2009, matexp2011} 
for a robust evaluation of matrix exponentials needed for such computations. More notes on the solution of \eqref{eq:linearcauchy} 
will be given in Section~\ref{sec:linearisedsol}.
%  and, if one were to choose to employ an approximate
% Jacobian with a block structure, it is natural to compute the exponentials of the smaller blocks, further 
% reducing the cost of the solver. 

We will set $\Qe(t;\ \Ss,\ t_n,\ \vec{Q}_n)$ to denote the analytical solution of \eqref{eq:linearcauchy}, and, 
as for $\Ss(\vec{Q};\ \Qs,\ \ts)$, the semicolon will separate the variable on which $\Qe$ and $\Ss$ 
depend continuously ($t$ or $\vec{Q}$) from the fixed parameters used to build the operators. 
The state vector at timestep $t_n$ is $\vec{Q}_n$, and the variable timestep size
is $\Delta t^n = t_{n+1} - t_n$.

\subsubsection{Timestepping}

Marching from a start time $t_0$ to an end time $t^\up{end}$ is carried out as follows.
First, an initial timestep size $\Delta t^0$ is chosen; 
% for this choice one may use in principle any value, 
% but a good estimate for the first step can be obtained by dividing a characteristic timescale arising
% from the analytical solution $\Qe(t;\ \Ss,\ t_0,\ \vec{Q}_0)$ obtained from an affine operator 
% $\Ss(\vec{Q};\ \vec{Q}_0,\ t_0)$ 
% based on the initial condition.
then, at each timestep, the state $\vec{Q}_{n+1}$ at the new time level $t_{n+1}$ 
can be computed by means of an iterative 
procedure which will terminate by yielding a value for
$\vec{Q}_{n+1}$, as well as a new timestep size ${\Delta t^{n+1} = t_{n+2} - t_{n+1}}$ based on an estimator
embedded in the iterative solution algorithm. There is also the possibility that, 
due to the timestep size $\Delta t$ being too large, the value of $\vec{Q}_{n+1}$ be flagged as not acceptable. 
In this case the procedure returns a reduced timestep size for the current timestep 
${\Delta t^{n} = t_{n+1} - t_{n}}$ and the solution for $\vec{Q}_{n+1}$ will be attempted again using this 
reduced timestep size. Specifically, in practice we set the new value of $\Delta t^{n}$
to be half of the one used in the previous attempt.

\subsection{Iterative computation of the solution during one timestep}
At each iteration of index $k$ an average state vector is computed as
\begin{equation}
\Qs_{n+1/2}^k = \frac{1}{2}\,\left(\vec{Q}_n + \Qs_{n+1}^{k-1}\right)
\end{equation}
and this midpoint state is formally associated with an intermediate time level
\begin{equation}
t_{n+1/2} = \frac{1}{2}\,\left(t_n + t_{n+1}\right).
\end{equation}
For the first iteration we provide a guess value for $\Qs_{n+1}^{k-1}$, with the simplest choice being given by
$\Qs_{n+1}^{k-1} = \vec{Q}_n$. Then the coefficients $\Cs_{n+1/2}^k$ are computed as
\begin{equation}
\Cs_{n+1/2}^k = \Cs_{n+1/2}^k(\Qs_{n+1/2}^k,\ t_{n+1/2}),
\end{equation}
and concurrently (that is, within the same code block), one can build the affine source operator
\begin{equation}
\Ss_{n+1/2}^k = \Ss_{n+1/2}^k(\vec{Q};\ \Qs_{n+1/2}^k,\ t_{n+1/2}),
\end{equation}
then one has to find the solution of the initial value problem for the linearised equations
\begin{equation} \label{eq:approxcauchy}
\left\{
\begin{aligned}
&\od{\vec{Q}}{t} = \Ss_{n+1/2}^k(\vec{Q};\ \Qs_{n+1/2}^k,\ t_{n+1/2}),\\[2mm]
& \vec{Q}(t_n) = \vec{Q}_n,\\
\end{aligned}
\right. 
\end{equation} 
yielding the updated state vector for the next iteration
\begin{equation}
\Qs_{n + 1}^k = \vec{Q}_\up{e}\left(t_{n+1};\ \Ss_{n+1/2}^k,\ t_n,\ \vec{Q}_n\right).
\end{equation}
It is then checked that the state vector $\Qs_{n+1}^k$ be physically admissible: in our case this means verifying that the \textcolor{black}{damage variable} is in the unit interval, i.e. $ \xi \in [0,1]$. 
% but for example
% one may require the positivity of mass density and internal energy, and any boundedness requirement 
% on the state can be implemented, 
Also, one can check that floating-point exceptions be absent.
Additionally, one must evaluate
\begin{equation}
\Cs_{n+1}^k = \Cs_{n+1}^k\left(\Qs_{n+1}^k,\ t_{n+1}\right),
\end{equation}
this vector of coefficients will not be employed for the construction of an affine source operator $\Ss_{n+1}^k$, 
but only for checking the validity of the solution obtained from \eqref{eq:approxcauchy}, 
by comparing the coefficients vector $\Cs_{n+1}^k$ to $\Cs_{n}$, as well as comparing the coefficients $\Cs_{n+1/2}^k$
used in the middle-point affine operator for the initial coefficients $\Cs_{n}$. At the end of the iterative procedure, 
one will set $\Cs_{n+1} = \Cs_{n+1}^k$, so that this can be reused as the new reference vector of coefficients 
for the next timestep.

The convergence criterion for stopping the iterative computation
is implemented by evaluating 
\begin{equation} \label{eq:convergencemetric}
r = \max\left(\frac{\abs{\Qs_{n+1}^{k} - \Qs_{n+1}^{k-1}}}{\abs{\Qs_{n+1}^k} + \abs{\Qs_{n+1}^{k-1}} + \epsilon_r}\right), 
\end{equation}
and checking whether $r \leq r_\up{max}$, with $r_\up{max}$ and $\epsilon_r$ given tolerances, 
or alternatively when the iteration count $k$ has reached a fixed maximum value $k_\up{max}$.
% In principle any norm may be used to compute the indicator \eqref{eq:convergencemetric}, as this is
% just a measure of the degree to which $\Qs_{n+1}^k$ was corrected in the current iteration.
Moreover, we
found convenient to limit the maximum number of iterations allowed, and specifically we set $k_\up{max} = 8$, or
alternatively one can decide to use $k_\up{max} = 3$ and flag the state vector $\Qs_{n+1}^{k}$ as not admissible, 
as if a floating-point exception had been triggered,
whenever the iterative procedure terminates by reaching the maximum iteration count.
When convergence is finally obtained, in order to test if the Cauchy problem
\begin{equation} \label{eq:exactcauchy}
\left\{
\begin{aligned}
&\od{\vec{Q}}{t} = \vec{S}(\vec{Q},\ t),\\[2mm]
& \vec{Q}(t_n) = \vec{Q}_n.\\
\end{aligned}
\right.
\end{equation}
is well approximated by its linearised version \eqref{eq:approxcauchy}, we compute two error metrics
\begin{align}
&\delta_{n+1/2} = \max\left(\frac{\abs{\Cs_{n+1/2} - \Cs_{n}}}{\abs{\Cs_{n+1/2}} + \abs{\Cs_{n}} + \epsilon_\delta}\right), \label{eq:deltaa}\\
&\delta_{n+1} = \max\left(\frac{\abs{\Cs_{n+1} - \Cs_{n}}}{\abs{\Cs_{n+1}} + \abs{\Cs_{n}} + \epsilon_\delta}\right), \label{eq:deltab}
\end{align}
and we verify if 
\begin{equation} \label{eq:testtol}
\delta = \max(\delta_{n+1/2},\ \delta_{n+1}) \leq \delta_\up{max}.
\end{equation}
Note that one should specify a tolerance $\delta_\up{max}$ for \eqref{eq:testtol} and a floor value $\epsilon_\delta$, 
which is used in order to prevent that excessive precision requirements be imposed in those situations
when all the coefficients are so small than even large relative variations expressed by Equations \eqref{eq:deltaa} and \eqref{eq:deltab}
do not affect the solution in a significant manner.
If $\delta \leq \delta_\up{max}$, then the state vector at the new time level is set to be
$\vec{Q}_{n+1} = \Qs_{n+1}^k$ and a new timestep size can be computed as
\begin{equation}
\Delta t_{n+1} = \lambda\,\frac{\delta_\up{max}}{\delta + \epsilon},\quad \text{with}\quad \lambda = 0.8,\quad\epsilon = 10^{-{14}},
\end{equation}
otherwise a new attempt at the solution of \eqref{eq:exactcauchy} is carried out, 
with a halved timestep size. The same will happens in case
at any time the admissibility test on $\Qs_{n+1}^k$ fails.

The resulting ODE integrator is second order accurate, that is, the error in the solution
decreases quadratically with the number of timesteps employed, and, thanks to the piecewise-exponential ansatz, yields by construction very good results compared to a standard second order pieciewise linear evolution of the state vector, especially in the integration steps in which the equations are only weakly nonlinear, i.e.   
away from rupture faults and cracks.

\subsubsection{Solution of the linearised problem}\label{sec:linearisedsol}

If the system of linearised equations is large, 
it might be infeasible to explicitly write a closed form expression for
the solution of the linearised initial value problem \eqref{eq:linearcauchy}, 
which can be written compactly as 
\begin{equation}
\vec{Q}(t) = \exp\left(\Js\,(t - t_n)\right)\,\left(\vec{Q}(t_n) + 
    \left(\Js\right)^{-1}\,\Bs - \Qs\right) - \left(\Js\right)^{-1}\,\Bs + \Qs, 
\end{equation}
but the evaluation of the matrix exponential and the computation of the inverse Jacobian $\left(\Js\right)^{-1}$
in general need to be carried out numerically, and in particular the inversion of $\Js$ may constitute an ill-conditioned problem. Specifically,  
the timescale associated with strain relaxation can be significantly different from the reaction speed of
material failure. This issue can be treated by detecting those rows of the Jacobian that 
have much smaller entries than a given global scale for the full system and removing them from the matrix, 
then inverting a reduced system of equations with better conditioning.

\begin{figure}[!b]
    \centering
        \includegraphics[width=1.0\textwidth]{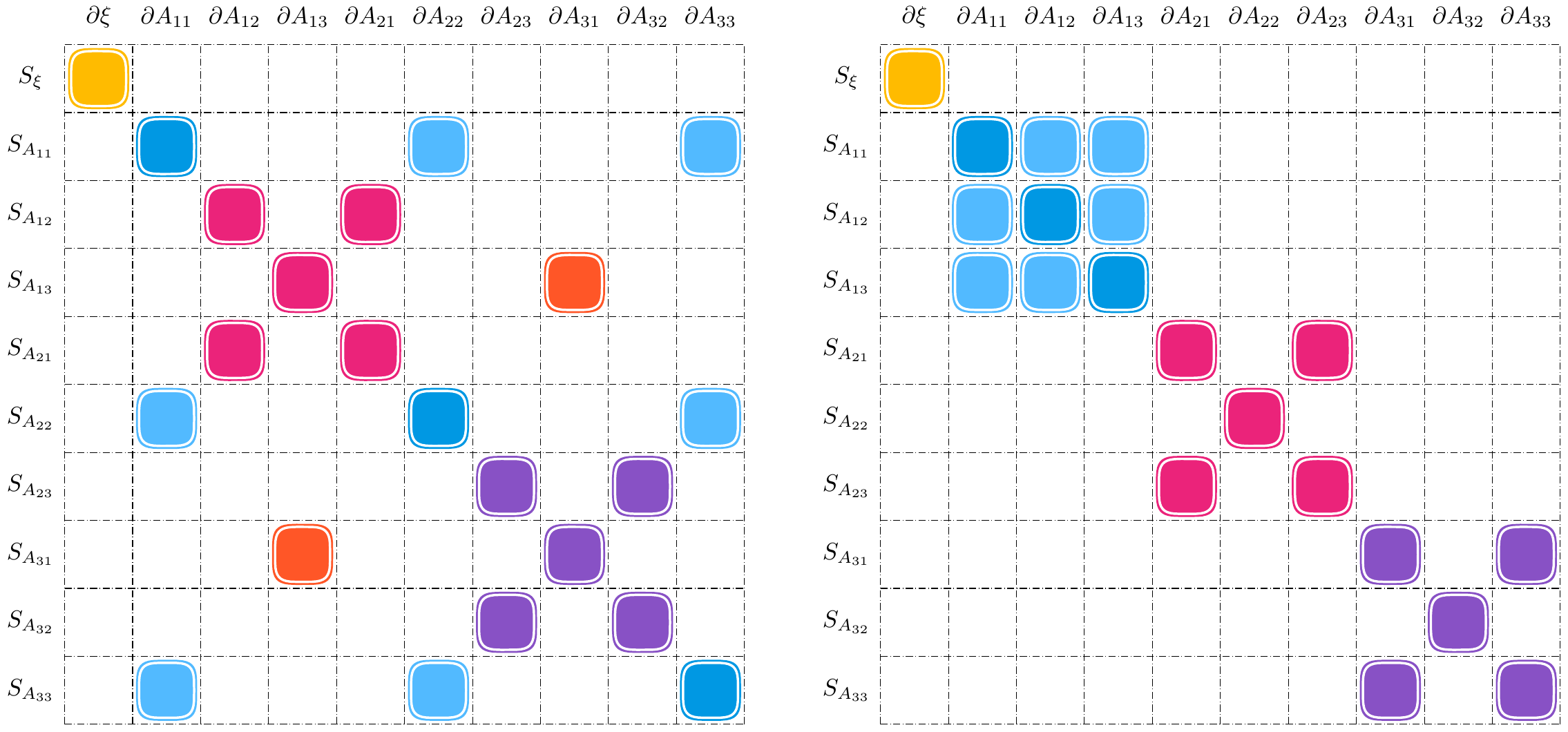}
    \caption{Structure of the two-block split (material failure and strain relaxation) of the source Jacobian (on the left) and of the 
    four-block split (material failure and three
small submatrices for strain relaxation) of the system (on the right).}
    \label{fig:Jacobiansplit}
\end{figure}

% The results presented here for example employ an approximate Jacobian
% which suppresses the interactions between the $\xi$ sub-system and the $\mathbf{A}$ sub-system and restores
% full coupling through the sub-timestepping scheme and by linearising the source in an iterative manner.
An alternative approach, used in this work, consists in constructing an approximate linearisation of
the source Jacobian so that it is structured in two independent blocks. One is a scalar equation for
$\xi$ obtained by suppressing all the off-diagonal entries of $\Js$. The other, analogously, can be
obtained from removing all dependencies on the \textcolor{black}{damage variable} of the strain-relaxation subsystem.
This way, one can compute the analytical solution for the two subsystems separately, and since the
two blocks are built so that they are independent on each other, these will constitute an
\emph{approximate} solution of \eqref{eq:linearcauchy}. Moreover, being the fracture kinetics
subsystem composed of a single scalar equation, the evaluation of the analytical solution can use
the standard scalar exponential and division operations, instead of the more delicate procedures for
computing matrix exponentials and especially inverting matrices. As an added benefit, if one
approximates the relaxation timescale $\tau_m$ with a constant value (to be recomputed iteratively),
the Jacobian for the strain relaxation subsystem can be evaluated analytically in a very efficient
manner, without resorting to numerical differentiation.

In some particular cases, another approximation step can be taken: whenever the off diagonal entries
of the distortion matrix are small in comparison to the diagonal ones, many of the entries of the
the strain-relaxation source Jacobian become negligible and one can exploit the sparsity of the
strain-relaxation Jacobian matrix to further split it into three blocks. In detail, we take the
material-damage/strain-relaxation source vector
\begin{equation}
    \vec{S}(\vec{Q}) = \trasp{\left(S_\xi,\ S_{A_{11}},\ S_{A_{12}},\ S_{A_{13}},\ S_{A_{21}},\ 
        S_{A_{22}},\ S_{A_{23}},\ S_{A_{31}},\ S_{A_{32}},\ S_{A_{33}}\right)}
\end{equation}
and rearrange it as four independent sources
\begin{align}
    &\vec{S}_\up{a} = \trasp{\left(S_\xi\right)},\\
    &\vec{S}_\up{b} = \trasp{\left(S_{A_{11}},\ S_{A_{22}},\ S_{A_{33}}\right)},\\
    &\vec{S}_\up{c} = \trasp{\left(S_{A_{12}},\ S_{A_{13}},\ S_{A_{21}}\right)},\\
    &\vec{S}_\up{d} = \trasp{\left(S_{A_{23}},\ S_{A_{31}},\ S_{A_{32}}\right)}.
\end{align}
Then of each one of the sources we compute the Jacobian with respect to only the variables in the
corresponding block, that is, for example $\vec{S}_\up{b}$ will be differentiated only with respect
to the diagonal components of $\vec{A}$ and all off-diagonal-block contributions to the global
Jacobian like $\partial S_{A_{11}}/\partial {A_{12}}$ will be assumed null. This approximation is
justified (under the assumption that off-diagonal components of $\vec{A}$ be small) for all
off-diagonal-block derivatives except for $\partial S_{A_{13}}/\partial {A_{31}}$ and $\partial
S_{A_{31}}/\partial {A_{13}}$, which maintain a large magnitude even when $\vec{A}$ is almost
diagonal. These elements can be suppressed regardless, as we already do for all the derivatives like
$\partial S_{\xi}/\partial {A_{ij}}$ and $\partial S_{A_{ij}}/\partial {\xi}$, relying on the
adaptive timestepping method and on the iterative relinearisation for the task of reintroducing the
lost coupling terms. A visual representation of the two-block split (material failure and strain
relaxation) and the four-block split (material failure and three small submatrices for strain
relaxation) of the source Jacobian is given in Figure~\ref{fig:Jacobiansplit}.

\subsubsection{Examples, validation against LSODA and comparison with the implicit Euler scheme} 
\label{sec.odevalidation}

Here we validate our exponential time integrator against available standard software for the
numerical integration of stiff ODE. In particular, we compare against the community standard
LSODA/ODEPACK, see \cite{lsode, odepack}. We also compare our exponential time integrator against a
simple implicit Euler time integration scheme, showing clearly that implicit Euler time stepping,
even with very fine substeps, can lead to significantly inaccurate results.

The benchmark consists in computing the stress-strain diagram for a given choice of spatially homogeneous material
parameters: we distort, with constant strain rate $\dot{\vec{\epsilon}}$, a homogeneous sample of material, which is initially unstressed ($\vec{A} = \vec{I}$) and undamaged ($\xi = 0$). The governing PDE system \eqref{eqn.A} therefore reduces to a simple ODE system. 
The strain rate tensor will have the form 
\begin{equation}
    \vec{\dot{\epsilon}}(t) = 
    \left( 
        \begin{array}{ccc}
            \dot{\epsilon}_{11}(t) & 0 & 0 \\
            0 & 0 & 0 \\
            0 & 0 & 0 \\
        \end{array}
    \right),
\end{equation}
and the resulting governing ODE system for the dynamics reads  
\begin{align}
    & \od{\xi}{t} = -\theta\,\frac{\partial E}{\partial \xi},\\
    & \od{\vec{A}}{t} = -\vec{A}\,\vec{\dot{\epsilon}}(t) - 
        \frac{3}{\tau_m}\,{\left(\det{\vec{A}}\right)}^{5/3}\,\vec{A}\,\dev{\vec{G}},
\end{align}

\paragraph{Brittle and ductile material behavior} 
The aim of this test is not only the validation of the new ODE solver with exponential time
integration, but also the capability of the GPR model with material failure to describe completely
different material behavior with one and the same set of governing equations, just changing the
model parameters appropriately. In particular, the proposed model is able to reproduce brittle
material as well as ductile material behavior.

 \begin{figure}[!b]
    \begin{center}
        \includegraphics[width=0.49\textwidth]{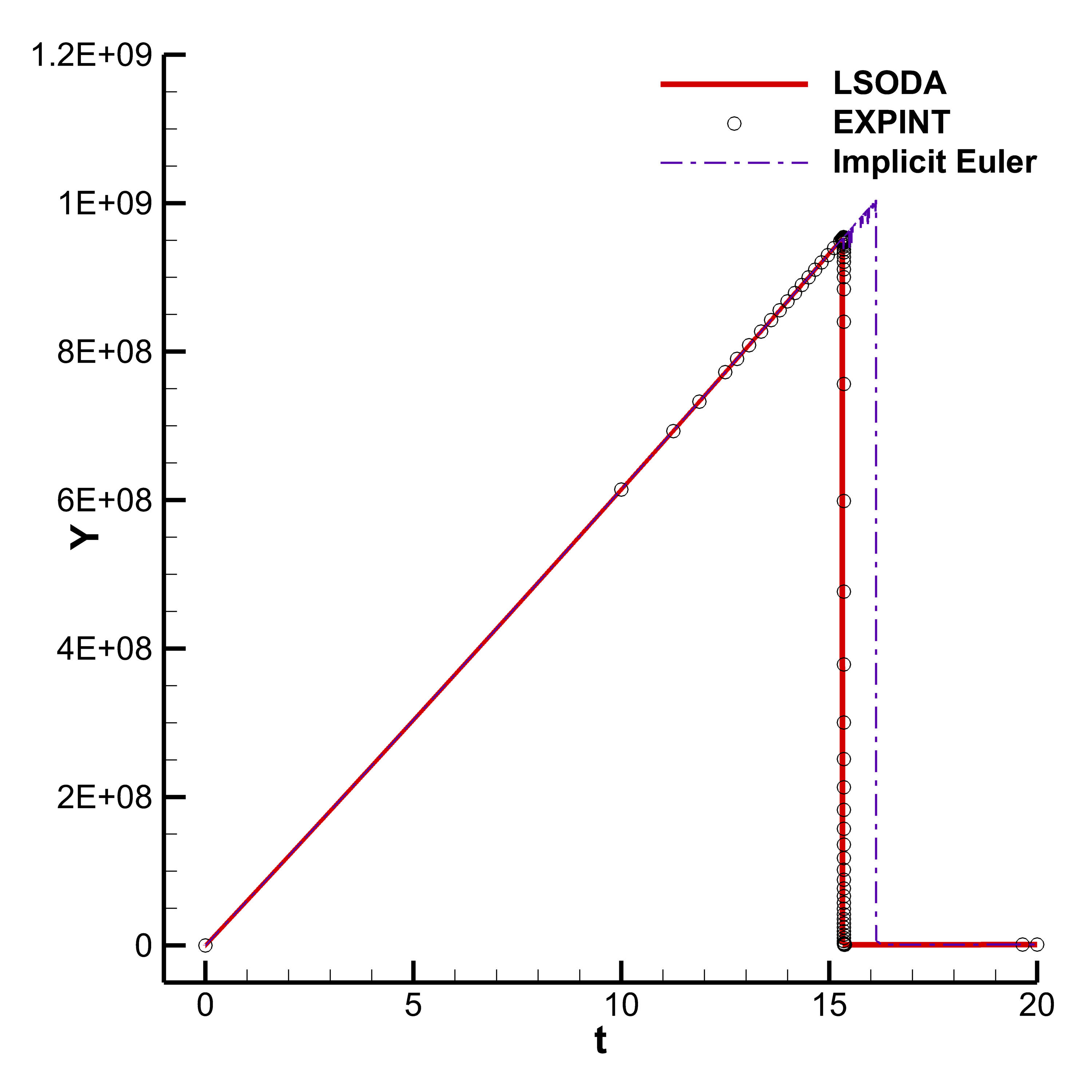}%  
        \includegraphics[width=0.49\textwidth]{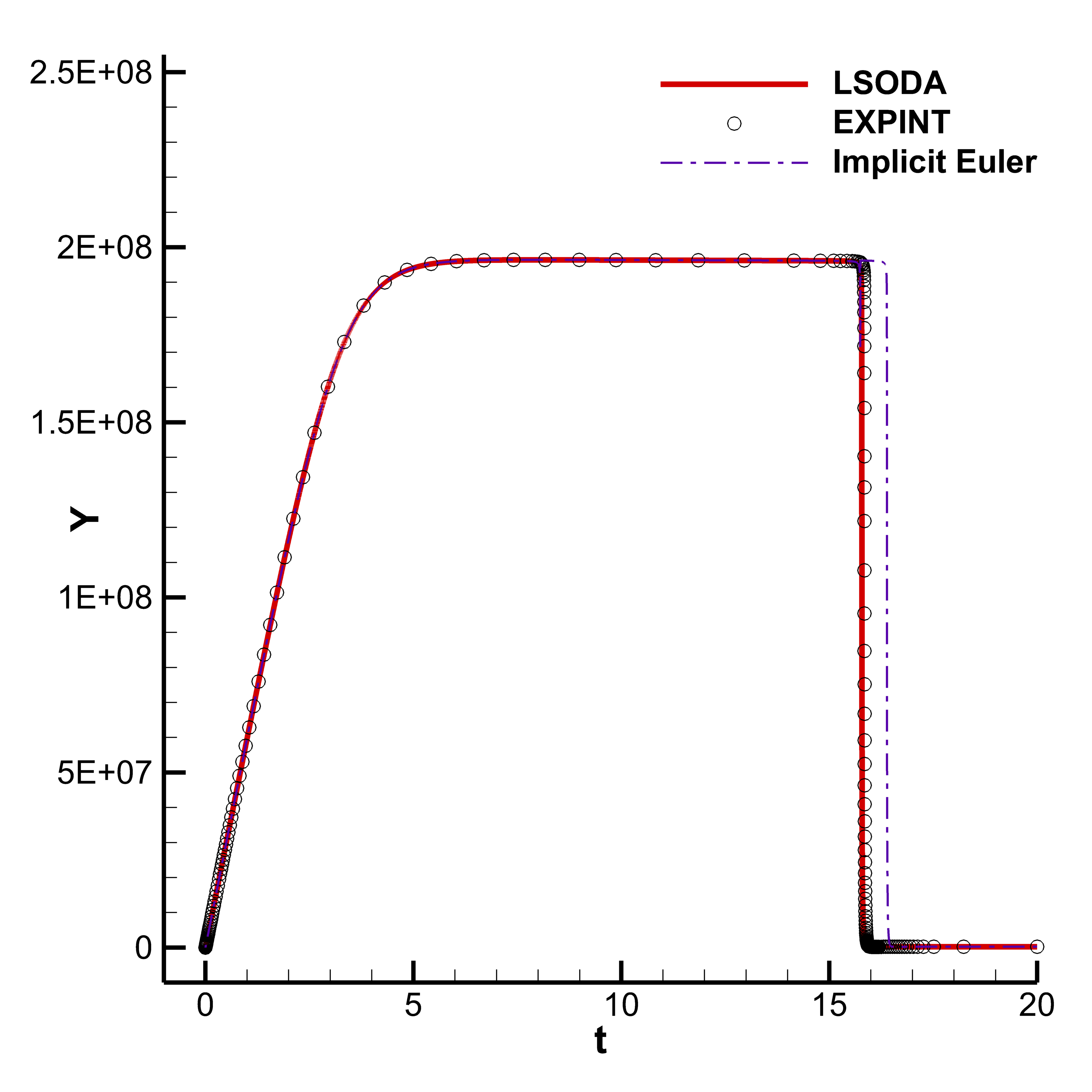} \\ 
        \includegraphics[width=0.49\textwidth]{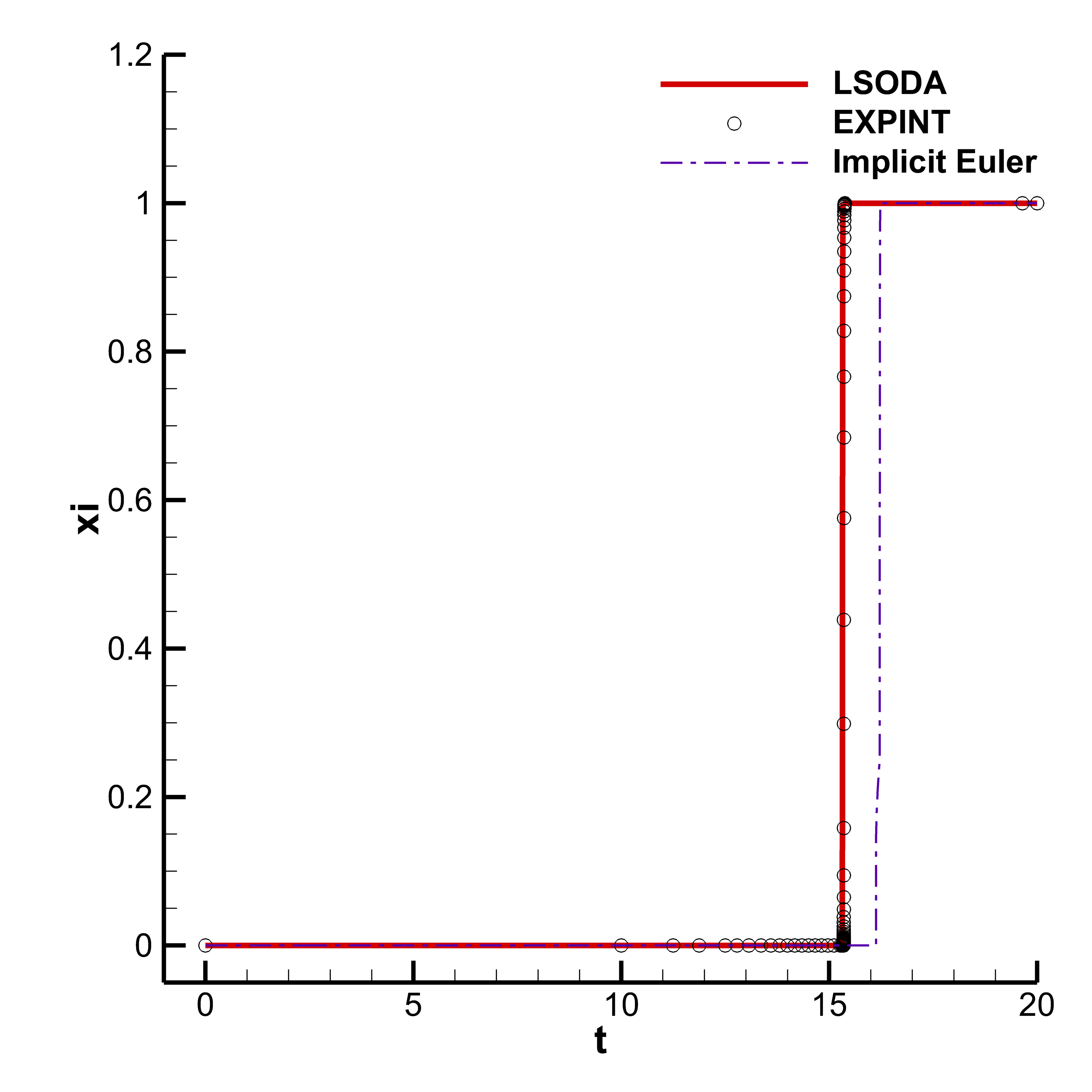}%
        \includegraphics[width=0.49\textwidth]{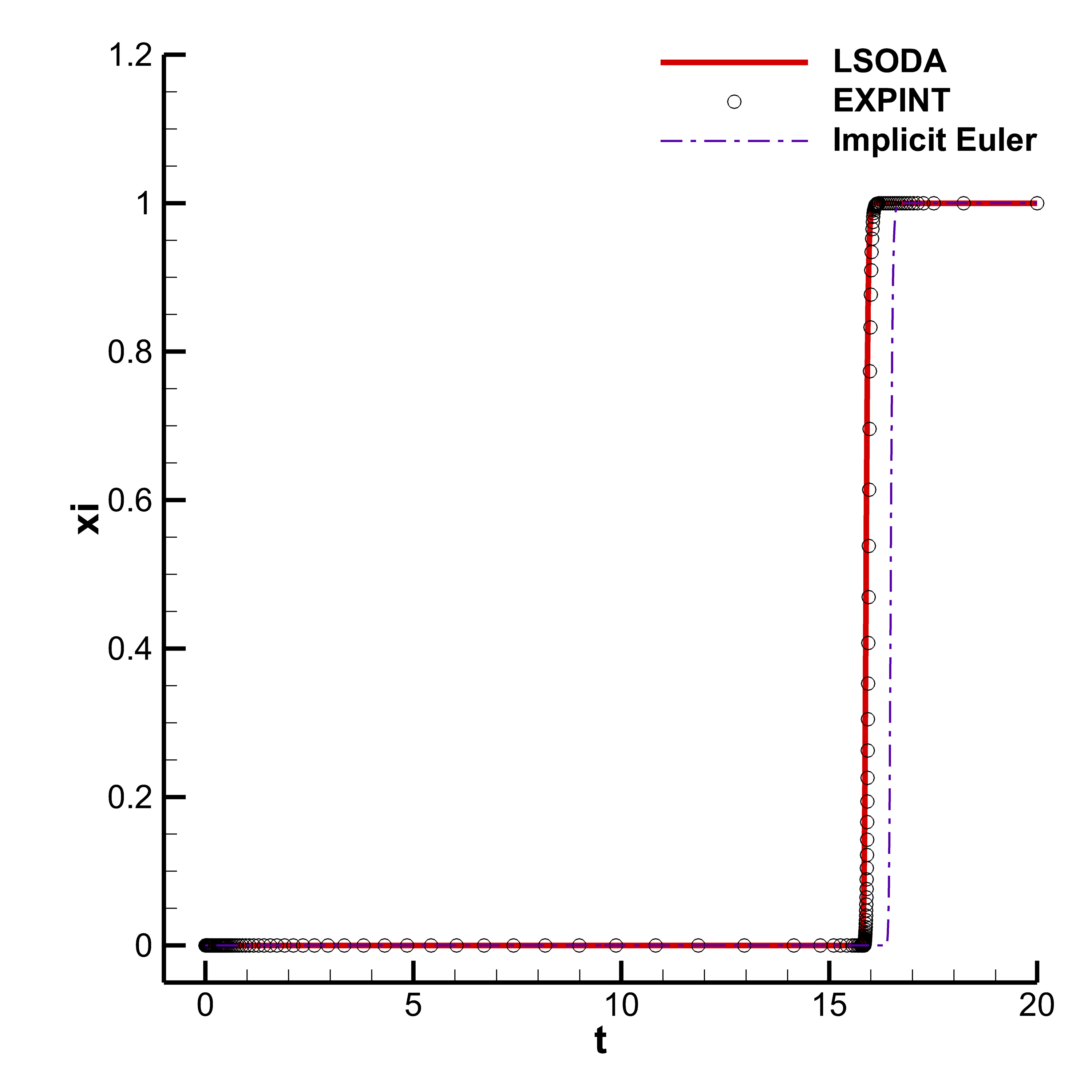} 
    \end{center}
    \caption{Stress-strain diagrams (top row) and time evolution of the \textcolor{black}{damage variable} $\xi$ 
    (bottom row) for an example material characterized by 
    brittle failure (left column) 
    and for a material with ductile failure mode (right column).
    We show that the solution obtained from the proposed exponential integration method matches
    the reference given by the LSODA library, while integration with the simple implicit
    Euler method can yield inaccurate solutions even with very small timesteps.}
    \label{fig.odesolver1}
\end{figure} 

We set the strain rate to $\dot{\epsilon}_{11}(t) = -0.001\,\up{s^{-1}}$ and compute the solution
until a fixed end time $t_\up{end} = 20\,\up{s}$, first with the new exponential integrator with
adaptive timestepping, and then with the $\texttt{scipy.integrate.solve\_ivp}$ Python library
routine, which contains a wrapper to the popular LSODA/ODEPACK Fortran code \cite{odepack}.
Furthermore, we also show the solution obtained by integrating the governing ODE in $10^4$ uniform
timesteps with the implicit Euler scheme.

In Figure~\ref{fig.odesolver1} we plot the von Mises stress $Y = \sqrt{3\,\tr{\left(\dev
\vec{\sigma}\,\dev \vec{\sigma}\right)/2}}$ against time, together with the time evolution of the
\textcolor{black}{damage variable} $\xi$. In a first run (first column of Figure~\ref{fig.odesolver1}) the material
parameters are chosen so that the material exhibits a characteristically brittle behavior (like
ceramics or glass), with linear elastic deformation until the failure point (clearly marked by the
jump in the \textcolor{black}{damage variable} $\xi$), while a second run is carried out with ductile material (like
for example most metals). In this latter case (second column of Figure~\ref{fig.odesolver1}), one
can clearly distinguish an initial linear elastic deformation, followed by a nonlinear transition
into ideal plastic flow until eventual failure.
Our example brittle material is obtained with the following choice of model parameters:
$\rho_0 = 3000\,\up{kg\,m^{-1}}$,
$\mu_I = 30\,\up{GPa}$,
$\mu_D = 30\,\up{MPa}$,
$\lambda_I = \lambda_D = 60\,\up{GPa}$,
$\tau_{I0} = 3 \times 10^{3}\,\up{s}$,
$\tau_{D0} = 3\,\up{s}$,
$\theta_0 = 8$,
$a = 32.5$,
$Y_0 = 1.4\,\up{GPa}$,
$Y_1 = 10\,\up{MPa}$,
$\alpha_I = 35$,
$\alpha_D = 35$,
$\beta_I = 2.2 \times 10^{-8}\, \up{Pa^{-1}}$,
$\beta_D = 2.2 \times 10^{-7}\, \up{Pa^{-1}}$, 
while the ductile behavior is given with
$\rho_0 = 3000\,\up{kg\,m^{-1}}$,
$\mu_I = 30\,\up{GPa}$,
$\mu_D = 30\,\up{MPa}$,
$\lambda_I = \lambda_D = 60\,\up{GPa}$,
$\tau_{I0} = 1 \times 10^{3}\,\up{s}$,
$\tau_{D0} = 1\,\up{s}$,
$\theta_0 = 1$,
$a = 1$,
$Y_0 = 8\,\up{TPa}$,
$Y_1 = 4\,\up{MPa}$,
$\alpha_I = 0$,
$\alpha_D = 30$,
$\beta_I = 2 \times 10^{-8}\, \up{Pa^{-1}}$,
$\beta_D = 1 \times 10^{-4}\, \up{Pa^{-1}}$.

The results from the exponential integrator are in perfect agreement with the LSODA reference
solution, while it is apparent that $10^4$ timesteps with the implicit Euler method are not
sufficient for adequately capturing the sudden onset of material failure.

\paragraph{Rate-dependent behavior and material fatigue}

In this test we want to show that our model is a so-called \textit{rate-dependent model}, that is,
the maximum stress that can be sustained by the material before the failure point is reached, can
vary as a function of the speed of deformation, i.e. the strain rate. It is quite common for example
that impacts or explosions allow a given material to achieve a higher equivalent stress than slow
(quasi-static) loads, before total breakdown.

 \begin{figure}[!b]
    \begin{center}
        \includegraphics[width=0.49\textwidth]{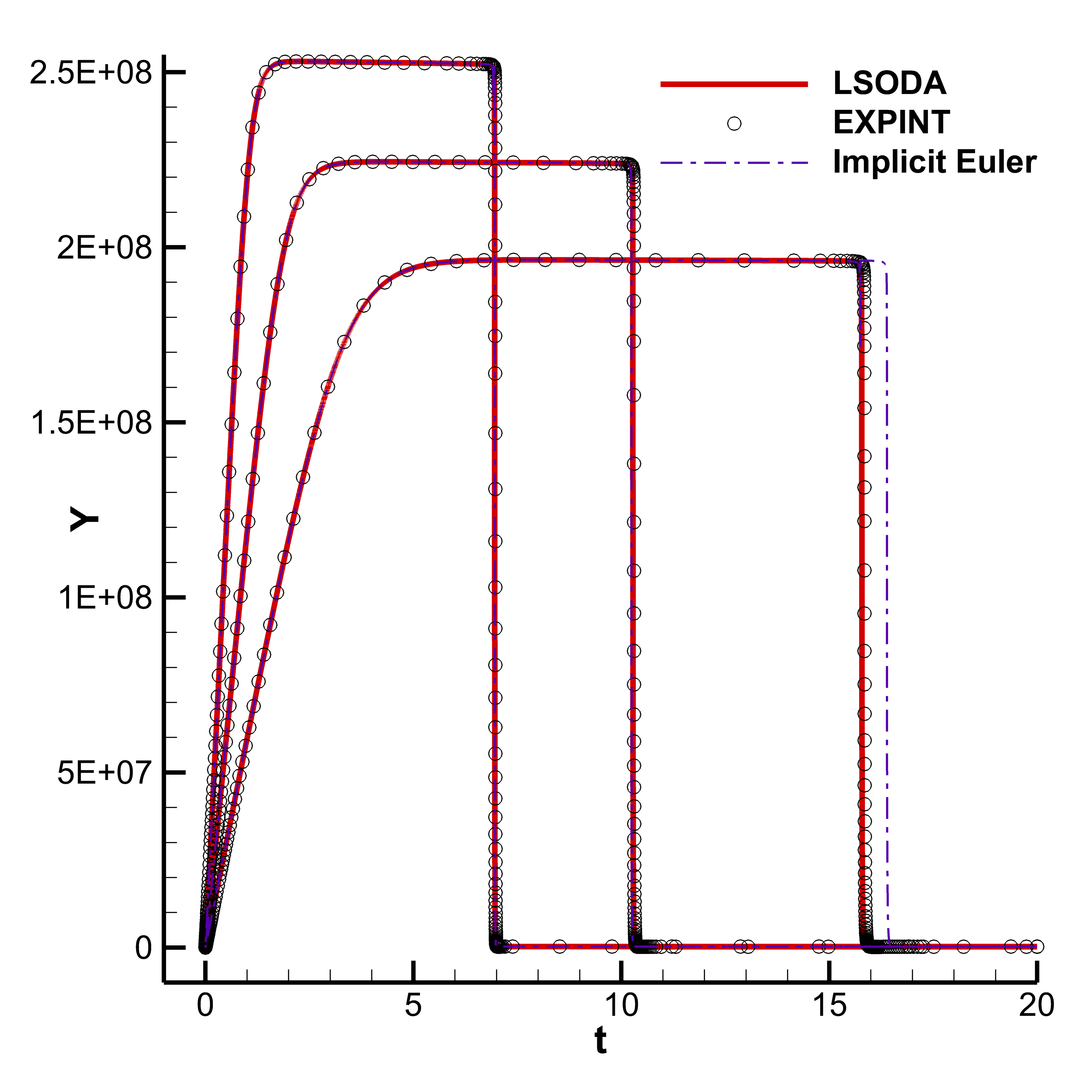}%
        \includegraphics[width=0.49\textwidth]{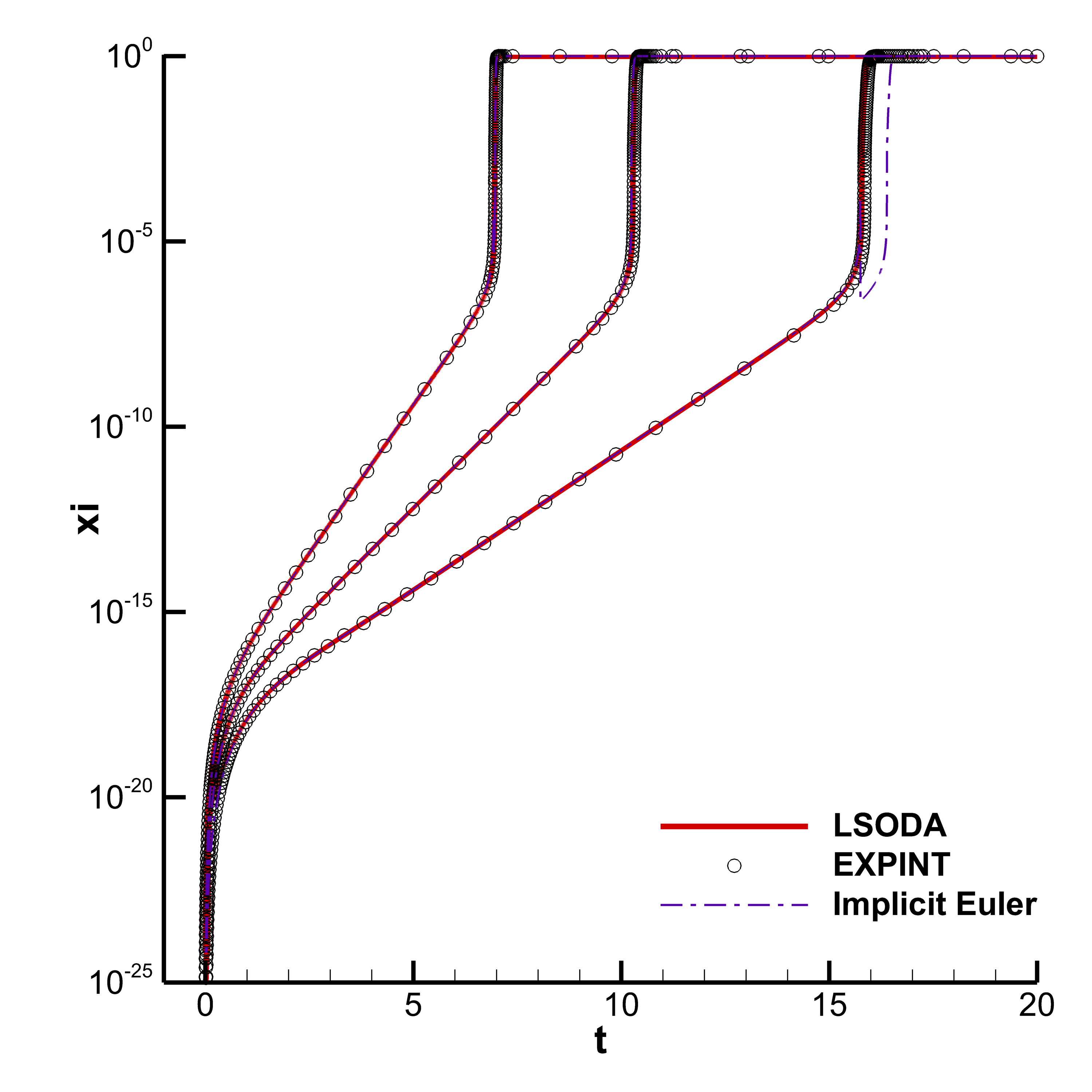}%
    \end{center}
    \caption{Rate-dependent stress-strain diagrams (left panel) and time evolution of 
    the \textcolor{black}{damage variable} (right panel) for a ductile material subjected to three different 
    constant strain rates.}
    \label{fig.odesolver2}
\end{figure} 

To show that effects of this type can be reproduced by our model, we repeat the previous test on the
ductile material, but by varying the strain rate from $\dot{\epsilon}_{11}(t) = -0.001\,\up{s^{-1}}$
to $\dot{\epsilon}_{11}(t) = -0.002\,\up{s^{-1}}$ and $\dot{\epsilon}_{11}(t) =
-0.004\,\up{s^{-1}}$, reaching higher maximum equivalent stress states as strain is applied faster.
This is shown in Figure~\ref{fig.odesolver2}, where we also report the time evolution of the damage
parameter $\xi$, with logarithmic scaling, in order to highlight that very small deviations of this 
variable from the perfectly intact state (of the order of $10^{-15}$ to $10^{-5}$) can indicate the
difference between the Hookian elastic response and the plastic regime.

Finally, we want to show that the proposed model can also describe \textit{fatigue} effects by allowing the
material to retain memory of previous stress states. In fact, thermodynamical consistency prescribes
that the source term associated with the \textcolor{black}{damage variable} always be positive (or possibly
vanishing), this means that any stress applied to the material will ever so slightly damage it and
eventually cause a deterioration of its mechanical properties. The extremely nonlinear response to
different intensity levels of the stress norm determines that appropriate choices of parameters can
describe materials with very late onset of fatigue effects, as well as properly governing the
deterioration process following a certain level of exposure to cyclic stress. 

 \begin{figure}[!b]
    \begin{center}
        \includegraphics[width=0.49\textwidth]{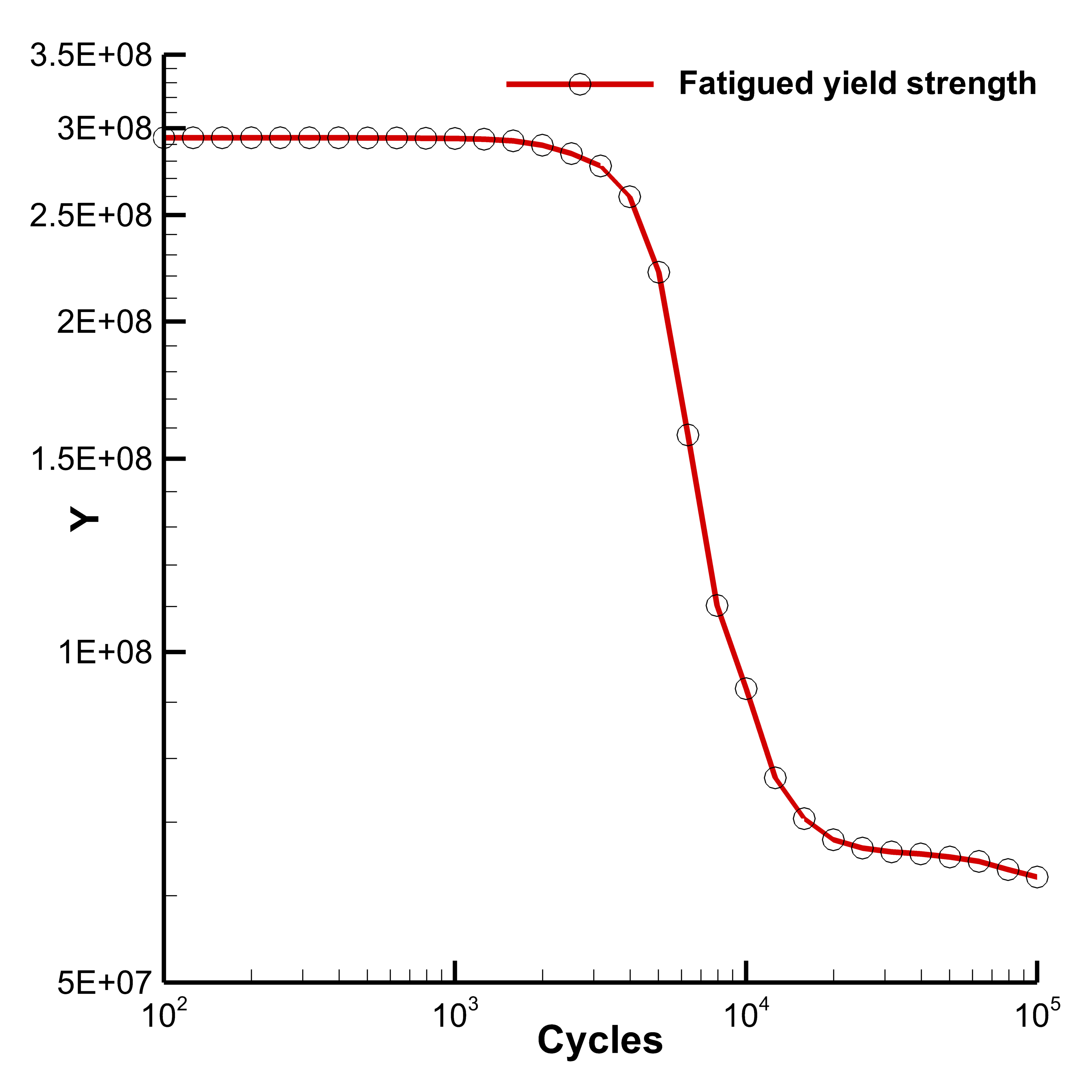}%
        \includegraphics[width=0.49\textwidth]{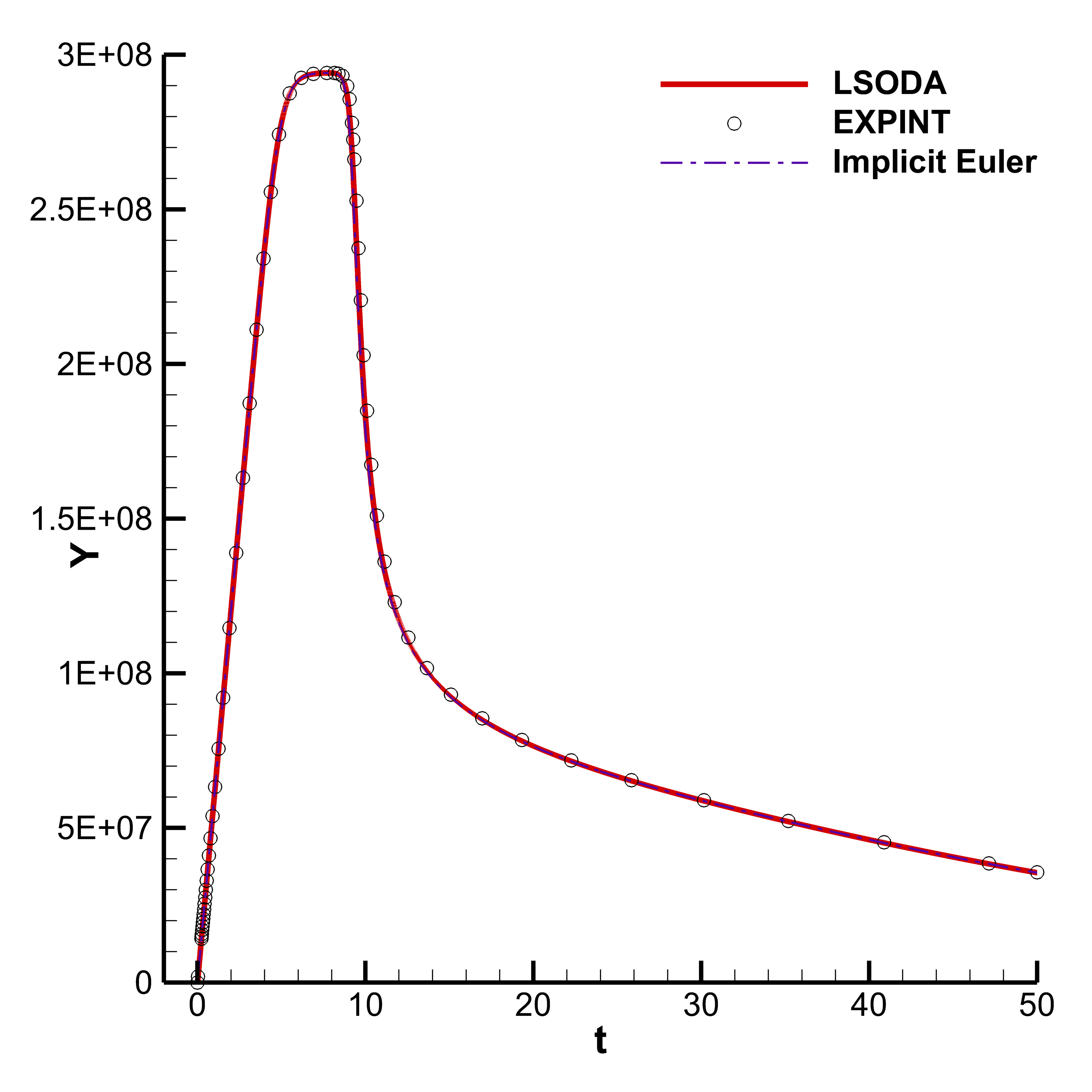} \\ 
        \includegraphics[width=0.49\textwidth]{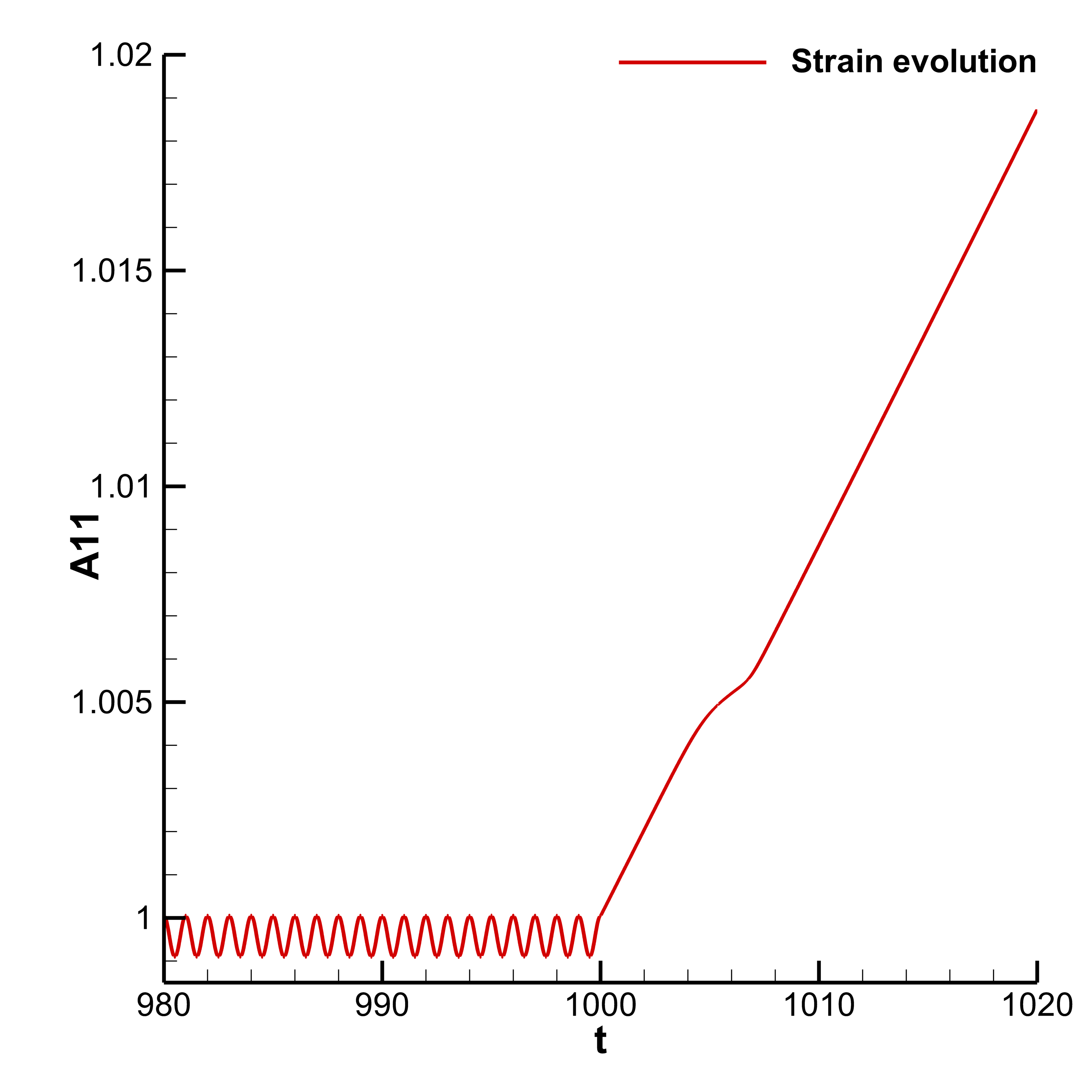}%
        \includegraphics[width=0.49\textwidth]{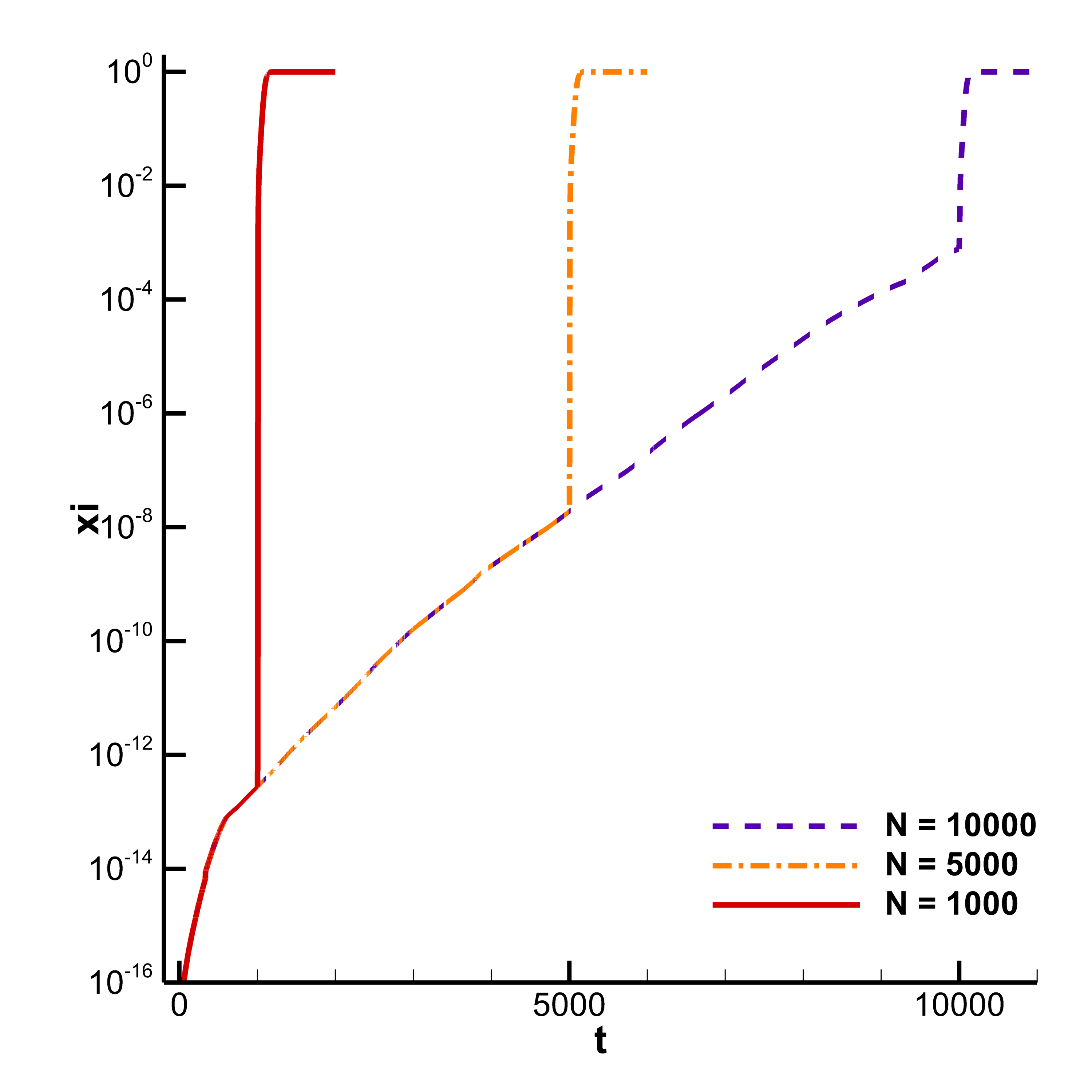} %  
    \end{center}
    \caption{Numerical results obtained for the fatigue behaviour. Top left panel: W\"ohler diagram, showing the weakening of the material strength as a function of the number of applied load cycles. Top right panel: stress-strain diagram of the test material for a single, quasi-static (slow) load cycle until rupture.           
    Bottom left panel: evolution of the $A_{11}$ component of the distortion field for the last elastic load cycles in the case $N=1000$ and final traction test until material failure.   
    Bottom right panel: temporal evolution of the \textcolor{black}{damage variable} $\xi$ for three different numbers of load cycles $N$. }
    \label{fig.odesolver3}
\end{figure} 

For this test we employ an elasto-plastic material characterized by the following choice of parameters:
$\rho_0 = 3000\,\up{kg\,m^{-1}}$,
$\mu_I = 30\,\up{GPa}$,
$\mu_D = 33\,\up{MPa}$,
$\lambda_I = \lambda_D = 60\,\up{GPa}$,
$\tau_{I0} = 2 \times 10^{5}\,\up{s}$,
$\tau_{D0} = 2 \times 10^{3}\,\up{s}$,
$\theta_0 = 1$,
$a = 1$,
$Y_0 = 8\,\up{TPa}$,
$Y_1 = 8\,\up{TPa}$,
$\alpha_I = 0$,
$\alpha_D = 0$,
$\beta_I = 3 \times 10^{-8}\,\up{Pa^{-1}}$,
$\beta_D = 0\,\up{Pa^{-1}}$.  
We setup a series of computations by first subjecting the material to a low intensity cyclic stress
(about 20\% of the limit for the elastic regime), determined by a variable strain rate of the form
$\dot{\epsilon}_{11}(t) = -0.001\,\sin{\left({2\,\pi\,t}\right)}\,\up{s^{-1}}$, and then measuring
the residual strength of the material by means of a quasi-static destructive test as those shown in
Figures~\ref{fig.odesolver1} and \ref{fig.odesolver2}. A plot of the material distortion caused by
this type of forcing (low intensity periodic deformations followed by a destructive test) is shown in Figure~\ref{fig.odesolver3}, together with the numerical stress-strain
diagram of an initially intact sample under slow (quasi-static) loading. The main result of this benchmark is the so-called \textbf{W\"ohler diagram}: for this purpose, in the top left panel of Figure~\ref{fig.odesolver3} we show a bilogarithmic plot of the material strength as a function of the number of applied cycles. The W\"ohler diagram obtained in our simulations shows the same qualitative behaviour as those obtained in experiments for real materials, with the classical division into three different regimes: short, intermediate and long-term durability. Our numerical results illustrate that the test material shows no signs of weakening at a low number of load cycles (up to about 1000). For more load cycles, one can note the typical decay of the material resistance, while for a very large number of cycles the the long-time durability limit is reached, which, however, still shows a small but steady degradation, as it is commonly found in aluminum-alloys. 
In order to give the reader a more detailed insight about how fatigue is accounted for in our model, we also show the temporal evolution of the \textcolor{black}{damage variable} $\xi$ (bottom right panel) for three of the tests (one at $N = 1000$ elastic cycles, one
at $N = 5000$ cycles, one at $N = 10000$ cycles), showing that the accumulation of fatigue effects is well captured by the exponential growth of the \textcolor{black}{damage variable}. 
 
At this point we would like to emphasize that the addition of the \textit{thermodynamically compatible} evolution equation \eqref{eqn.xi} of a simple scalar $\xi$ to the GPR model \cite{PeshRom2014,DumbserPeshkov2016} is enough to model such a complex behaviour as material fatigue.

\clearpage 

%===========================================================================
%==========         N U M E R I C A L         R E S U L T S 
%
%\subsection{Fitting of the model parameters with experimental stress-strain diagrams}  
%\label{sec.stress.strain}

\section{Numerical results}       
\label{sec:results}

\textcolor{black}{If not stated otherwise, throughout the entire section all units are in SI units. In the case of adaptive mesh refinement, the AMR refinement factor between two adjacent grid levels is always $\mathfrak{r}=3$, see \cite{AMR3DCL,Zanotti2015a,Reinarz2020} for details on AMR. }

\subsection{Stiff inclusion test case}  
\label{sec.stiff.inc}
The first test that we show here was originally proposed by LeVeque \cite{LeVeque:2002a} in Chapter 22.7 for linear elasticity. The test consists in a p-wave that travels through a heterogeneous medium  with free surface boundary conditions on the top and bottom boundaries of the computational domain.
In the original test, an outer solid medium $\Omega_{out}=[-1,1]\times [-0.5,0.5]$ contains another stiff material placed in $\Omega_{in}=[-0.5,0.5]\times [-0.1,0.1]$. The Lam\'e parameters for the two materials are $(\lambda, \mu, \rho)_{out}=(2,1,1)$ and $(\lambda, \mu, \rho)_{out}=(200,100,1)$ so that the resulting p- and s- wave speeds in the stiffer material are 10 times larger than in the other material.

\begin{figure}[!bp]
    \begin{center}
        \includegraphics[width=0.7\textwidth]{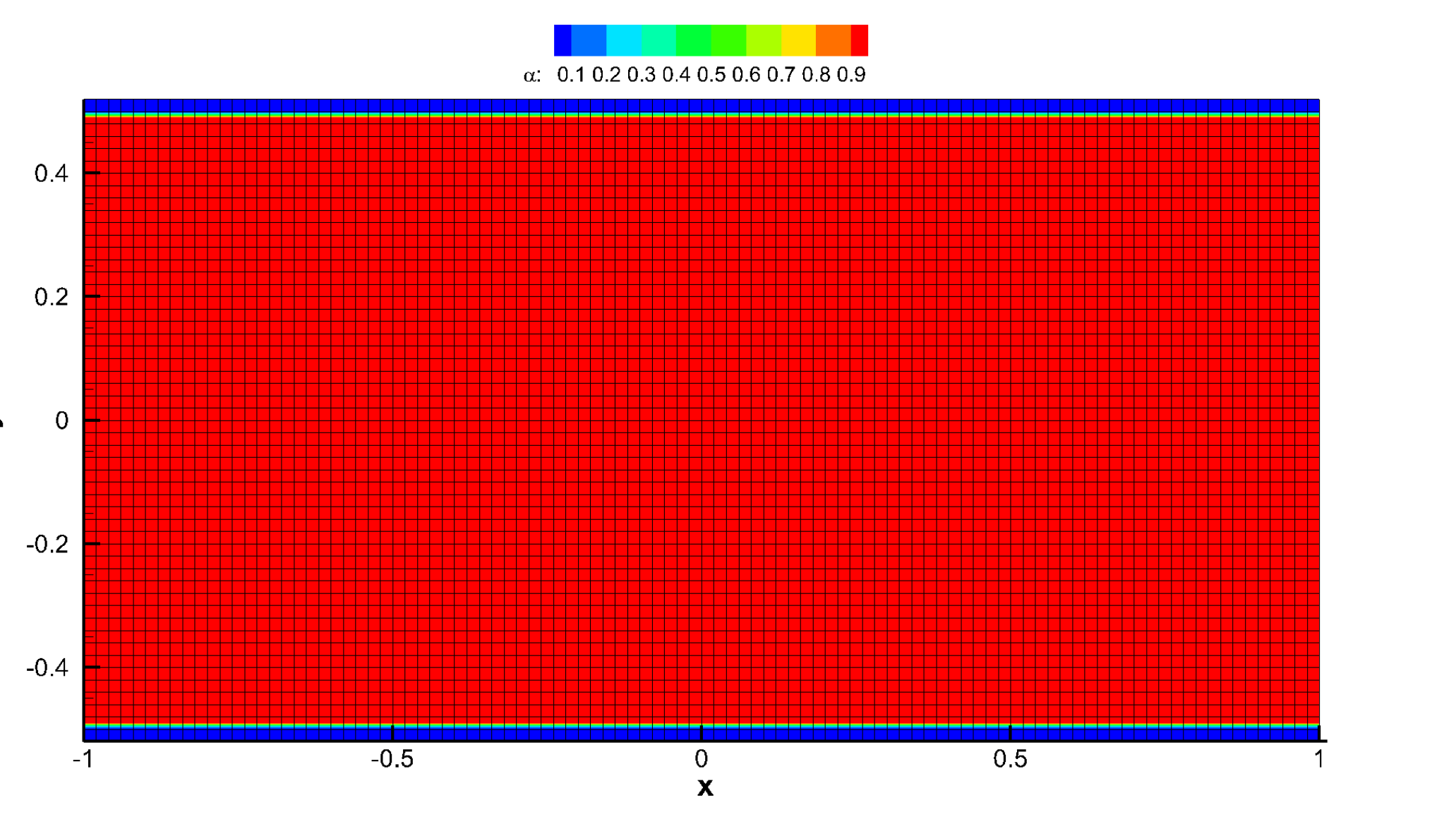}
    \end{center}
    \caption{Shape of the diffuse interface, represented by the solid volume fraction function $\alpha$.} 
    \label{stiff_alpha}
\end{figure}

\textcolor{black}{The main purpose of this test problem in the context of nonlinear hyperelasticity is to verify if the numerical method generates spurious pressure oscillations at the \textit{moving} interface where the Lam\'e parameters jump by two orders of magnitude, or not. It is well known from conservative numerical methods applied to the compressible Euler equations that spurious pressure oscillations are generated when solving the Euler equations with a spatially variable ratio of specific heats, see e.g. \cite{ToroAnomalies,AbgrallQuasiCons} for a more detailed discussion. Note that in our model the PDE for the transport of the Lam\'e constants and the parameter $Y_0$ is written in \textit{non-conservative} form, see eqn. \eqref{eqn1b}. }

 The initial condition for the travelling p-wave is given by
 \begin{eqnarray}
    \boldsymbol{\sigma} &=&\epsilon R_\sigma \cdot \mbox{exp}\left({-\frac{1}{2}\frac{(x-x_0)^2}{\delta^2}}\right) \nonumber \\
    \vec{v} &=&\epsilon R_{\vec{v}} \cdot \mbox{exp}\left({-\frac{1}{2}\frac{(x-x_0)^2}{\delta^2}}\right)
 \end{eqnarray}
where $\boldsymbol{\sigma} = (\sigma_{xx},\sigma_{yy},\sigma_{zz},\sigma_{xy},\sigma_{yz},\sigma_{xz})$, 
$R_\sigma=(4,2,2,0,0,0)$, $R_{\vec{v}}=(-2,0,0)$, $x_0=-0.8$ and $\delta=0.01$ represents a 
p-wave of Gaussian shape travelling in the $x$-direction. Since the test regards a linear elastic medium, we chose $\epsilon=10^{-4}$ in order to have small deformations and negligible contributions from the nonlinear convective part. The initial condition for the matrix $\mathbf{A}$ is obtained by fixing a direction angle $\theta=0$ and then use the procedure $\mathbf{A}=\mathbf{A}(\boldsymbol{\sigma},\theta)$ as described in \ref{sec.init.A}. The thermal impulse is set to 
$\mathbf{J}=0$ and furthermore in this test we set $c_h=0$. 

Finally the free surface boundary conditions are obtained by enlarging the computational domain to $\Omega=[-1.0,1.0]\times [-0.52,0.52]$ and then use the diffuse interface parameter $\alpha$ in order to identify the location of the solid medium, i.e. setting $\alpha=1$ for $\mathbf{x} \in \Omega_{out}$ and $\alpha=0$ elsewhere. In order to prevent spontaneous ruptures, both materials are assumed to be unbreakable, 
i.e. the yield stress is set to the very large value $Y_0=10^{22}\,\up{Pa}$. The domain is covered using $100 \times 52$ elements and we use fourth order polynomials $N=M=4$ in space and time.

The inclusion of the constraint $\rho / \rho_0 =|A|$ is crucial in order to obtain a solution due to the heterogeneity of the medium. 

\begin{figure}[!bp]
    \begin{center}
        \includegraphics[width=0.49\textwidth]{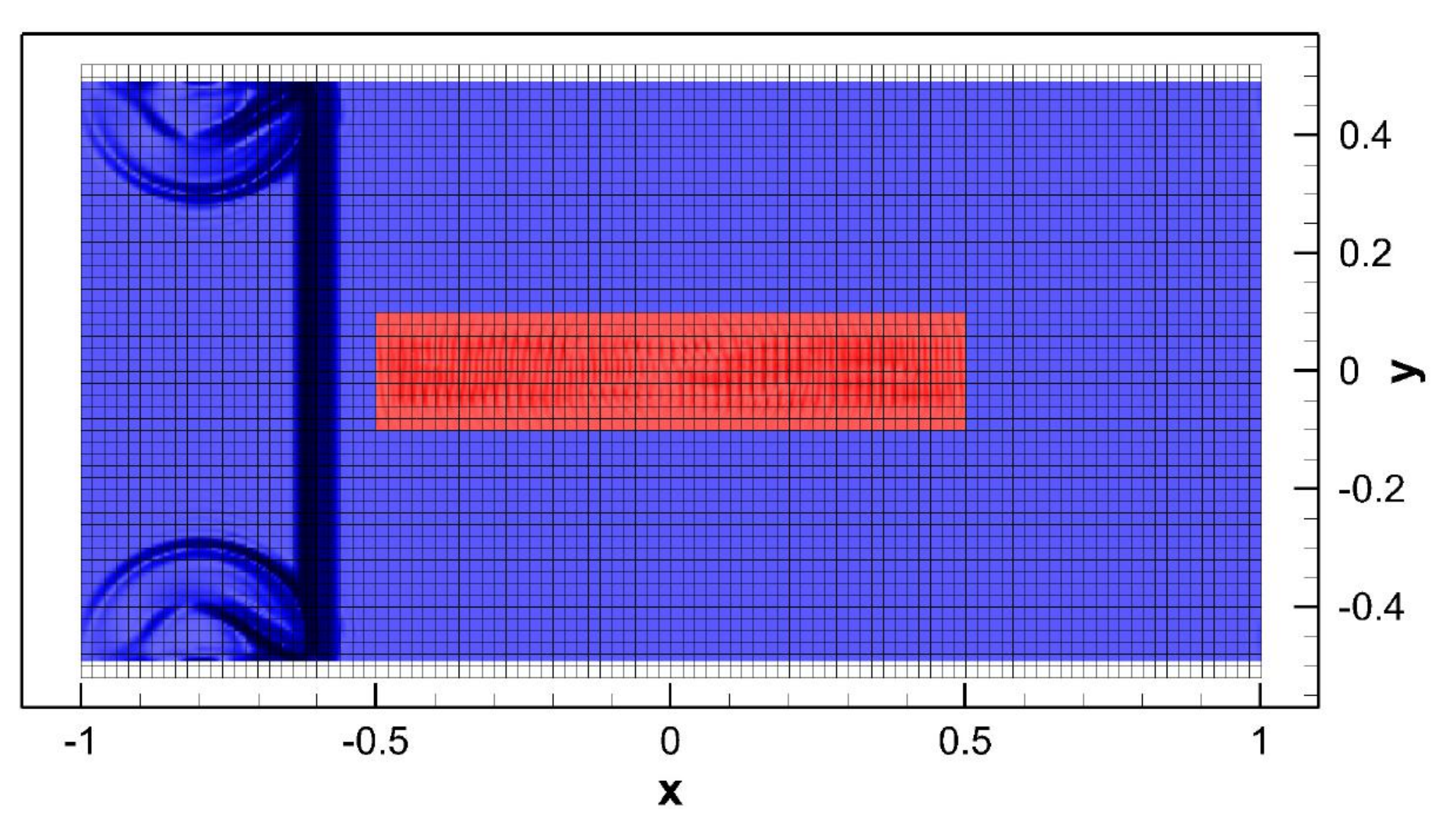}
        \includegraphics[width=0.49\textwidth]{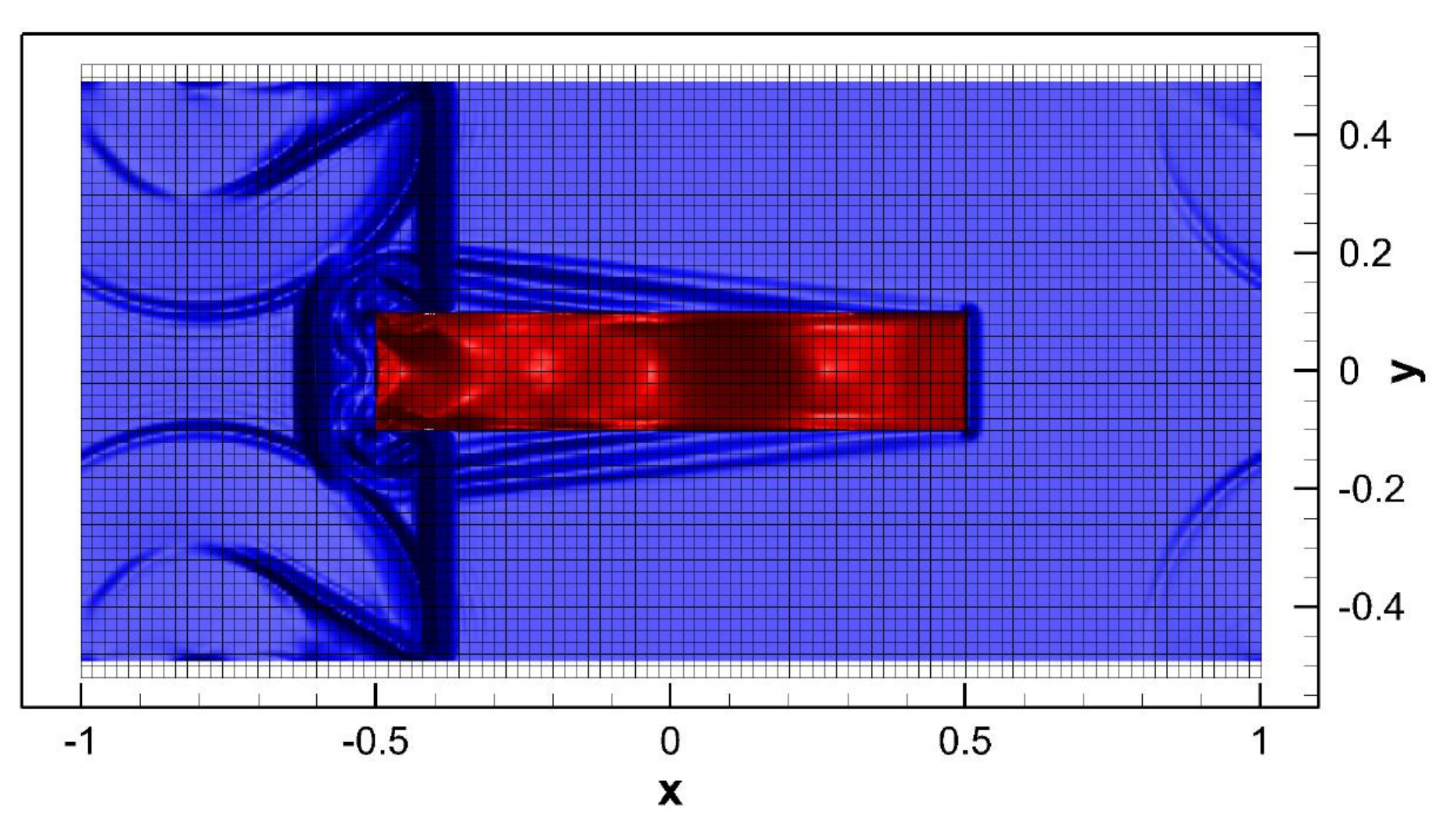}
        \includegraphics[width=0.49\textwidth]{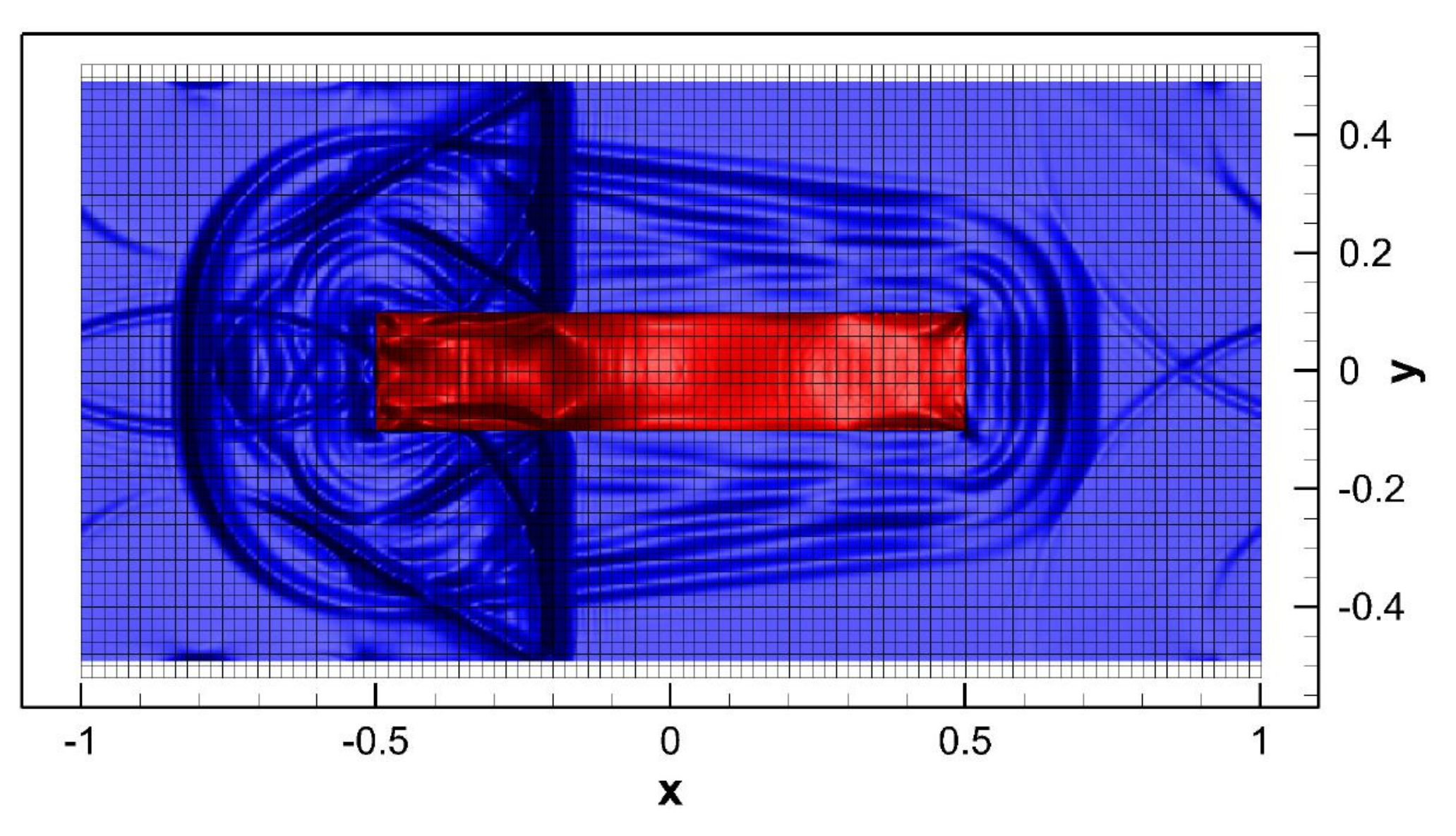}
        \includegraphics[width=0.49\textwidth]{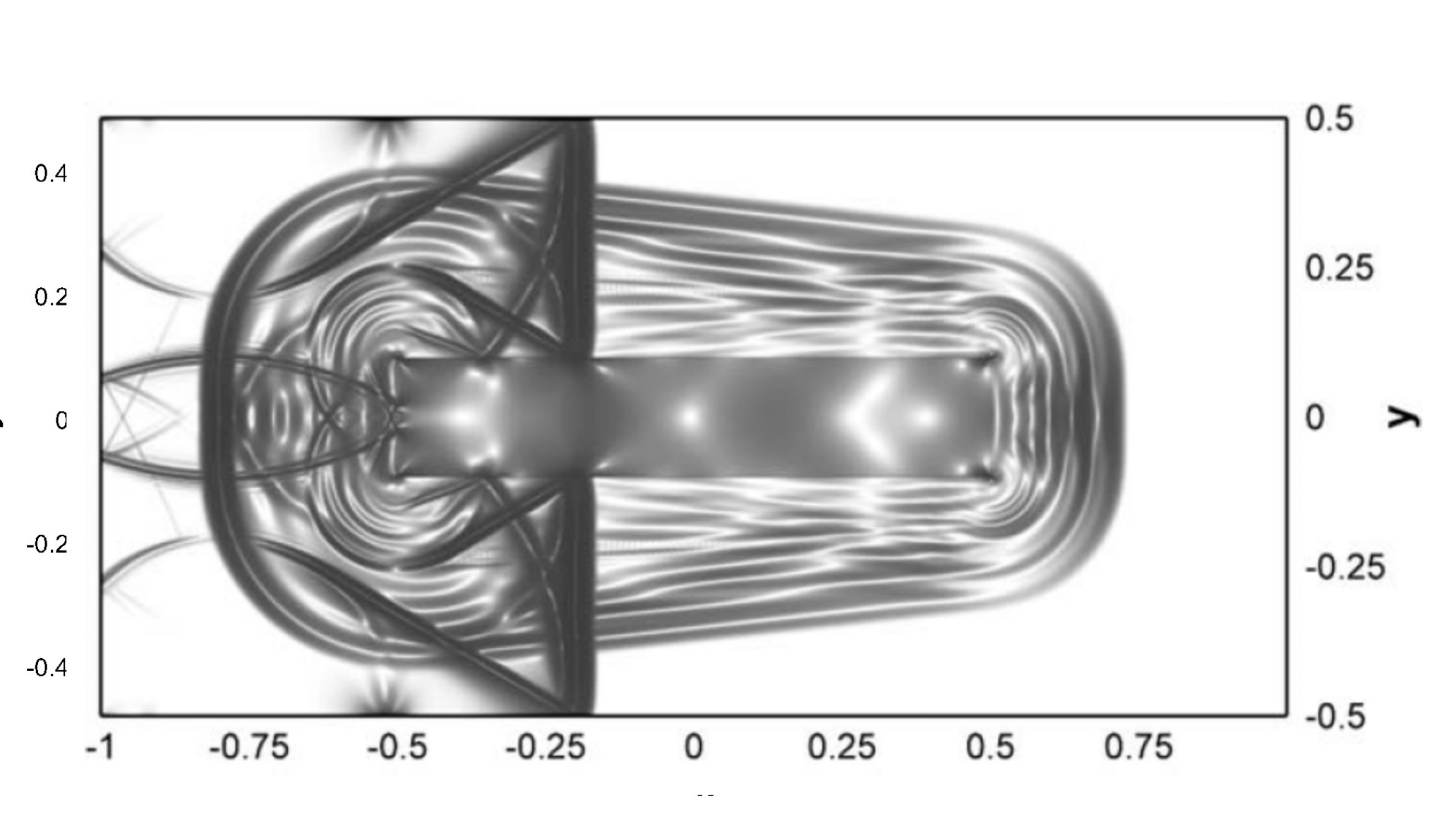}
    \end{center}
    \caption{Evolution of the stress component $\sigma_{xx}$ at times $t=0.1,0.2$ and $t_{end}=0.3$ colored 
    with $\lambda$. In the bottom right figure we report the numerical solution obtained in \cite{gij1} 
    with an ADER-DG scheme of fourth order on $4 \cdot 200 \times 100$ triangular elements.}
    \label{stiff_timeevol}
\end{figure}

The resulting profile of the diffuse interface $\alpha$ is reported in Fig. \ref{stiff_alpha}, 
while the time evolution of the stress component $\sigma_{xx}$ is reported in Fig. \ref{stiff_timeevol}. 
In Fig. \ref{stiff_timeevol} one clearly can observe the influence of the free surface in the 
propagation of the p-wave, before it hits the stiff inclusion. Note that in our diffuse interface 
approach, the free surface boundary condition is fully taken into account by the governing PDE system 
and the scalar volume fraction function $\alpha$. No explicit boundary conditions need to be applied, 
see also \cite{Tavelli2019}. Then the pattern generated after the hit with the stiff material results 
very similar to the one obtained in  \cite{gij1} with an ADER-DG scheme of fourth order of accuracy 
on a boundary-fitted mesh of $4 \cdot 200 \times 100$ triangular elements. 

\textcolor{black}{From the obtained computational results one clearly observes that the produced wavefield is free from spurious pressure oscillations at the interface and that therefore also jumps in the Lam\'e constants by two orders of magnitude can be properly handled by our numerical method. }

\subsection{Fracture generation after high speed impact of a copper plate onto a pyrex glass}  
\label{sec.copper.pyrex}

\begin{figure}[!bp]
    \begin{center}
        \includegraphics[width=0.7\textwidth]{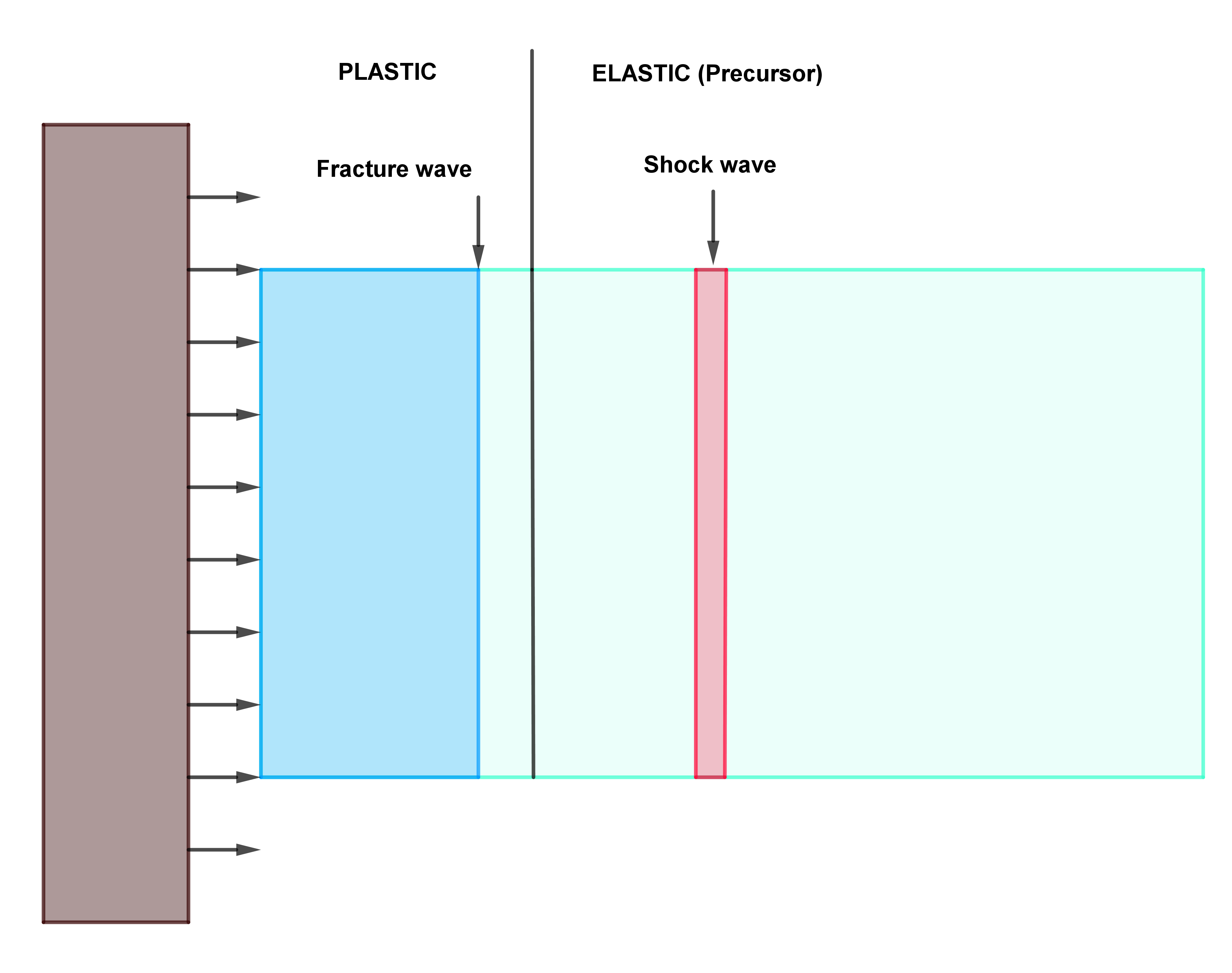}
    \end{center}
    \caption{Sketch of the experimental setup. }
    \label{fig.exp1}
\end{figure}

Here we want to study the wave and rupture propagation driven by the impact of a flying copper plate that impacts on a Pyrex glass brick, following the setup proposed in \cite{Resnyansky2003}.  We consider first a $1D$ impact problem of a copper plate onto a pyrex glass. A sketch of the experimental setup is depicted in Fig. \ref{fig.exp1}. 
As discussed in \cite{Resnyansky2003}, a proper equivalent stress $Y$ 
has to include both shear and pressure contributions. For this numerical experiment we use the following linear definition of the equivalent stress:
\begin{equation} \label{eq.yieldcriterion}
    Y=A\,Y_s + B\,|Y_p|
\end{equation}
\textcolor{black}{where, $Y_s = \sqrt{3\,\tr{\left(\dev
\vec{\sigma}\,\dev \vec{\sigma}\right)/2}}$ is the von Mises shear stress norm, 
$Y_p = \tr{\vec{\Sigma}}/3$ accounts for the pressure contribution, 
and the weighting factors are set to $A=0.9$ and $B=0.05$.}

As numerical parameters we take a computational domain $\Omega=[-0.001,0.021]\times[-0.001,0.001]$
covered with $108 \times 8$ elements, polynomial approximation degree $N=M=3$ and periodic boundaries in the $y$-direction.
 We use a multi-material approach where the mechanical properties inside the domain are set as 
 unstressed Pyrex, while copper is located at the left boundary (see Table~\ref{tab.materials} for material parameters), where also an impact 
 velocity $v_B$ is assumed.  For this test we use two different impact velocities of 
 $v_1=250\,\up{m/s}$ and $v_2=530\,\up{m/s}$, as done in the experimental setup, see \cite{Bourne1994}. 
We take the time series of the stress, computed as $\sigma=|\sigma_{xx}|$, in $x=2.5\,\up{mm}$.

\begin{figure}[!bp]
    \begin{center}
        \includegraphics[width=0.7\textwidth]{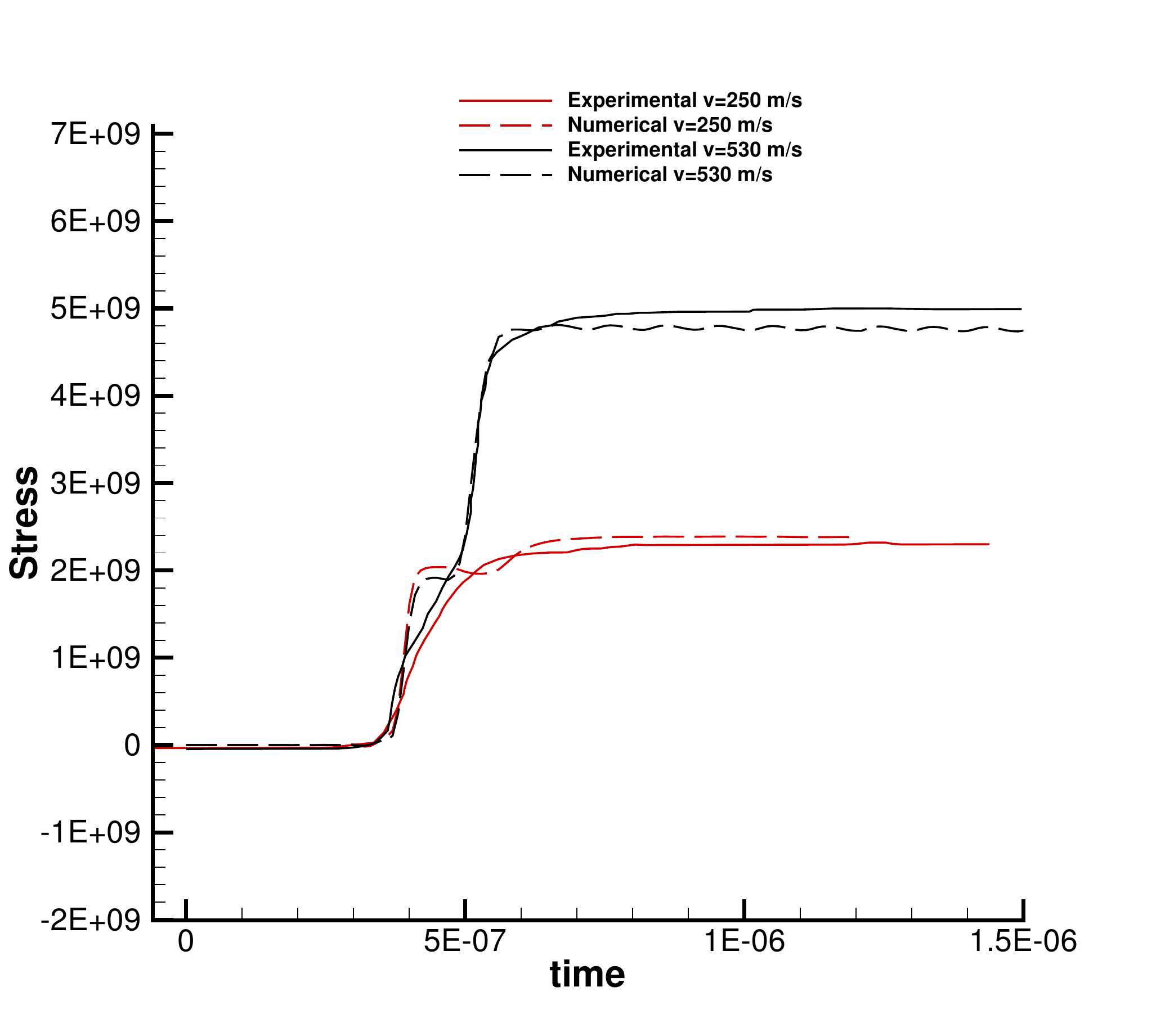}
    \end{center}
    \caption{Time series of the computed stress defined as $\sigma=|\sigma_{xx}|$ compared with 
    experimental results \cite{Resnyansky2003,Bourne1994} for two different impact velocities, namely $v_1=250\,\up{m/s}$ and $v_2=530\,\up{m/s}$. }
    \label{fig.exp1_2}
\end{figure}

The comparison between numerical and experimental stress, properly shifted according to 
the $p-$wave velocity, is reported in Figure $\ref{fig.exp1_2}$. Note here that the use of 
multi-material is crucial in order to obtain the proper velocity inside the target object. 
The resulting penetration velocity for the two cases are, respectively 
$v_{p,1} \approx 426\,\up{m/s}$ and $v_{p,1} \approx 187\,\up{m/s}$ and strictly depend on the mechanical properties of the bullet. 
We now move to a two dimensional experiment where a hemispherical projectile hits a 
Pyrex glass block. According to \cite{Resnyansky2003} we need to generalize the equivalent stress definition as follows
\begin{equation}
Y=B |Y_p|+A \frac{1}{2}\left(Y_s-s_0 \right)\left[1+\textnormal{erf}\left(\frac{Y_s-s_0}{\epsilon}\right)\right], 
\end{equation}

\begin{figure}[!bp]
    \begin{center}
        \includegraphics[width=0.32\textwidth]{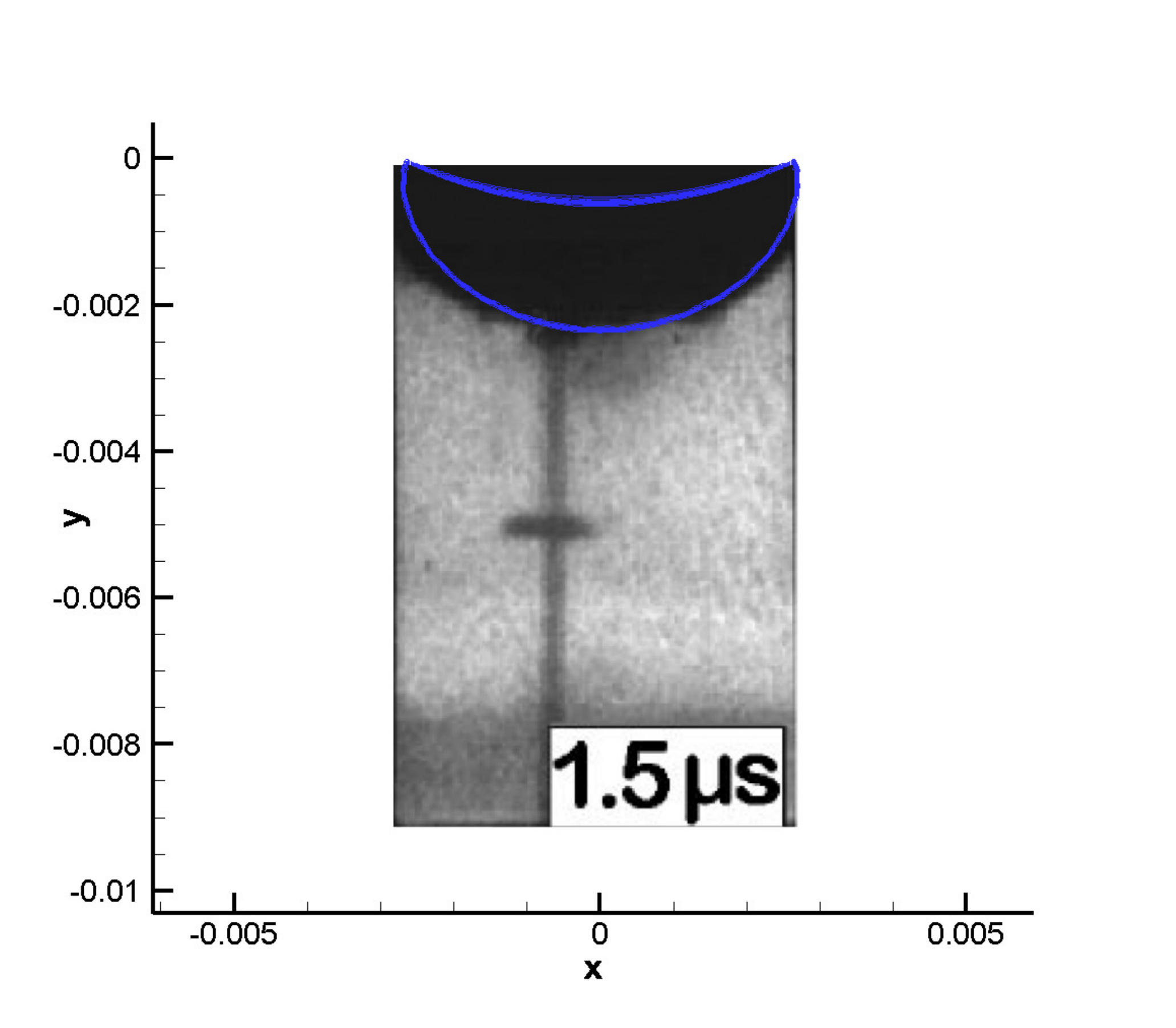}
        \includegraphics[width=0.32\textwidth]{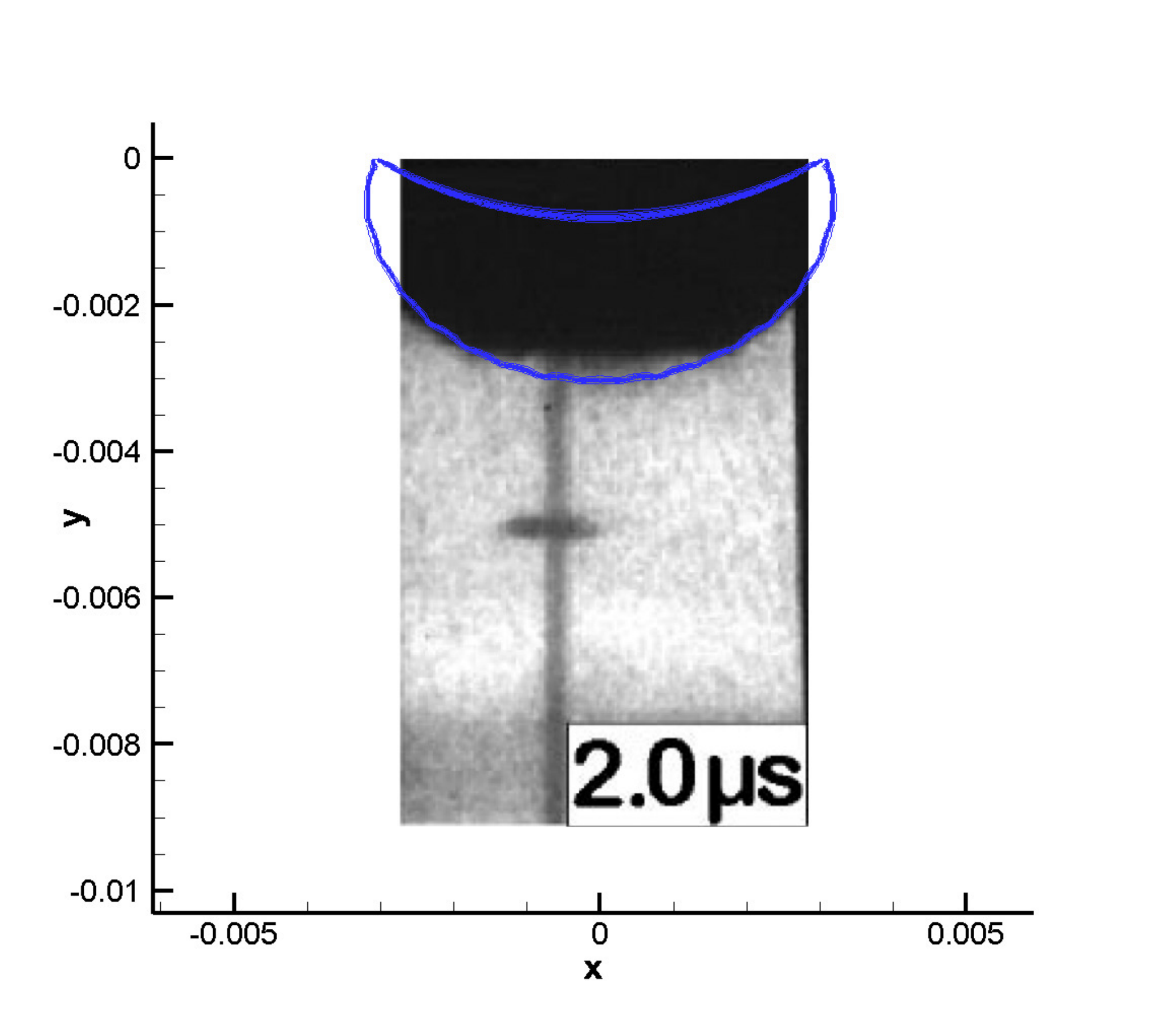}
        \includegraphics[width=0.32\textwidth]{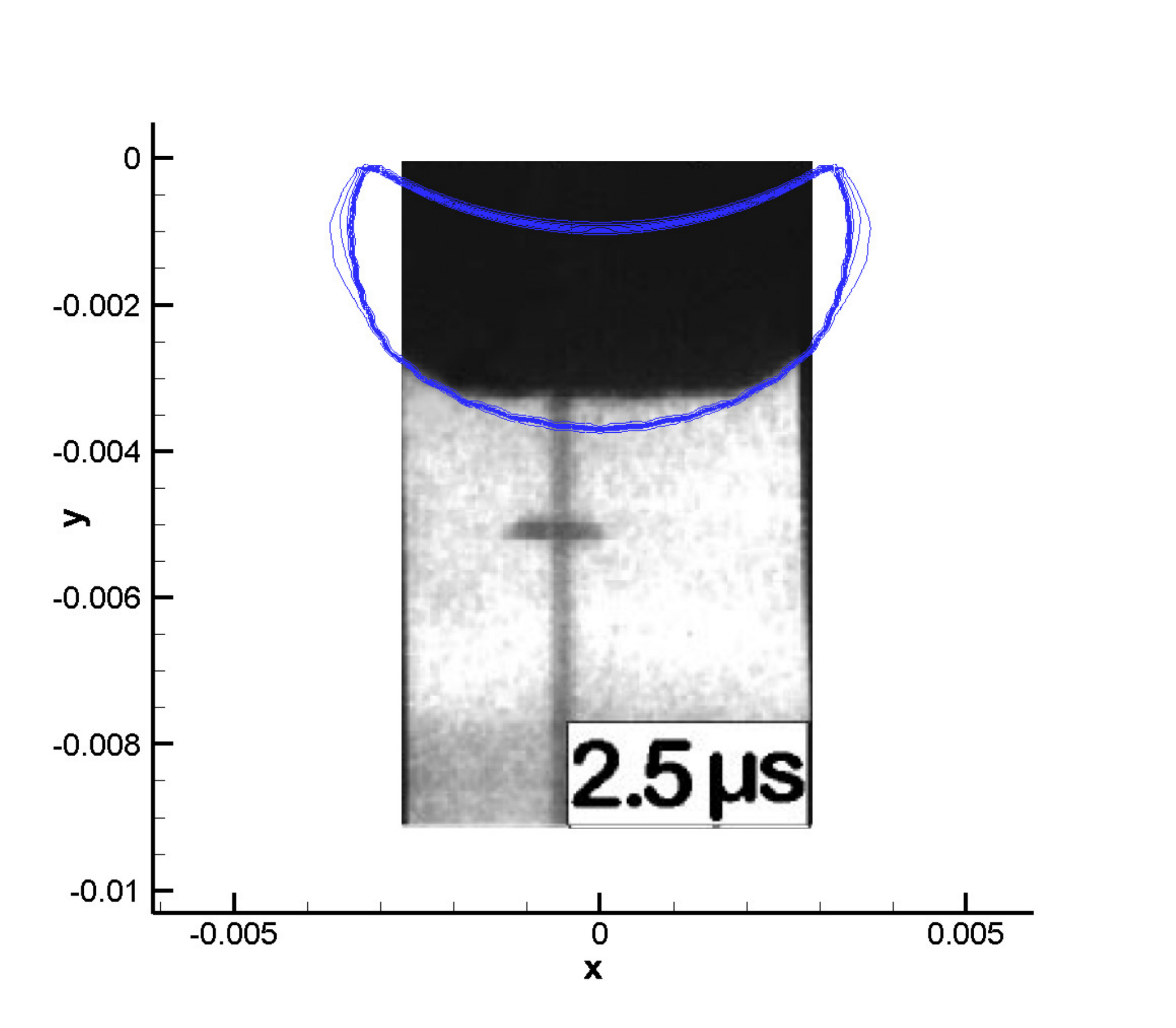}
        \includegraphics[width=0.32\textwidth]{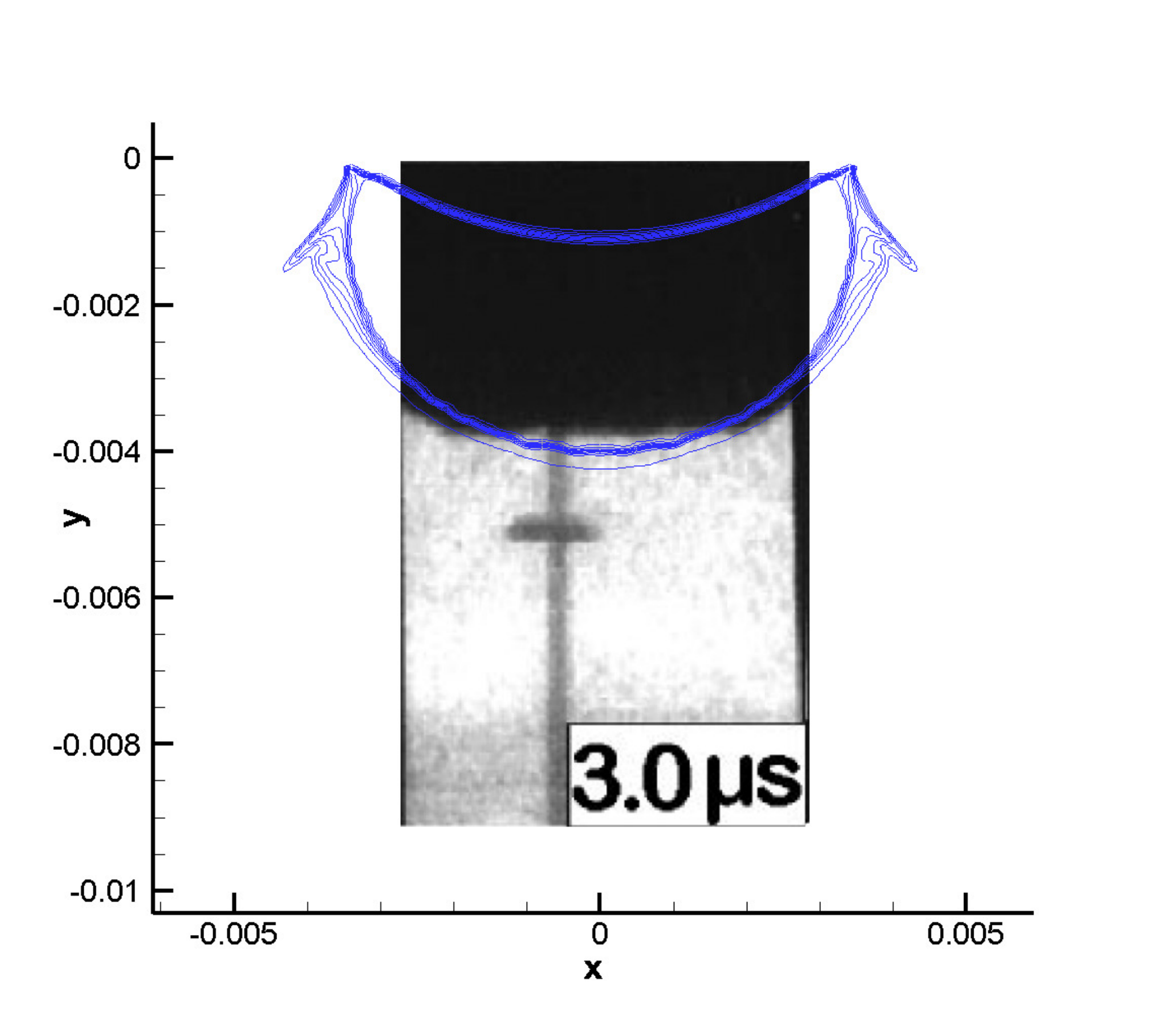}
        \includegraphics[width=0.32\textwidth]{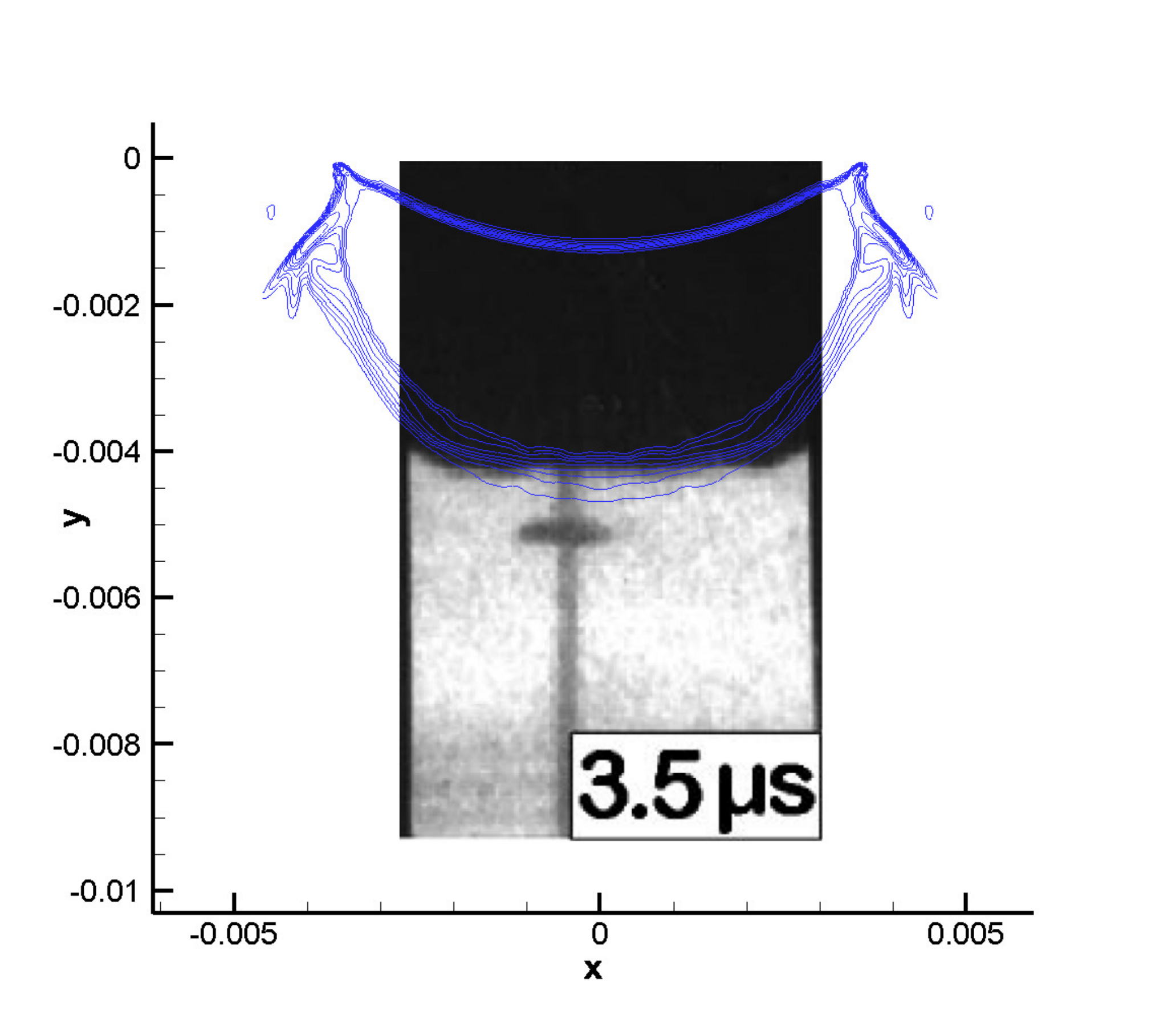}
    \end{center}
    \caption{Comparison between numerical fracture front (blue lines) with experimental results at several times.}
    \label{fig.exp2D_1}
\end{figure}

\begin{figure}[!bp]
    \begin{center}
        \includegraphics[width=0.32\textwidth]{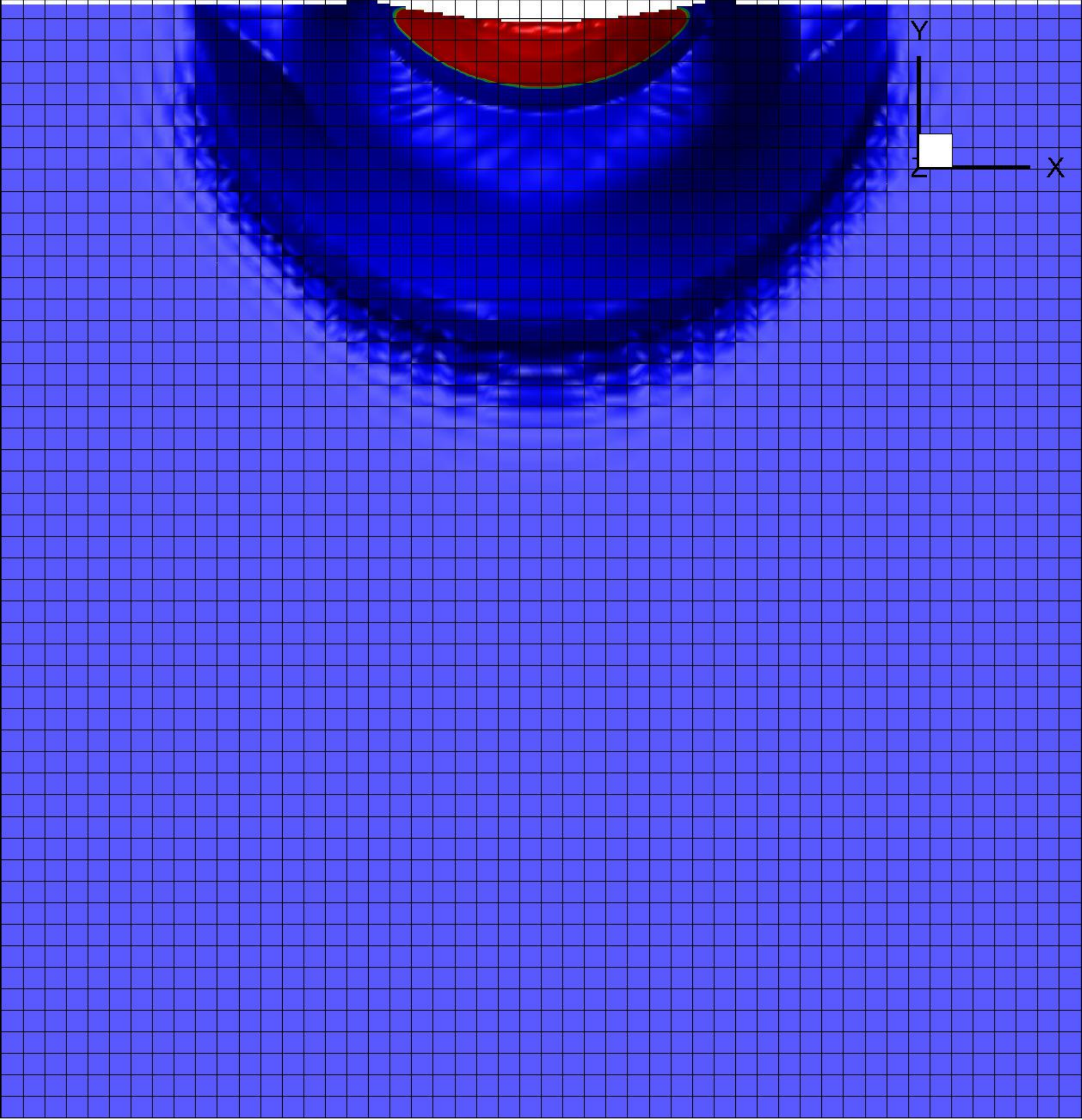}
        \includegraphics[width=0.32\textwidth]{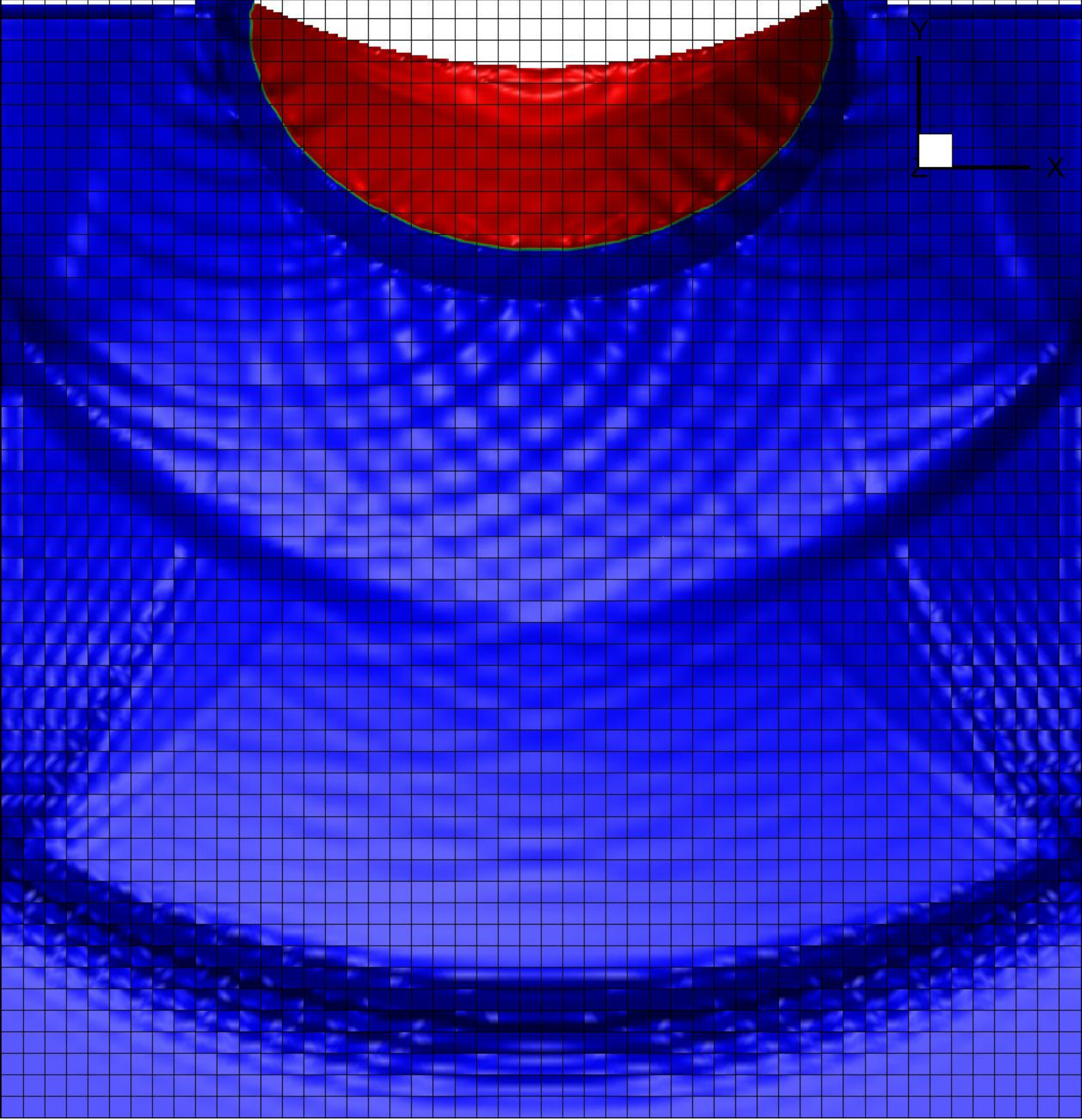}
        \includegraphics[width=0.32\textwidth]{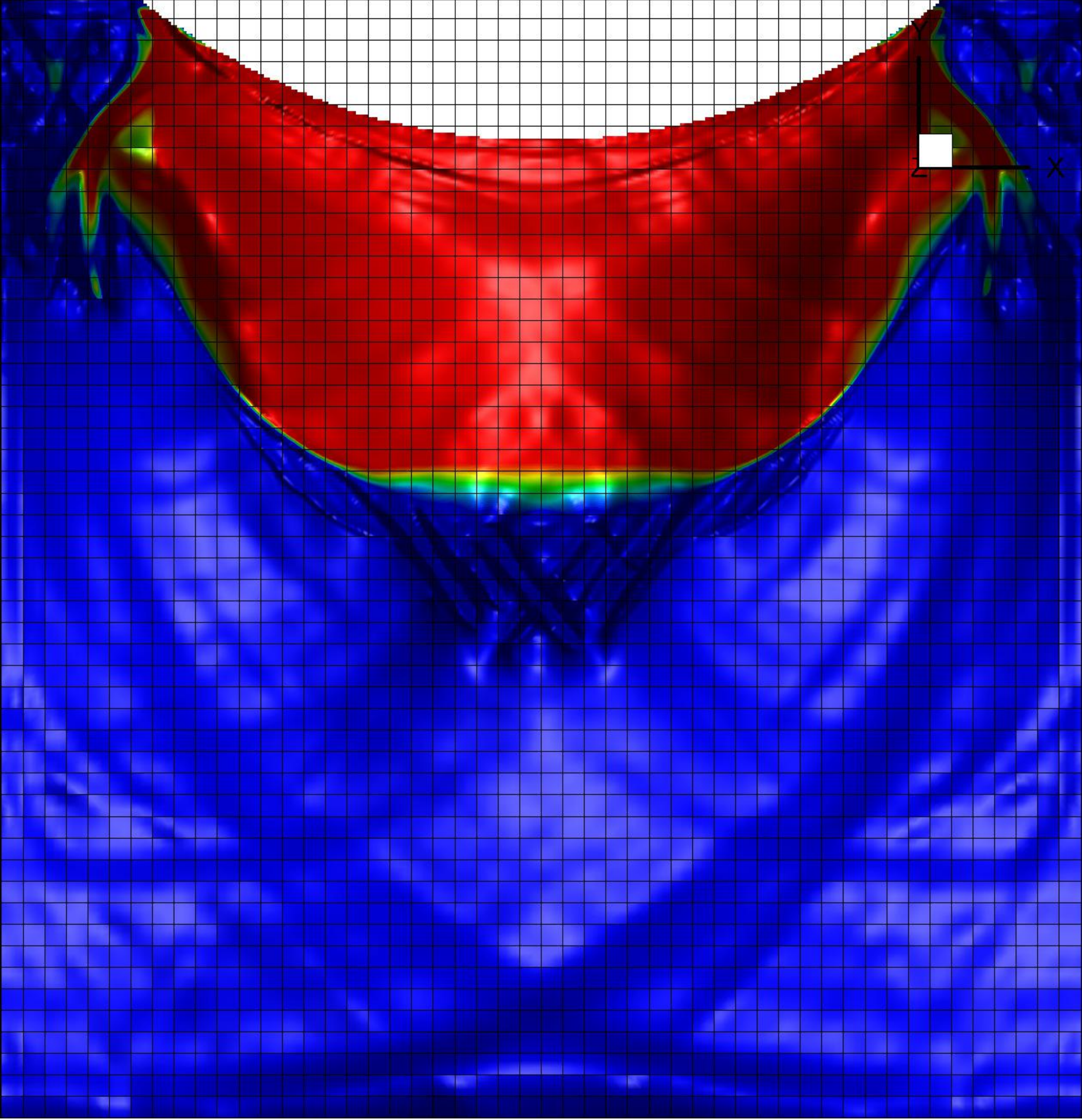}
        
    \end{center}
    \caption{Two dimensional crack propagation and wave pattern at times $t=0.5, 1.5$ and $3.8\,\up{\mu s}$. }
    \label{fig.exp2D_2}
\end{figure}

\begin{figure}[!bp]
    \begin{center}
        \includegraphics[width=0.32\textwidth]{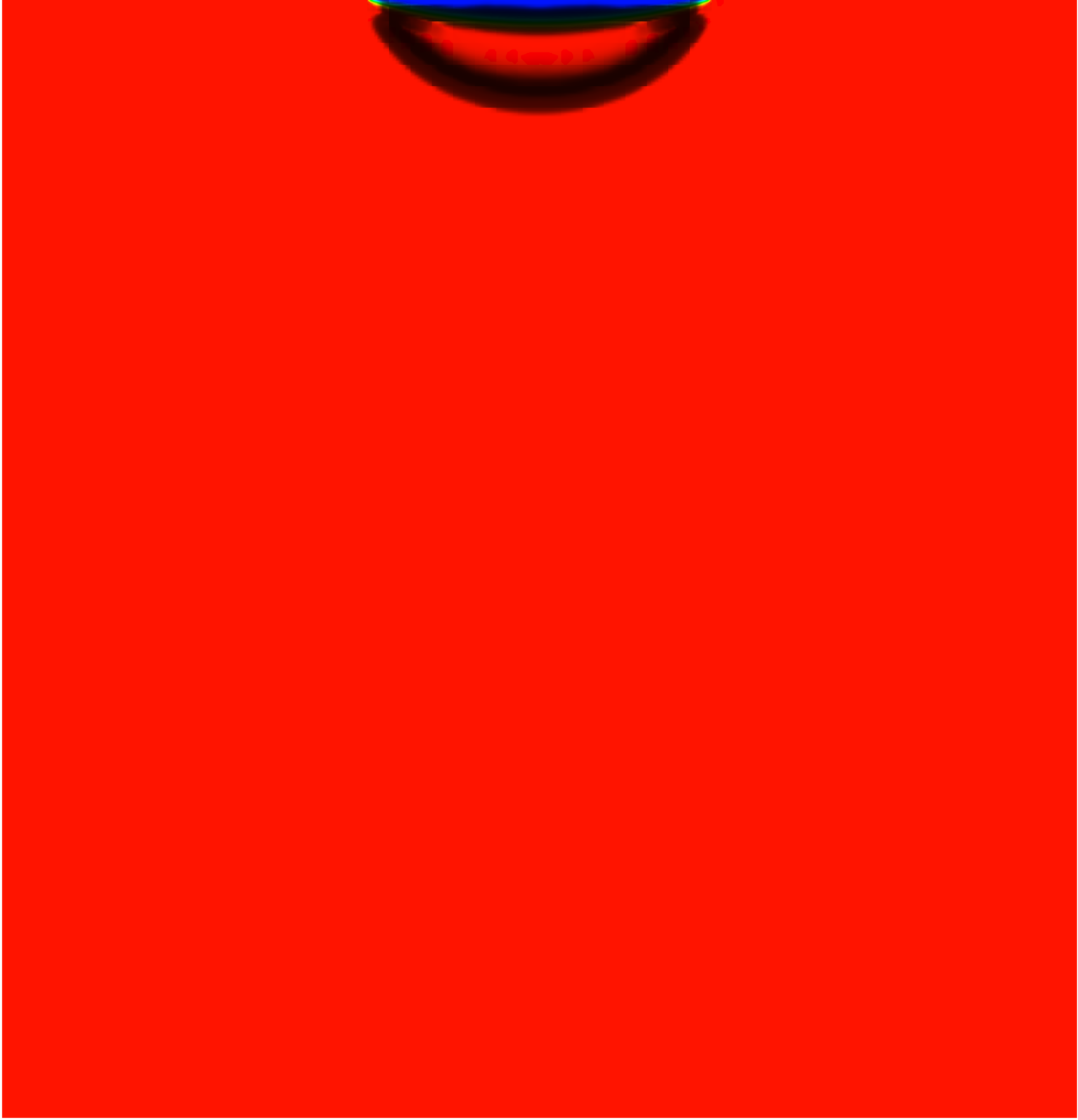}
        \includegraphics[width=0.32\textwidth]{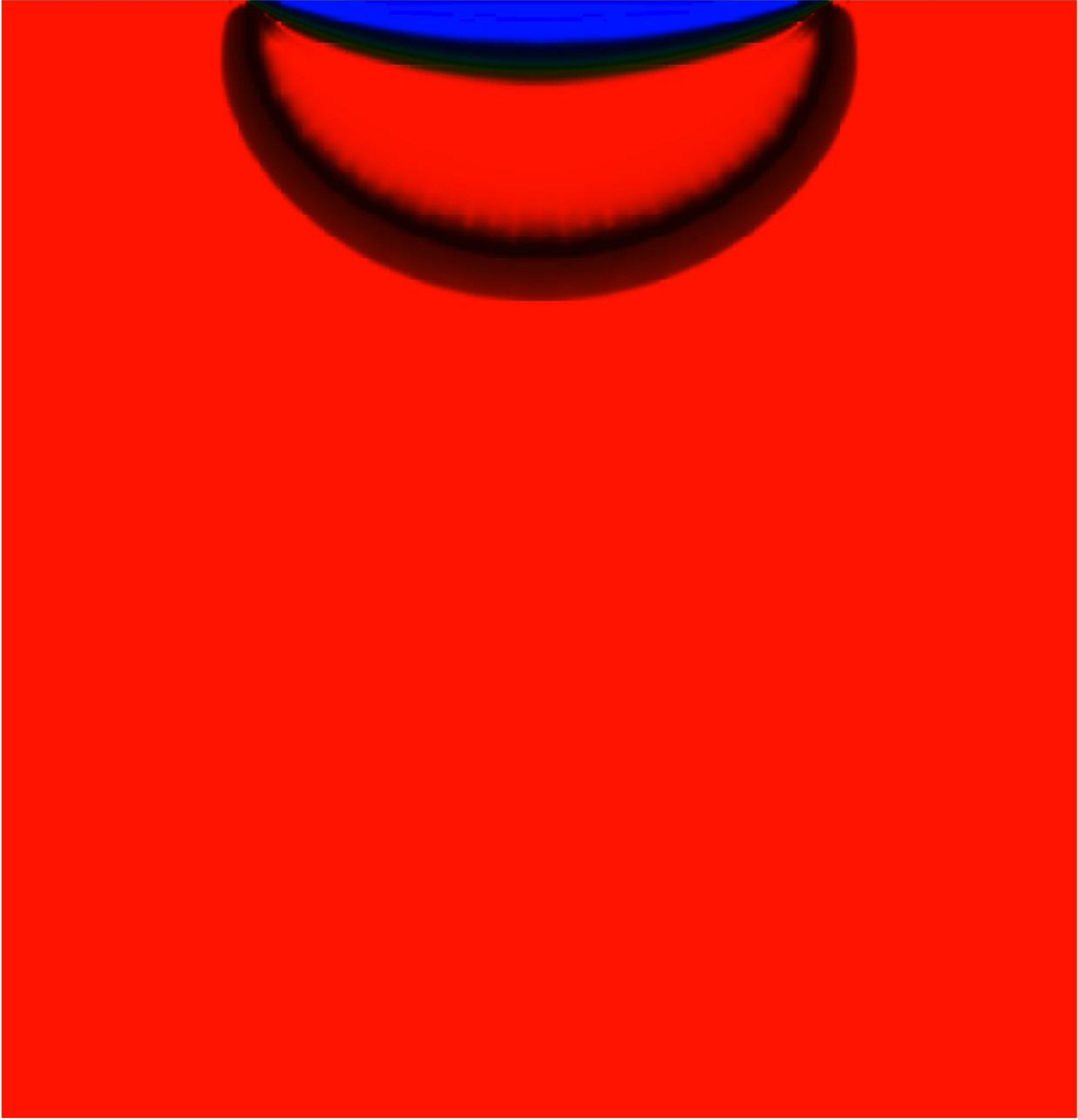}
        \includegraphics[width=0.32\textwidth]{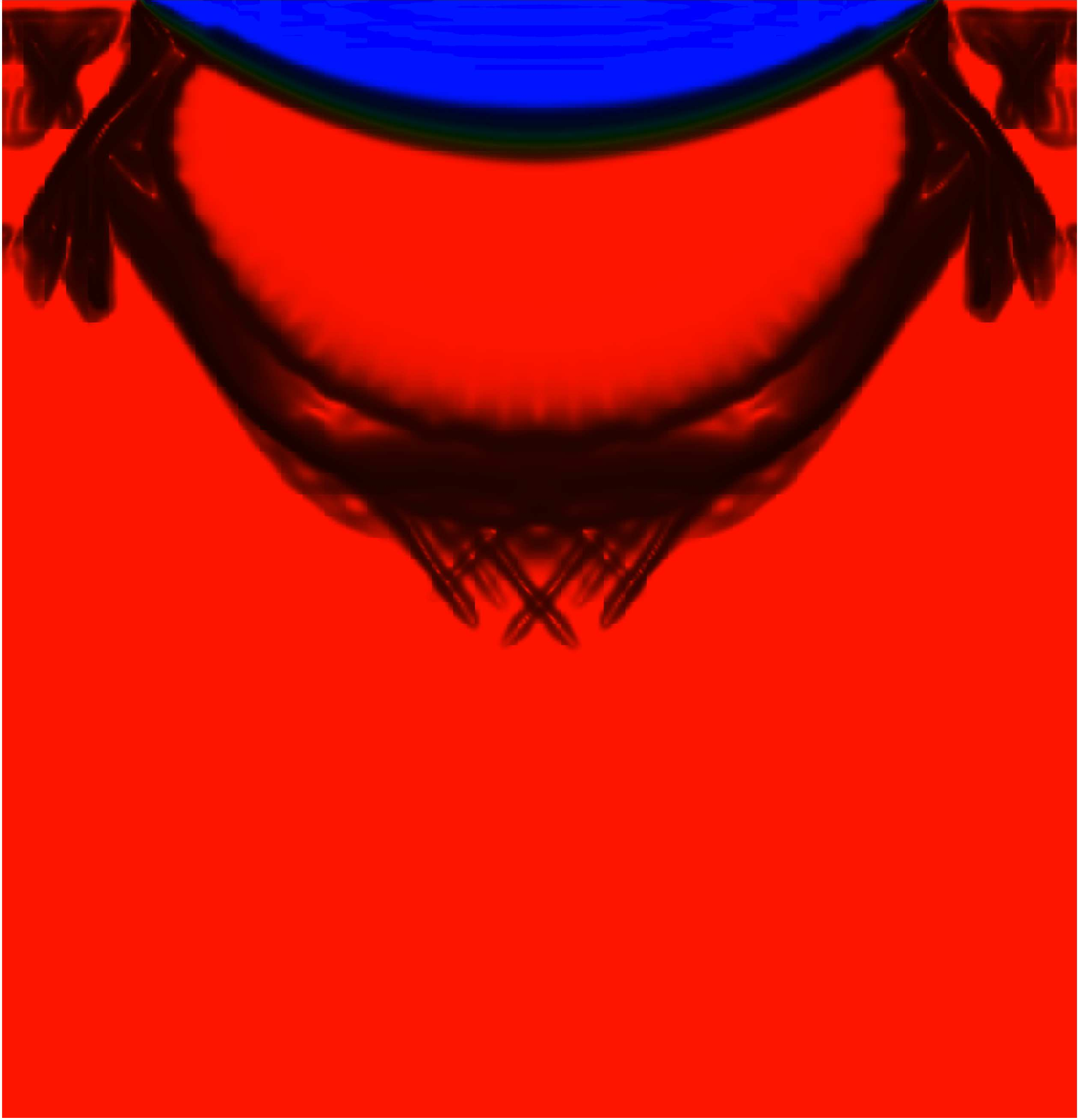}
    \end{center}
    \caption{Penetration pattern of the unbreakable copper (blue zone) in the Pyrex material (red zone) at times $t=0.5, 1.5$ and $3.8\,\up{\mu s}$. The evolution of the rupture front is also reported (shadow zone) and shows the transport of the rupture parameter according to the bullet penetration. }
    \label{fig.exp2D_3}
\end{figure}

where $\epsilon=10^{-3}\, Y_0$ is a smoothing factor and $s_0=1.8\,\up{GPa}$, $A=0.3$ and $B=0.33$. The computational 
domain is taken as $\Omega=[-0.005,0.005]\times [-0.0104,0.0004]$  covered with $50 \times 54$ elements 
 and using a polynomial approximation degree of $N=M=2$. The hemispherical projectile is represented 
 through time dependent boundary conditions. The impact velocity is set as $v=536\,\up{m/s}$, while the penetration velocity, that here changes the impact area, is extrapolated from \cite{Resnyansky2003} where also the bullet dynamics was considered. The resulting  evolution of the damaged front is compared with the experimental pictures taken by Bourne in 1997, see  \cite{Bourne1997}, and is reported in Figure $\ref{fig.exp2D_1}$. The damage front propagation is close to the one obtained in the experiments and is slower than the elastic precursor wave as shown in Figure $\ref{fig.exp2D_2}$. The evolution of the \textcolor{black}{damage variable} $\xi$ is presented in Figure $\ref{fig.exp2D_2}$, showing the evolution of the fracture wave in the Pyrex block. Furthermore, secondary cracks can be observed in Figure $\ref{fig.exp2D_1}$, similarly to the one obtained in \cite{Resnyansky2003}. 
The time evolution of the bullet penetration is reported in Fig. $\ref{fig.exp2D_3}$ where the rupture 
front is also depicted. Note that the \textcolor{black}{damage variable} $\xi$ is transported during the bullet impact (blue area) and no rupture in the copper is observed due to the strongest material properties.
We can finally observe how for long times the shape of the rupture is different in the numerical 
case and is more close to the case of a flat-nose rod, see again  \cite{Resnyansky2003}. Nevertheless, 
the experimentally observed propagation velocity of the rupture front is very well reproduced in our numerical simulations. 


\begin{table}[!bp]
    \caption{Set of parameters for the materials.}
    \centering
        \begin{tabular}{llllllll}
            \hline
                                         & Rock 1    & Rock 2    & Rock 3    & Rock 4    & Pyrex                 & Copper    \\
            \hline
            $\rho\,[\up{kg\,m^{-3}}]$    & 2670      & 2620      & 2670      & 2670      & 2230                  & 8930      \\
            $\mu_I\,[\up{GPa}]$          & 32.04     & 21.44     & 32.04     & 32.04     & 30.36                 & 48.27     \\
            $\mu_D\,[\up{GPa}]$          & 27.46     & 0.15008   & 27.46     & $0.03204$ & $0.1518$              & 41.38     \\
            $\lambda_I\,[\up{GPa}]$      & 32.04     & 21.44     & 32.04     & 32.04     & 20.90                 & 105.79    \\
            $\lambda_D\,[\up{GPa}]$      & 35.10     & 0.15008   & 35.10     & 53.38     & 30.97                 & 110.39    \\
            % $c_s\,[\up{m\,s^{-1}}]$    & 3464      &           & 3464      & 3464      & 3690                  & 2325      \\
            % $c_l\,[\up{m\,s^{-1}}]$    & 6000      &           & 6000      & 6000      & 6050                  & 4760      \\
            $\theta_0\,[-]$              & 10.0      & 1.0       & 0.2       & 1.0       & 1.0                   & 0.0       \\
            $Y_0\,[\up{MPa}]$            & $180$     & 10        & $240$     & $9.0$     & $1200$                & $10^{16}$ \\
            $Y_1\,[\up{MPa}]$            & $10.0$    & $10^{-6}$ & $10.0$    & $10.0$    & $10.0$                & $10^{16}$ \\
            $a\,[-]$                     & 42.5      & 60.0      & 42.5      & 52.5      & 32.5                  & 1.0       \\
            $\alpha_I\,[-]$              & 36.25     & 0.0       & 36.25     & 36.25     & 36.25                 & 40.0      \\
            $\alpha_D\,[-]$              & 36.25     & 0.0       & 36.25     & 36.25     & 34.8                  & 40.0      \\
            $\beta_I\,[\up{Pa^{-1}}]$    & 0.0       & 0.0       & 0.0       & 0.0       & $22.31\times10^{-9}$  & 0.0       \\
            $\beta_D\,[\up{Pa^{-1}}]$    & $10^{-6}$ & 0.0       & $10^{-6}$ & $10^{-6}$ & $223.07\times10^{-9}$ & 0.0       \\
            % $\mu_{1\rightarrow 2}$     & 6/7       &           & 6/7       & 0.001     & 0.005                 & 6/7       \\
            % $\lambda_{1\rightarrow 2}$ & 2/3       &           & 2/3       & 2/3       & 1/3                   & 2/3       \\
            $\tau_{I0}\,[\up{s}]$        & $3.2\times10^6$ & $10^{ 5}$ & $3.2\times10^6$ & $3.2\times10^6$ & $3.0\times10^6$  & $4.8\times10^6$  \\    
            $\tau_{D0}\,[\up{s}]$        & $2.75\times10^6$ & $10^{-3}$ & $2.75\times10^6$ & $3.2\times10^3$ & $1.5\times10^4$  &  $4.1\times10^6$ \\    
            \hline
        \end{tabular}
    \label{tab.materials}
\end{table}

\clearpage 

\subsection{Fracture generation in a piece of pre-damaged rock}  
\label{sec.rock}

\begin{figure}[!bp]
    \begin{center}
        \includegraphics[width=0.6\textwidth]{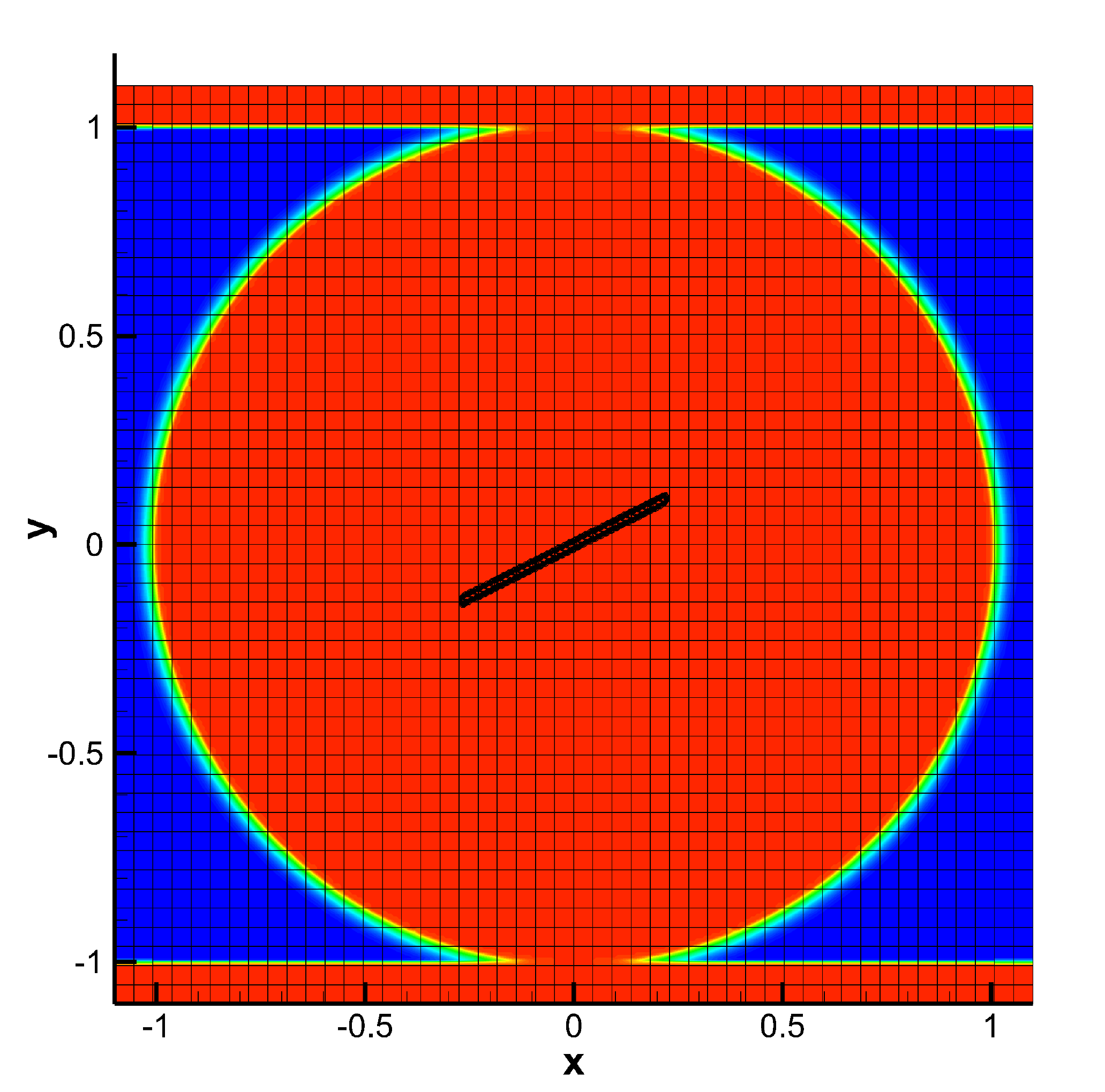}
    \end{center}
    \caption{Representation of the pre-damaged disc in the diffuse interface approach via the solid volume fraction. The darkest area shows the pre-damaged zone.} 
    \label{Disc_alpha}
\end{figure}

In order to show the applicability of our new algorithm to realistic experiments that are frequently 
carried out in civil engineering and in geophysics, we consider the crack propagation in a pre-damaged 
rock-like disc (so-called Brazilian test)\textcolor{black}{, with the aim of comparing to the experimetal results given
in \cite{DiskExp2014}. The disc is described with the moving diffuse interface 
approach, and its boundary is thus represented by means of the solid volume fraction $\alpha$. 
An inclined slit is represented as predamaged zone. In order to match the experimental data, we rotated
the disc with respect to the clamping apparatus in such a way that the slit has an inclination angle
of $35^\circ$ with respect to the $x$-axis.
We impose a velocity pointing towards the disc 
on the upper and bottom boundary, with a magnitude given by the gaussian profile
$|v| = 4\,\exp(-(25\,x)^2)$.
We take our domain to be $\Omega=[-1.1,1.1]^2$ and employ a fourth order ADER-DG 
scheme with \textit{a-posteriori} subcell limiting. 
The distribution of the volume fraction $\alpha$ of 
the diffuse interface method is depicted in Fig. \ref{Disc_alpha}. 
The material properties for the disc are those associated to Rock 2 in Table \ref{tab.materials}, 
while in the top and bottom clamps the material is modified by setting $Y_0 = Y_1 = 100\,\up{TPa}$, 
effectively rendering these regions unbreakable.
For the equivalent stress we take a Drucker--Prager-like combination of shear stress and 
pressure defined as}
\begin{equation}
% Y=A |Y_s|+B \frac{1}{2}|Y_p|\left[1+\textnormal{erf}\left(-\frac{Y_p}{\epsilon}\right)\right]
    \textcolor{black}{Y = A\,|Y_s| + B\,Y_p + C}
\end{equation}

\begin{figure}[!bp]
    \begin{center}
        \includegraphics[width=0.99\textwidth]{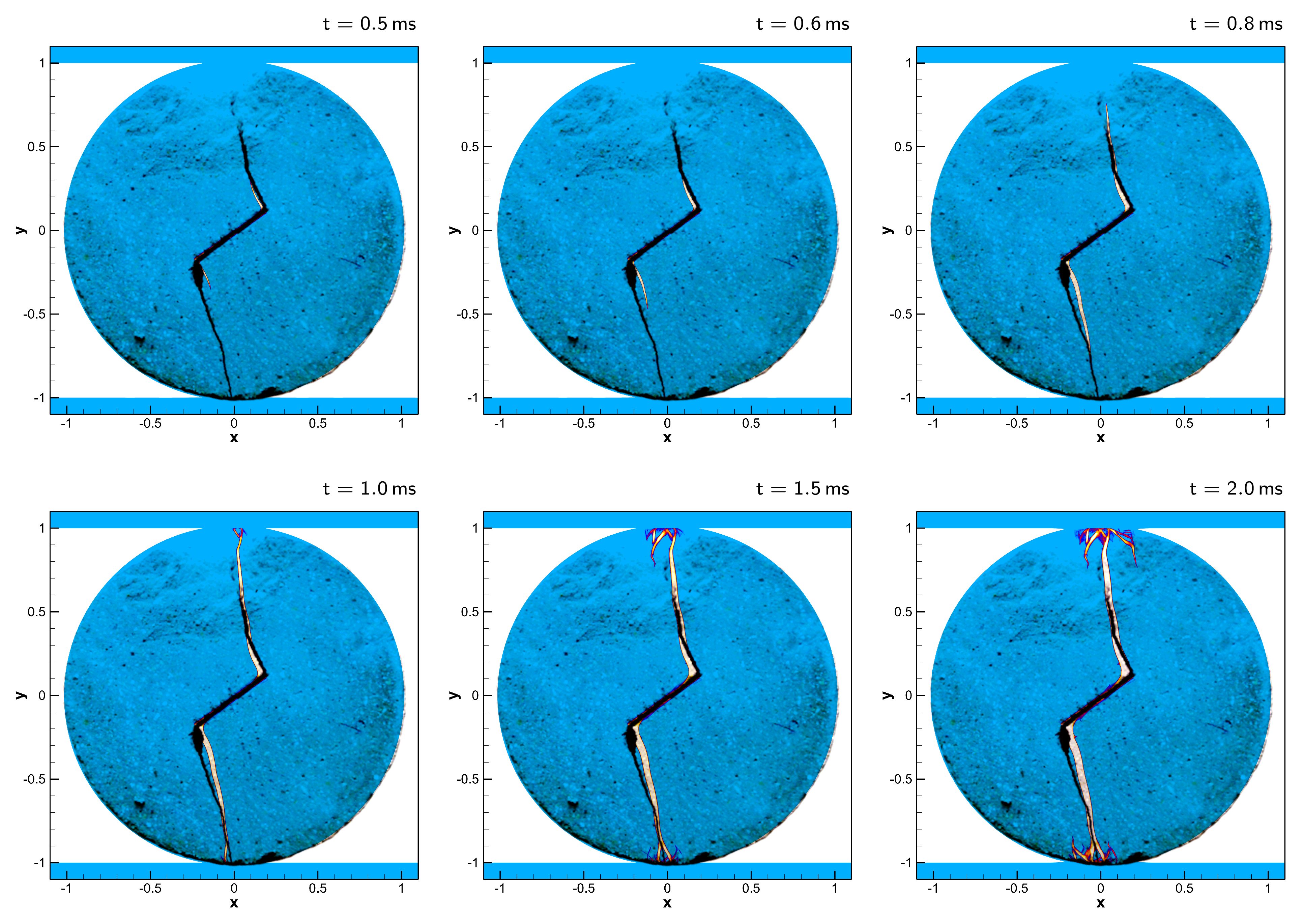}
    \end{center}
    \caption{Evolution of the rupture front (color contours of the damage field $\xi$) compared with experiments 
    at times $t=[0.5,0.6,0.8,1.0,1.5,2.0]\,\up{ms}$. The computations are carried out with a fourth order
    ADER-DG scheme with \textit{a-posteriori} subcell limiting on a uniform Cartesian grid of $256 \times 256$ elements. 
    Only regions with $\alpha > 0.5$ are shown. 
    The experimental picture is taken from \cite{DiskExp2014}.}
    \label{Disc2}
\end{figure}

\begin{figure}[!bp]
    \begin{center}
        \includegraphics[width=0.99\textwidth]{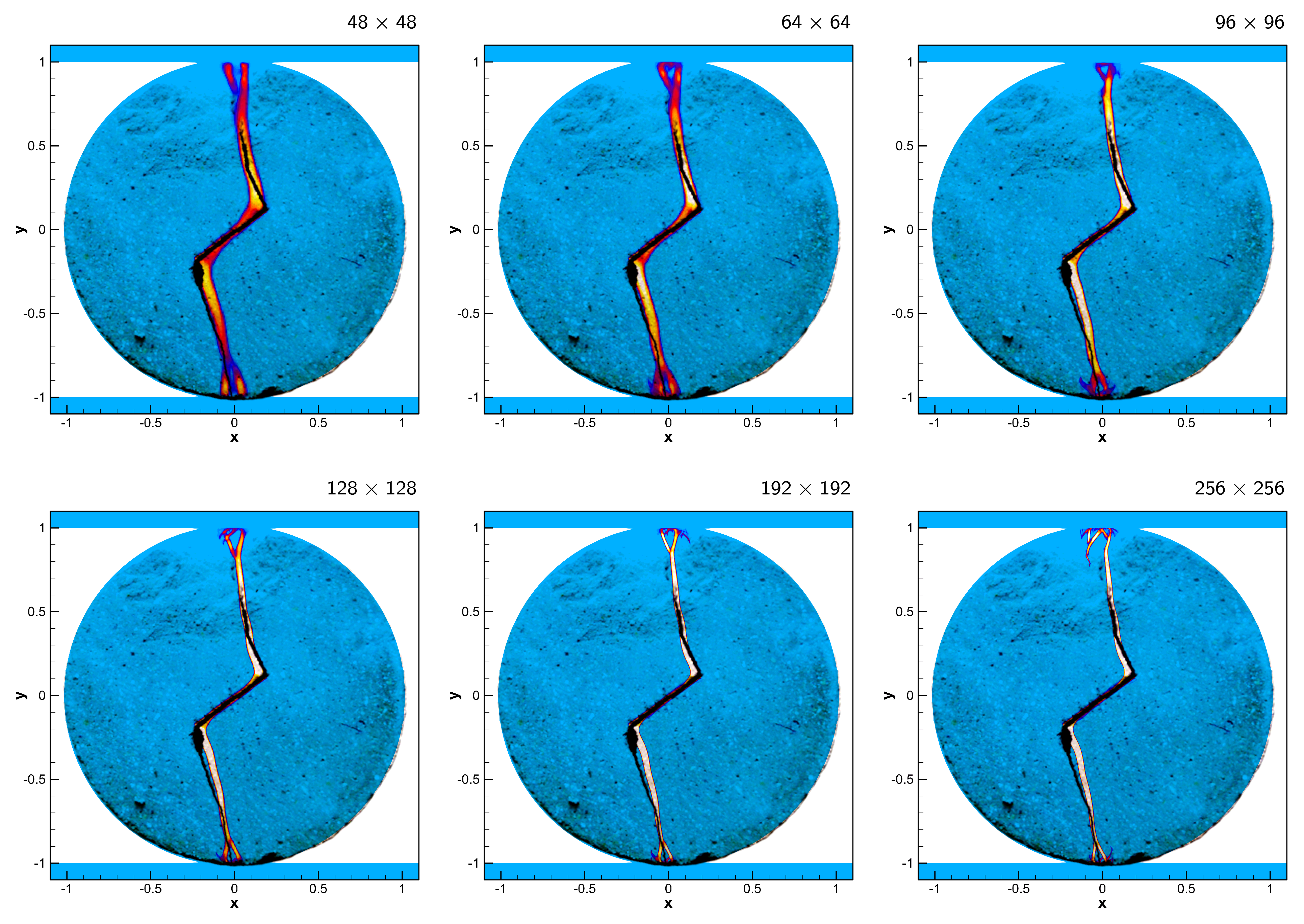}
    \end{center}
    \caption{Mesh convergence study (color contours of the damage field $\xi$) for the simulation 
    of a pre-damaged rock disc. The computations employ a fourth order ADER-DG scheme with MUSCL-Hancock
    \textit{a posteriori} subcell finite volume limiting and show mesh convergence of the main cracks for
    six grid spacings ranging from $48 \times 48$ to $256 \times 256$ Cartesian cells. Only regions with $\alpha > 0.5$ are shown.
    The experimental picture has been taken from \cite{DiskExp2014}.}
    \label{Disc3}
\end{figure}

\textcolor{black}{where $A=1.0$, $B=1.5$ and $C = -2.0\,\up{MPa}$}
The resulting evolution of the rupture is shown in Fig. \ref{Disc2}. The main crack propagates 
starting from the corners of the pre-damaged zone, following the experimental results very 
closely to the top and bottom of the disc. In agreement with the experimental observations, 
we have also a damaged zone due to shear on the upper and lower part of the disc, while the rupture 
inside is driven by traction. \textcolor{black}{Last but not least, we add a careful mesh 
convergence analysis for this nontrivial test problem, in order to check whether the position 
and shape of the crack converges with mesh refinement. The result of this detailed mesh 
refinement study carried out on a sequence of uniform Cartesian grids with mesh resolution 
from $48 \times 48$ elements to $256 \times 256$ elements is shown in Figure \ref{Disc3}. One 
can observe that the main cracks properly converge with mesh refinement, in particular 
the crack becomes thinner on finer meshes. Only the secondary cracks generated at the upper 
and lower boundaries of the disc differ between one mesh and the other.}

%\clearpage 

\subsection{Out of plane rupture with secondary and tertiary crack generation}  
\label{sec.rupture2d}
\begin{figure}[!b]
    \begin{center}
        \includegraphics[width=0.9\textwidth]{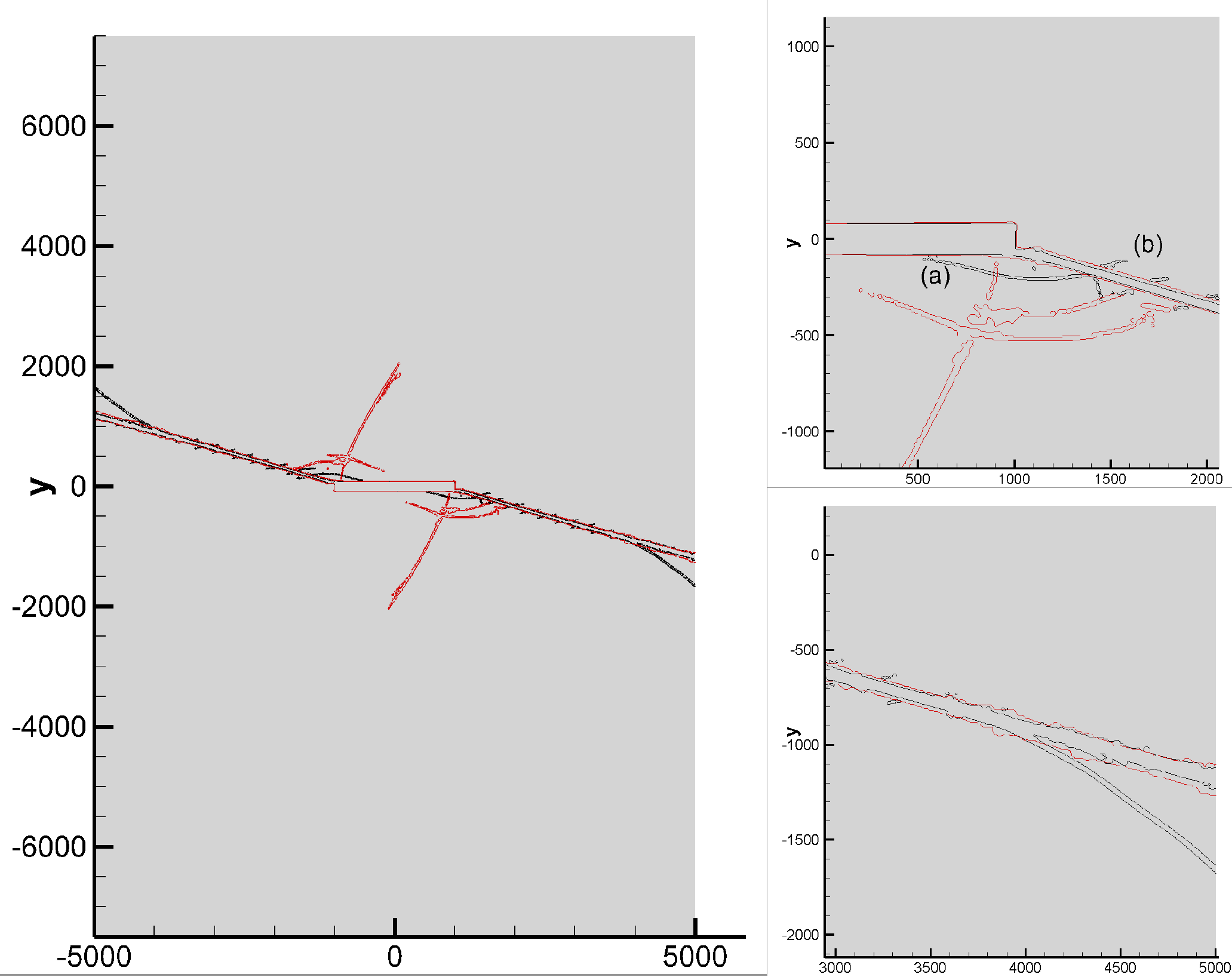}
    \end{center}
    \caption{Out of plane crack propagation. Comparison of the crack profile for a coarse 
    mesh (red line) and a finer one (dark line) at $t=5.0\,\up{s}$. The main crack 
    direction is the same, independently from the resolution, while small secondary 
    cracks are more present on the finer mesh (b). }
    \label{2DOPR}
\end{figure}

Here we consider the generation and propagation of a crack in a two dimensional domain $\Omega=[-5000,5000]\times[-7500,7500]$. The 
material is pre-damaged on the line $|x|<1000,|y|<66$, which is mesh aligned and that is loaded with
 an initial stress of  $\sigma_{yy}=-120\,\up{MPa}$ and $\sigma_{xy}=70\,\up{MPa}$. 
 Here we simply use the von Mises stress for the definition of the equivalent stress. 
The considered material is ROCK 1 (see Table \ref{tab.materials} for the defining parameters). 
In the numerical simulation we consider two AMR refinement levels that are activated 
dynamically according to the rupture propagation and that are needed to capture the 
small-scale structures that appear during the dynamic rupture process. 
We use a polynomial approximation degree of $N=M=3$ and on the coarsest level we consider 
two different meshes, a coarse one with $30 \times 24$ elements and a finer one with 
twice the resolution in each space dimension. The  subcell FV limiter is activated 
according to \cite{Dumbser2014}, as using an additional physical detection criterion 
based on the ratio of the von Mises stress and the yield stress of the material.  
A comparison between the crack propagation using the coarse and the fine mesh at time 
$t=5.0\,\up{s}$ is reported in Figure \ref{2DOPR}. Here we can observe how the direction 
of the main crack is the same for both grid resolutions, while with higher resolution we 
generate also several secondary and tertiary cracks that bifurcate from the main one. 
This is the expected behaviour due to the better resolution of shear bands forming. 
%Alice: add? (mesh-dependent spatial strain/strain rate accumulation).
In this test one secondary crack prevents energy accumulation close to the corner points, 
see top right panel in Figure \ref{2DOPR}.

\begin{figure}[!b]
    \begin{center}
        \includegraphics[width=0.32\textwidth]{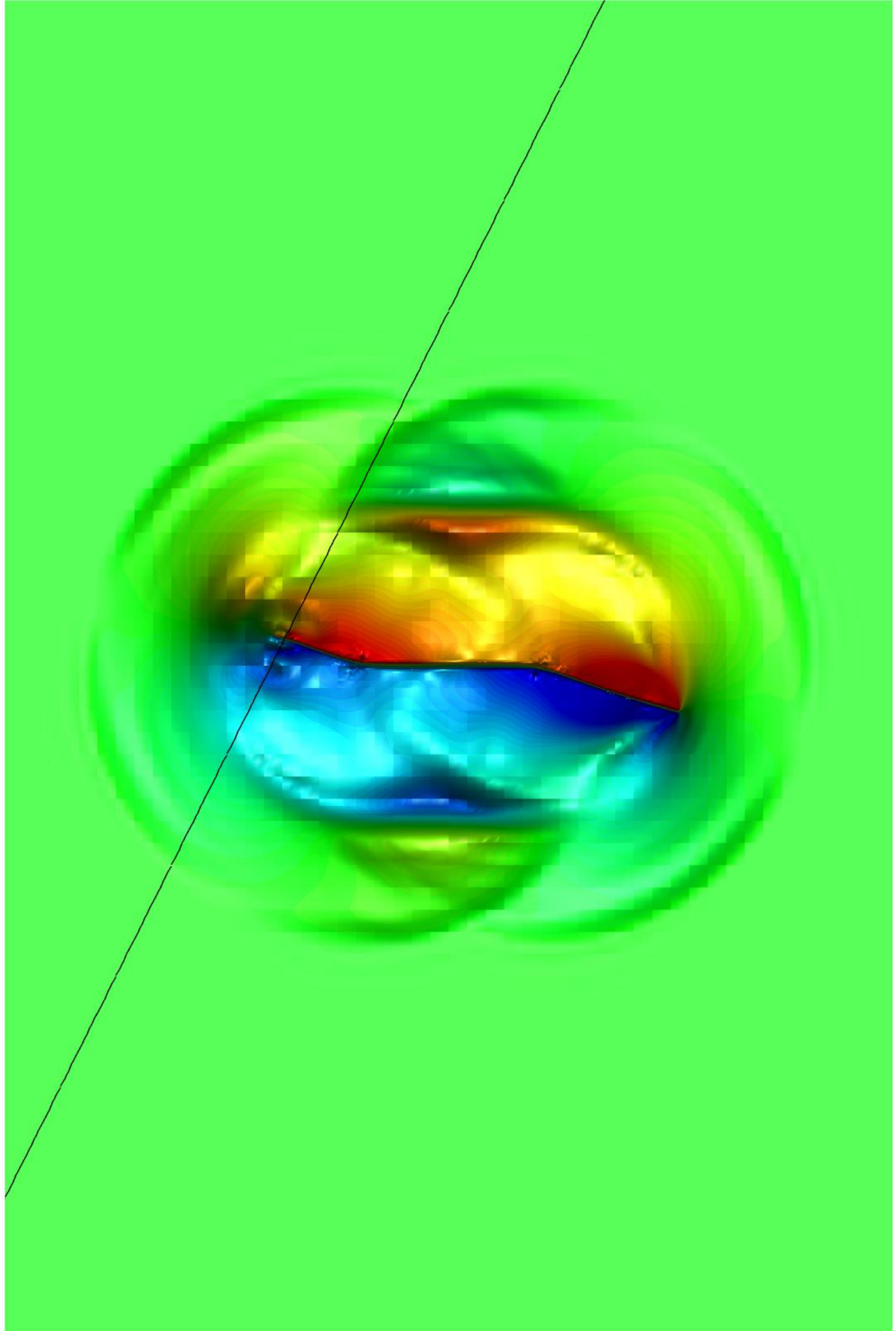}
        \includegraphics[width=0.32\textwidth]{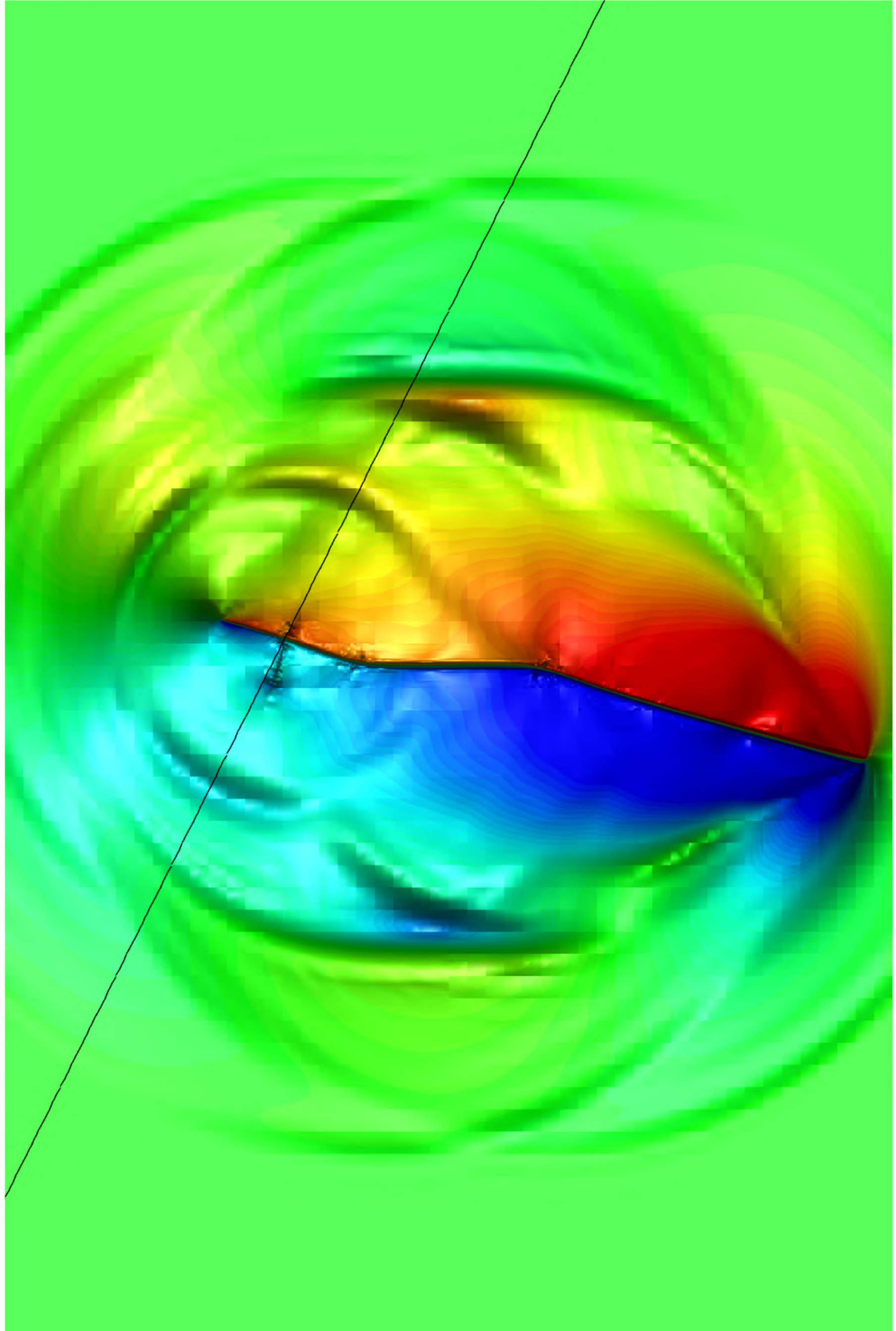}
        \includegraphics[width=0.32\textwidth]{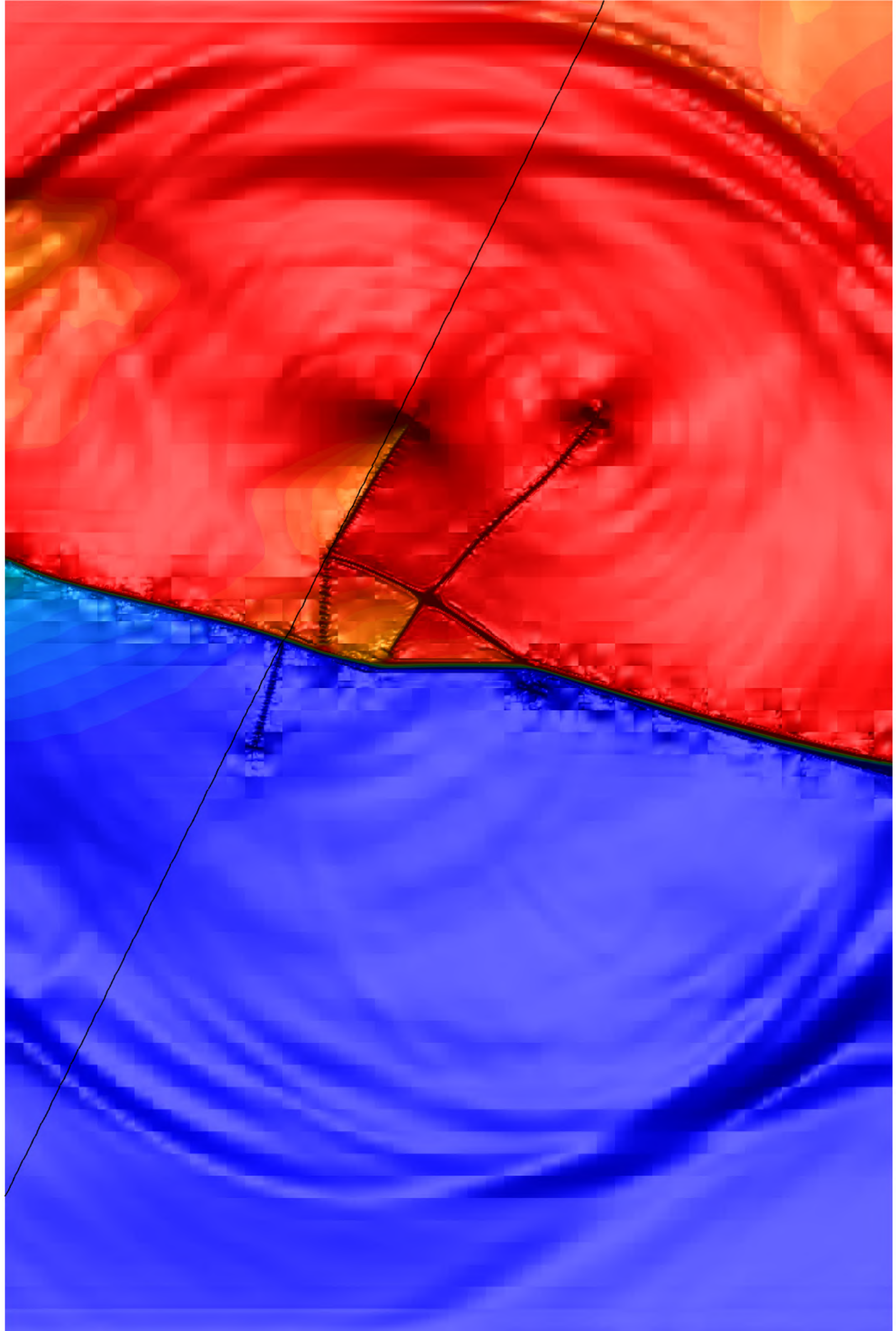}
    \end{center}
    \caption{Evolution of the out of plane rupture in a heterogeneous material 
    at times $t=[0.5,0.9,5.0]$. The separation line is reported and the two material are characterized by different 
    critical stress $Y_0$, i.e. $Y_0=180\,\up{MPa}$ and $Y_0=240\,\up{MPa}$ respectively 
    for the first and the second material. Super shear can be observed in the right crack 
    propagation while crack tends naturally to propagate parallel to the separation line of the two materials  }
    \label{2DOPR_ETMAT}
\end{figure}

\begin{figure}[!b]
    \begin{center}
        \includegraphics[width=0.32\textwidth]{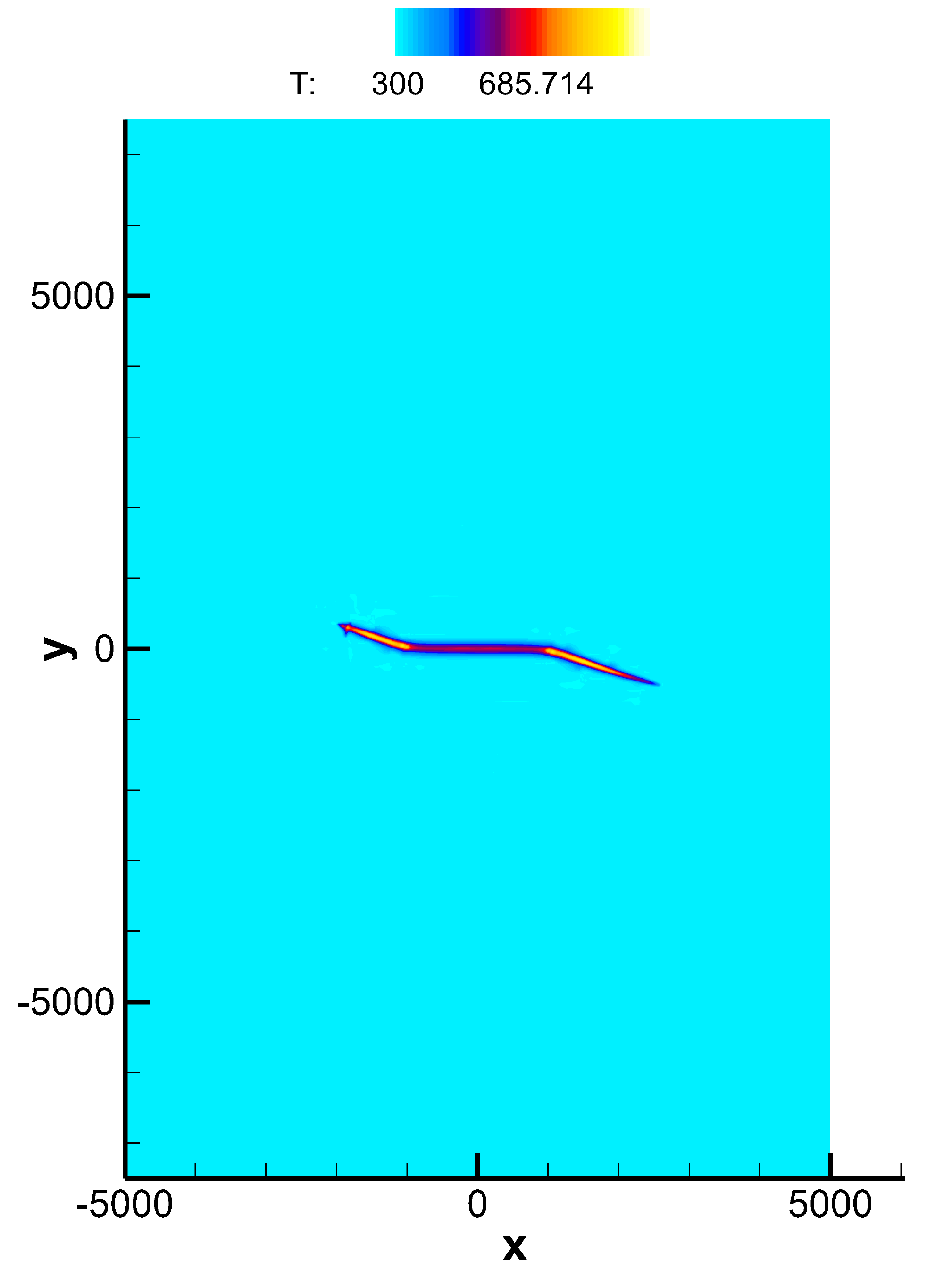}
        \includegraphics[width=0.32\textwidth]{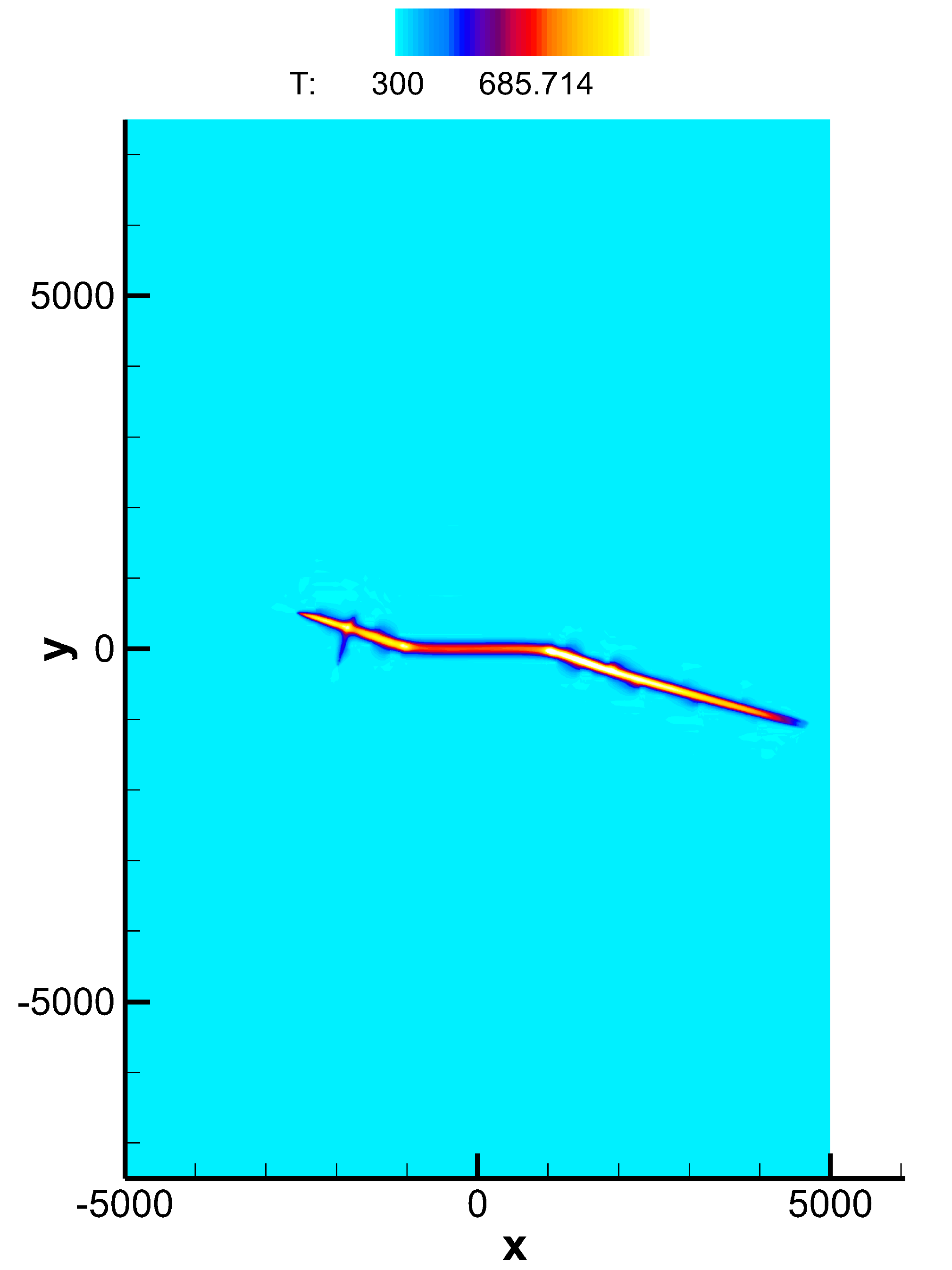}
        \includegraphics[width=0.32\textwidth]{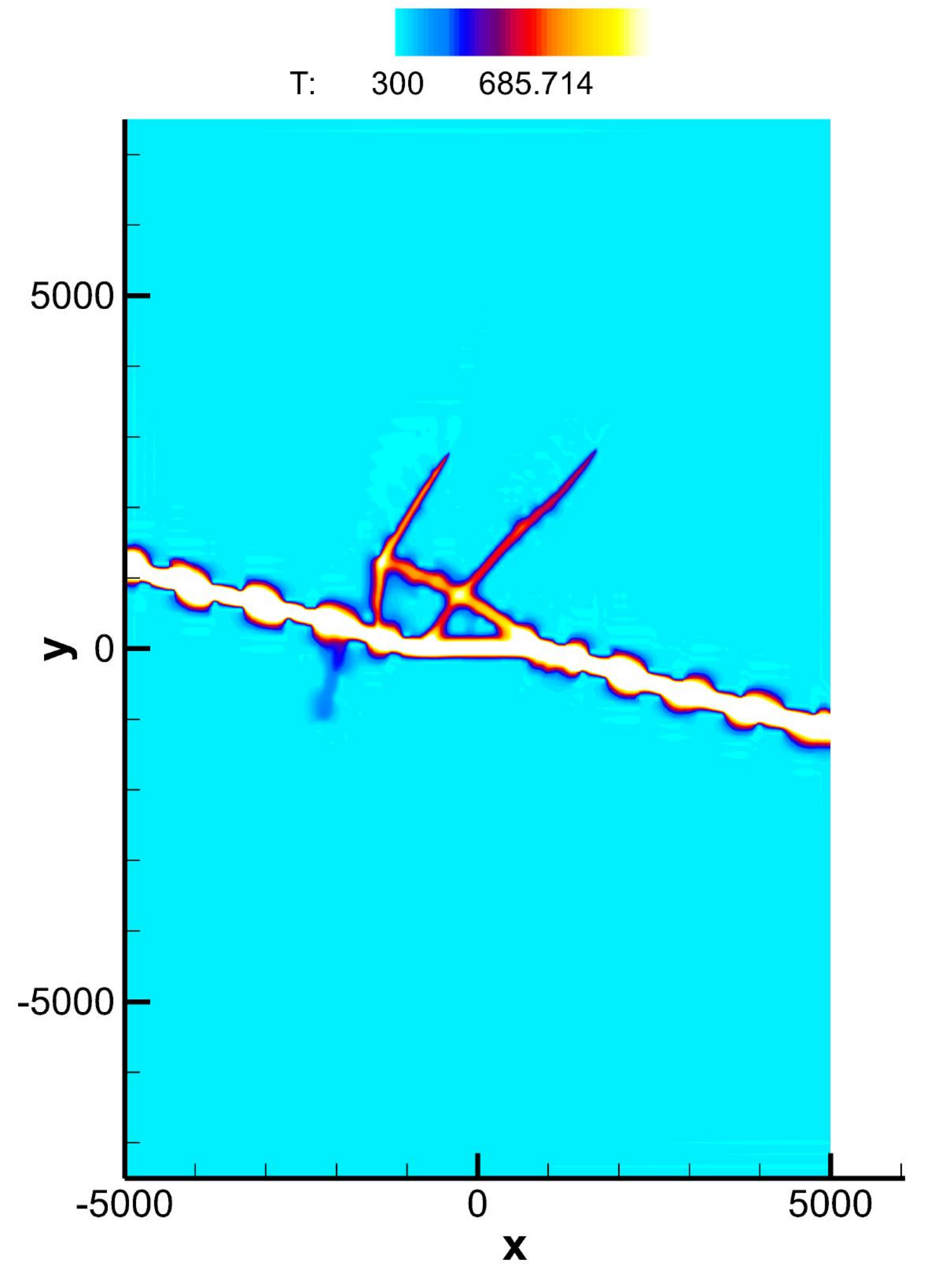} \\
        \includegraphics[width=0.9\textwidth]{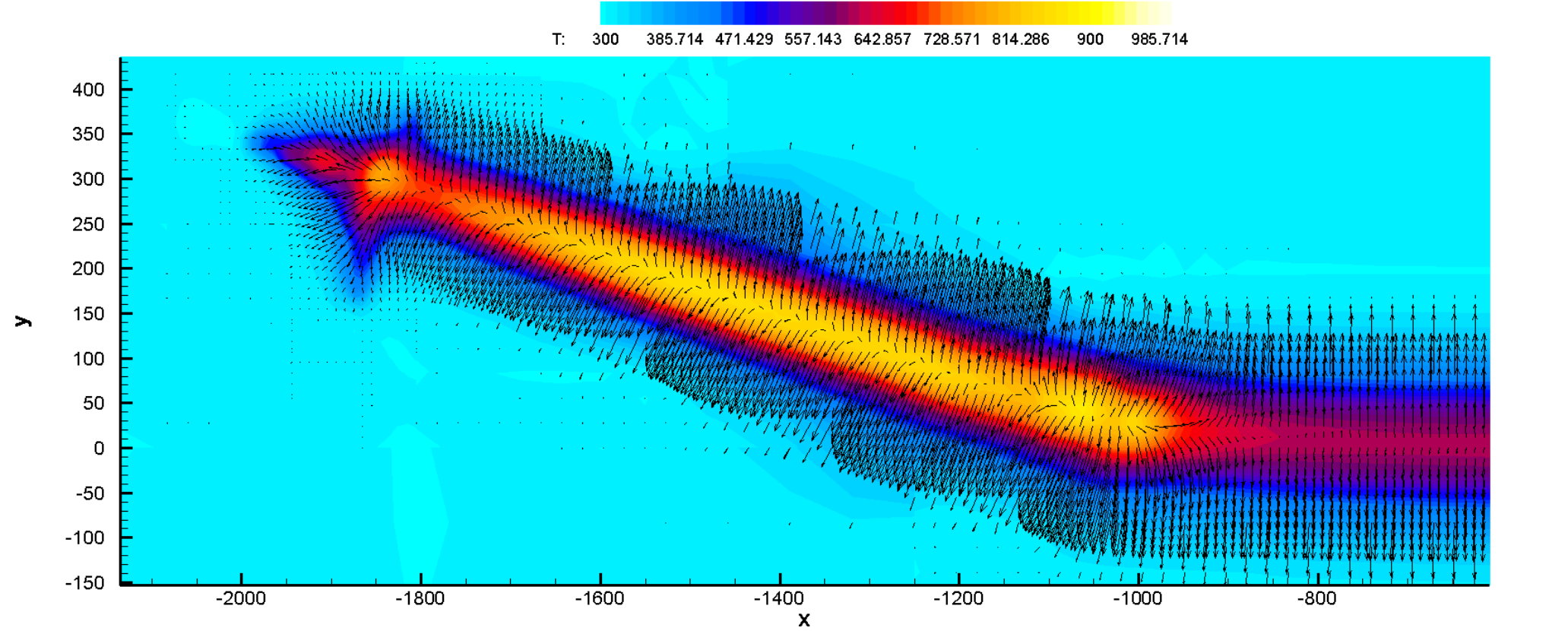}
    \end{center}
    \caption{Evolution of the thermal trace generated by friction during the sliding 
    process (top figures); representation of the heat flux vector $(J_x,J_z)$ at $t=0.9$ in the 
    bottom figure, showing the heat propagation during the process.}
    \label{2DOPR_ETMAT_2}
\end{figure}

In order to show the effect of heterogeneous materials in the crack propagation we consider 
the same domain using $64 \times 100$  elements with one AMR refinement level 
and two different materials, namely ROCK 1 and ROCK 3, that are characterized by a different 
yield stress. The resulting time evolution of crack lines at times $t=[0.5,0.9,5.0]$ are 
reported in Figure \ref{2DOPR_ETMAT}, where also the horizontal velocity component is depicted. 
We can note how the presence of the second material, whose interface is represented by the black
 line in  \ref{2DOPR_ETMAT}, breaks the symmetry of the rupture propagation that is decelerated in 
 the ROCK 3. Furthermore, the secondary cracks seem to naturally align with the interface between 
 the two materials. Finally, the thermal trace is reported in $\ref{2DOPR_ETMAT_2}$ and shows 
 the energy that is converted into thermal energy due to friction and localized on the sliding lines.

\subsection{Double coin crack in three space dimensions} 
\label{sec.doublecoin} 
\begin{figure}[!b]
    \begin{center}
        \includegraphics[width=0.45\textwidth]{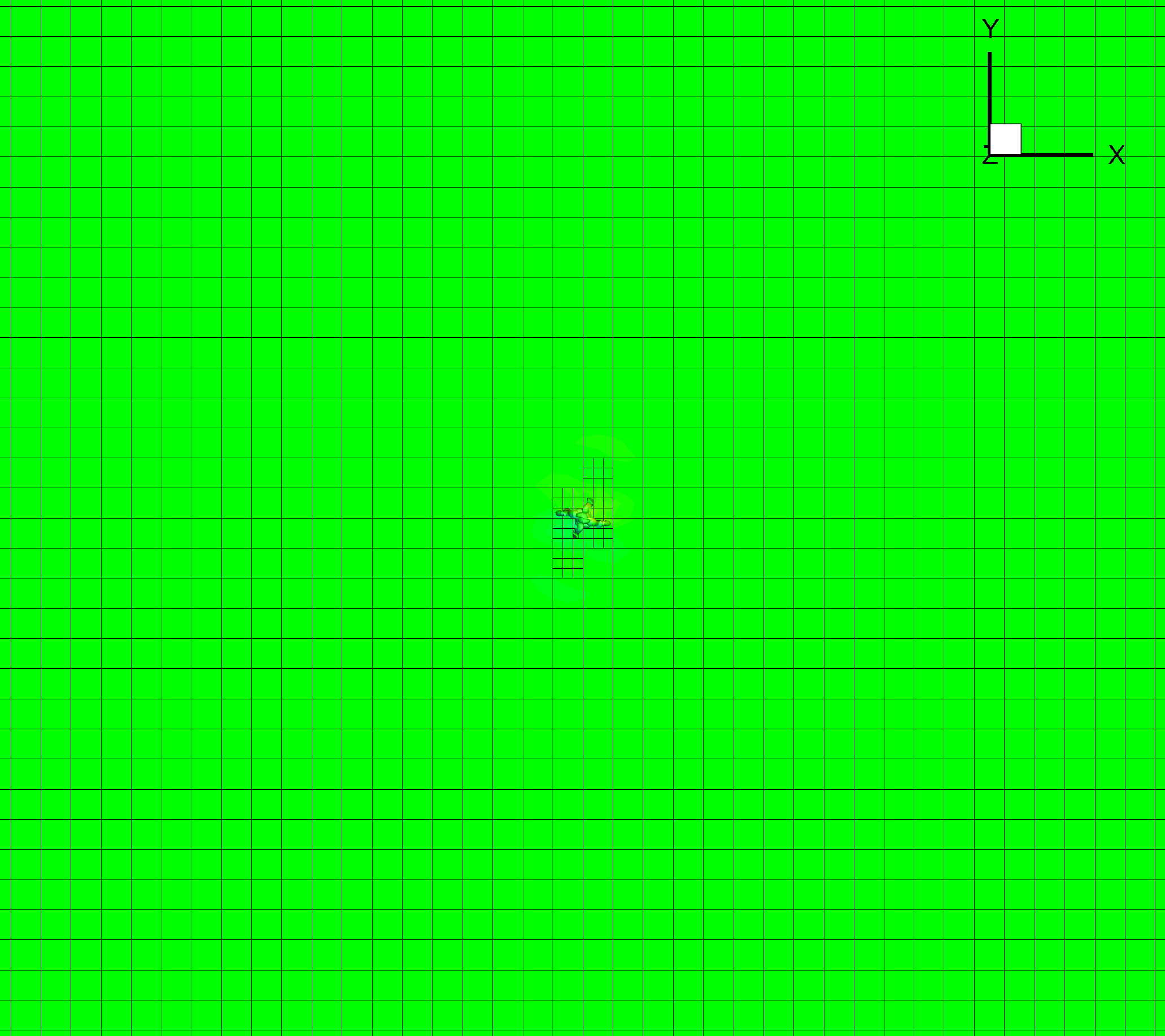}
        \includegraphics[width=0.45\textwidth]{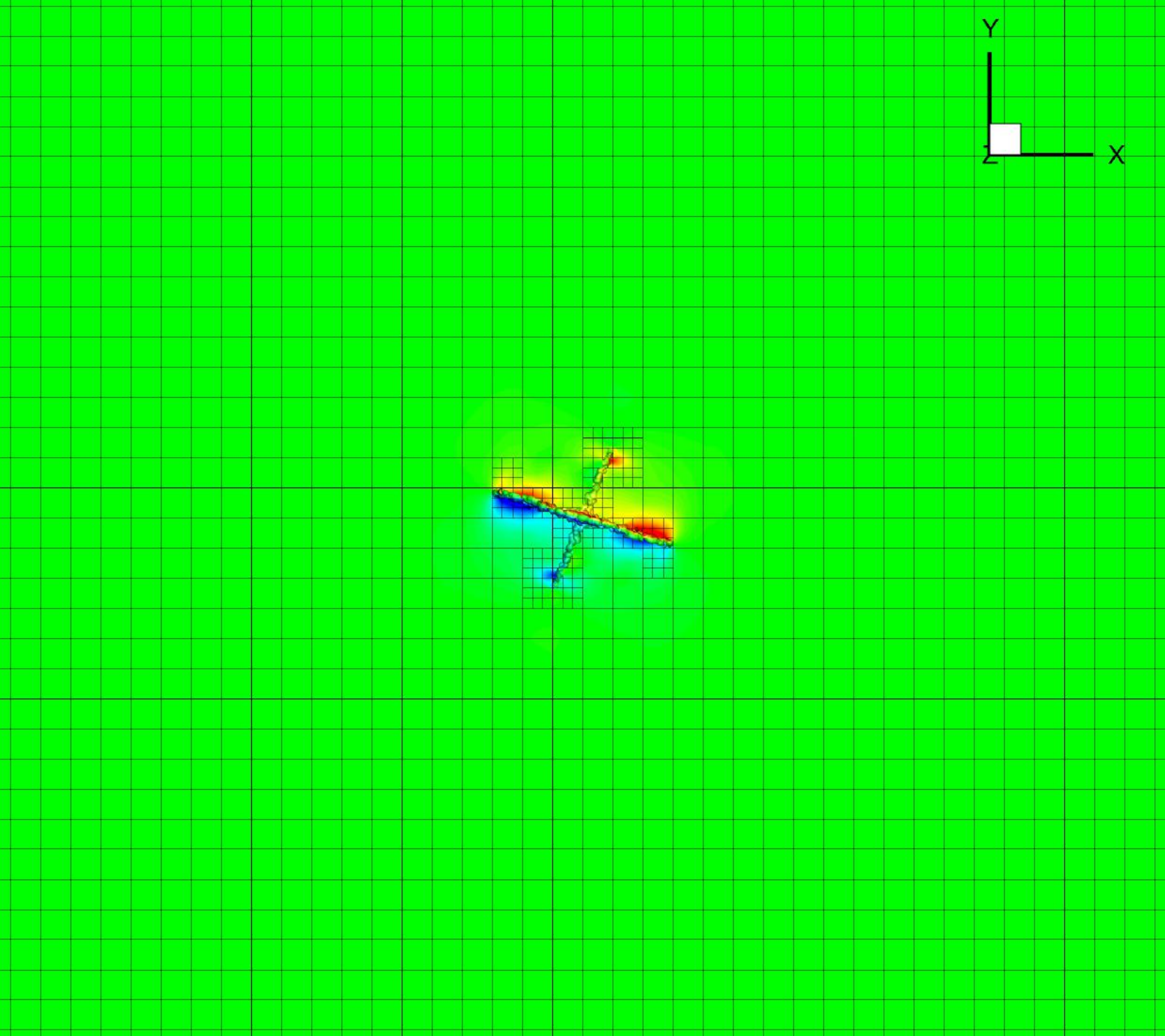}\\[1mm]
        \includegraphics[width=0.45\textwidth]{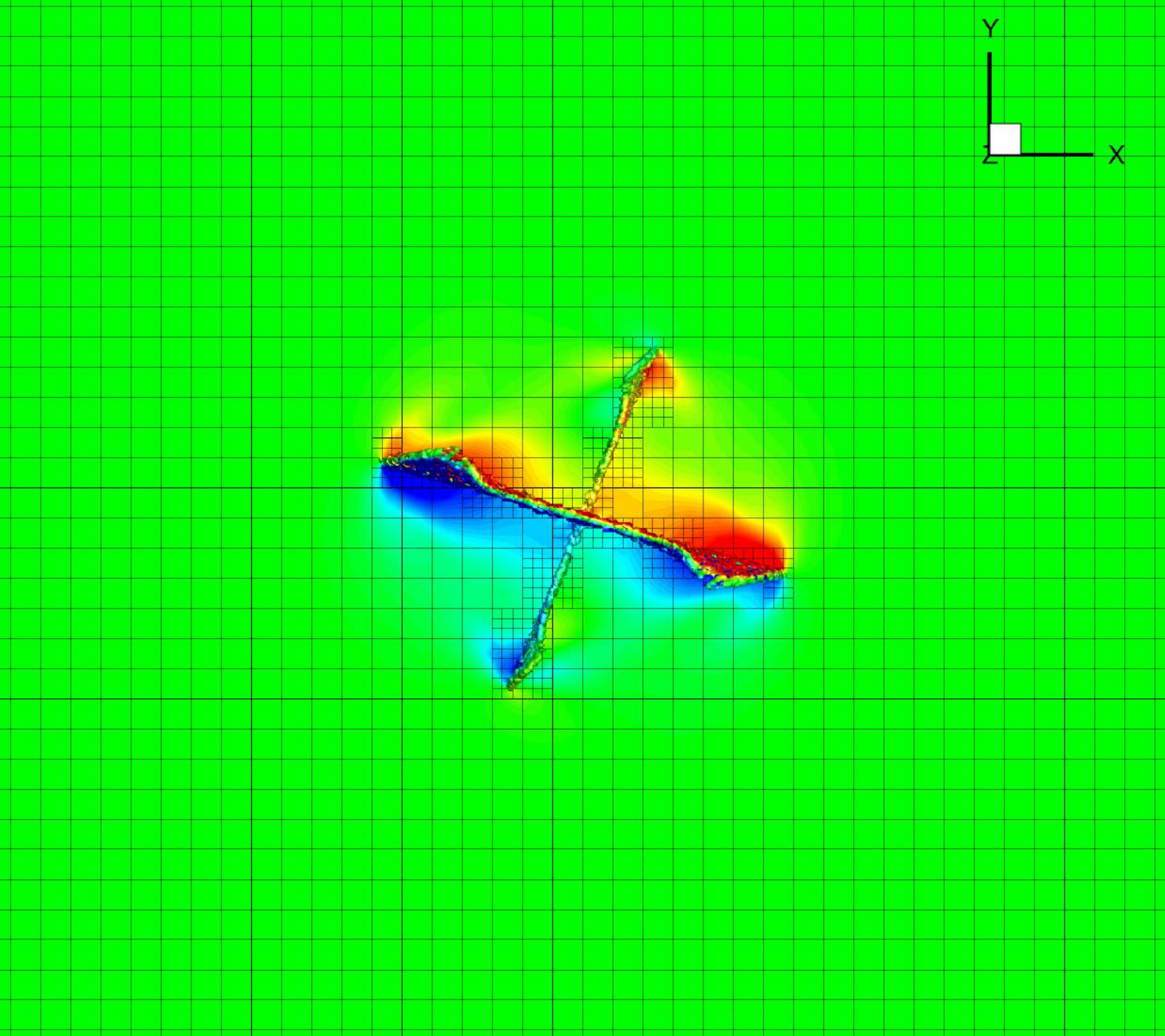}
        \includegraphics[width=0.45\textwidth]{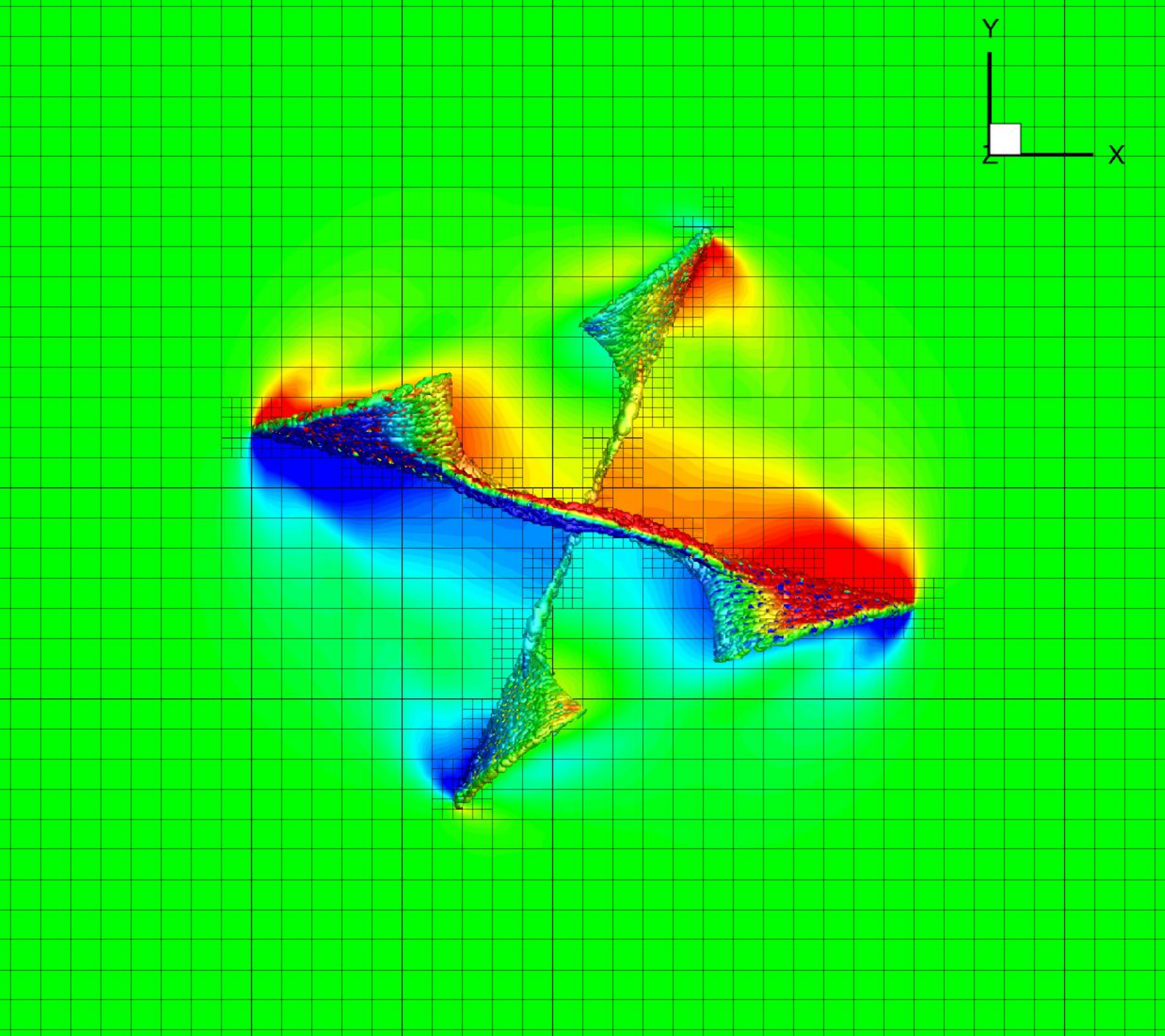}
    \end{center}
    \caption{Resulting rupture front for the $3D$ out of plane rupture colored using the horizontal velocity $u$. }
    \label{3Drupture1}
\end{figure}

We can obviously apply the same algorithm also in three space dimensions. 
Also in this case we take advantage of dynamic adaptive mesh refinement and 
the \textit{a posteriori} subcell FV limiter. The setup of the test is the same 
as in the previous section, but on a three  dimensional domain
 $\Omega=[-5,5]\,\up{km}\times[-5,5]\,\up{km}\times[-2.5,2.5]\,\up{km}$ covered 
 with $40 \times 40 \times 20$ elements on the coarsest mesh level plus one AMR 
refinement level and polynomial approximation degree $N=M=2$. 
 We impose an  initial load of $\sigma_{yy}=-120\,\up{MPa}$ and $\sigma_{xy}=70\,\up{MPa}$ and 
 we pre-damage a zone of $33^3\,\up{m^3}$ in the origin characterized by $\xi=0.5$. 
 The material properties are those associated to a weaker 
 version of ROCK 1, where $Y_0=1.75\,\up{GPa}$. The rupture propagation at 
 times $t=[0.2,0.4,0.6,0.8]$ is shown in Fig. \ref{3Drupture1} and shows initially the formation
  of two planes aligned with the main eigenvectors of the initial stress. 
  Then we can observe also in this test that some three-dimensional bifurcations start. 
\label{sec.rupture3d}

\subsection{Dynamic rupture in complex geometry}  
\label{sec.complex3d}
\begin{figure}[!bp]
    \begin{center}
        \includegraphics[width=0.4\textwidth]{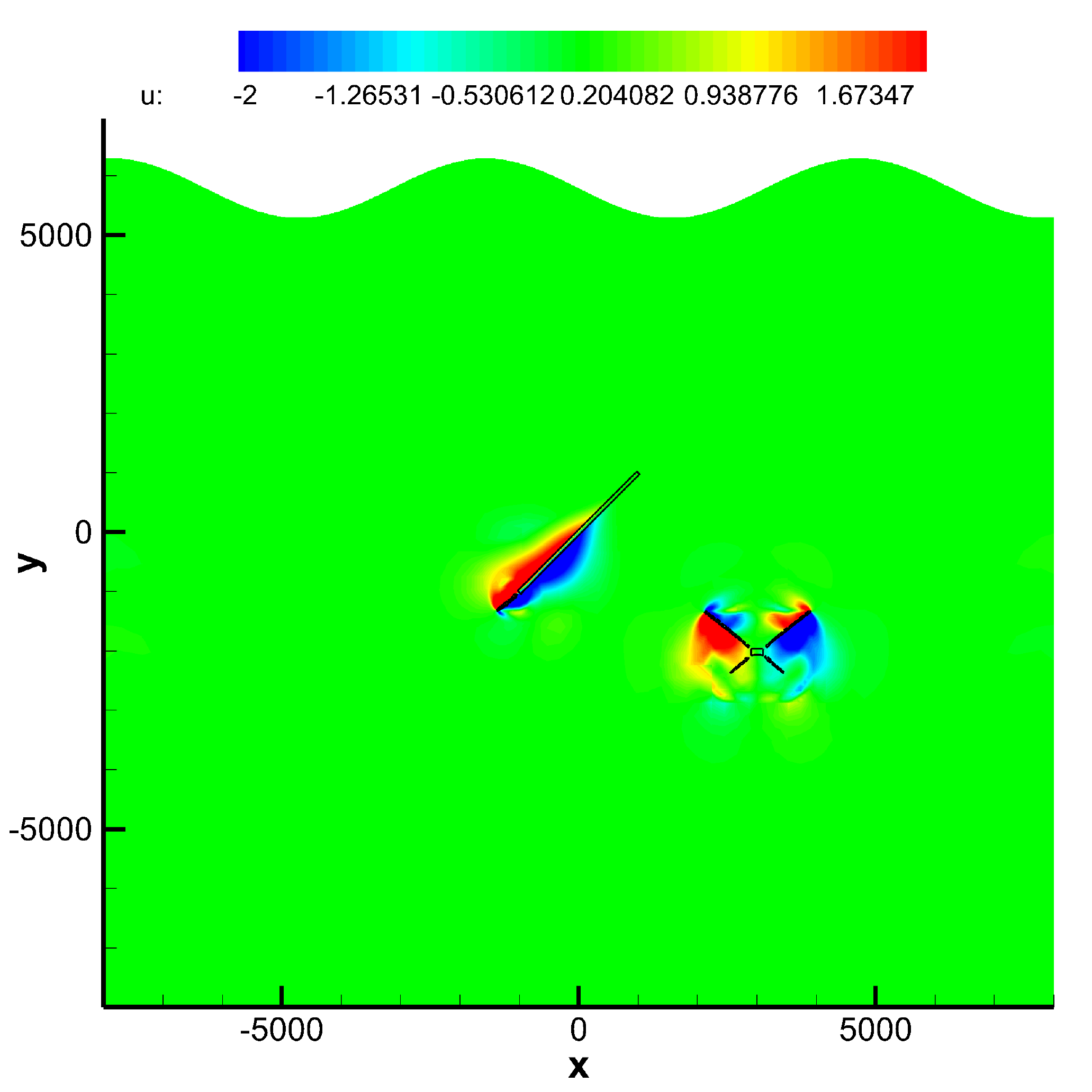}
        \includegraphics[width=0.4\textwidth]{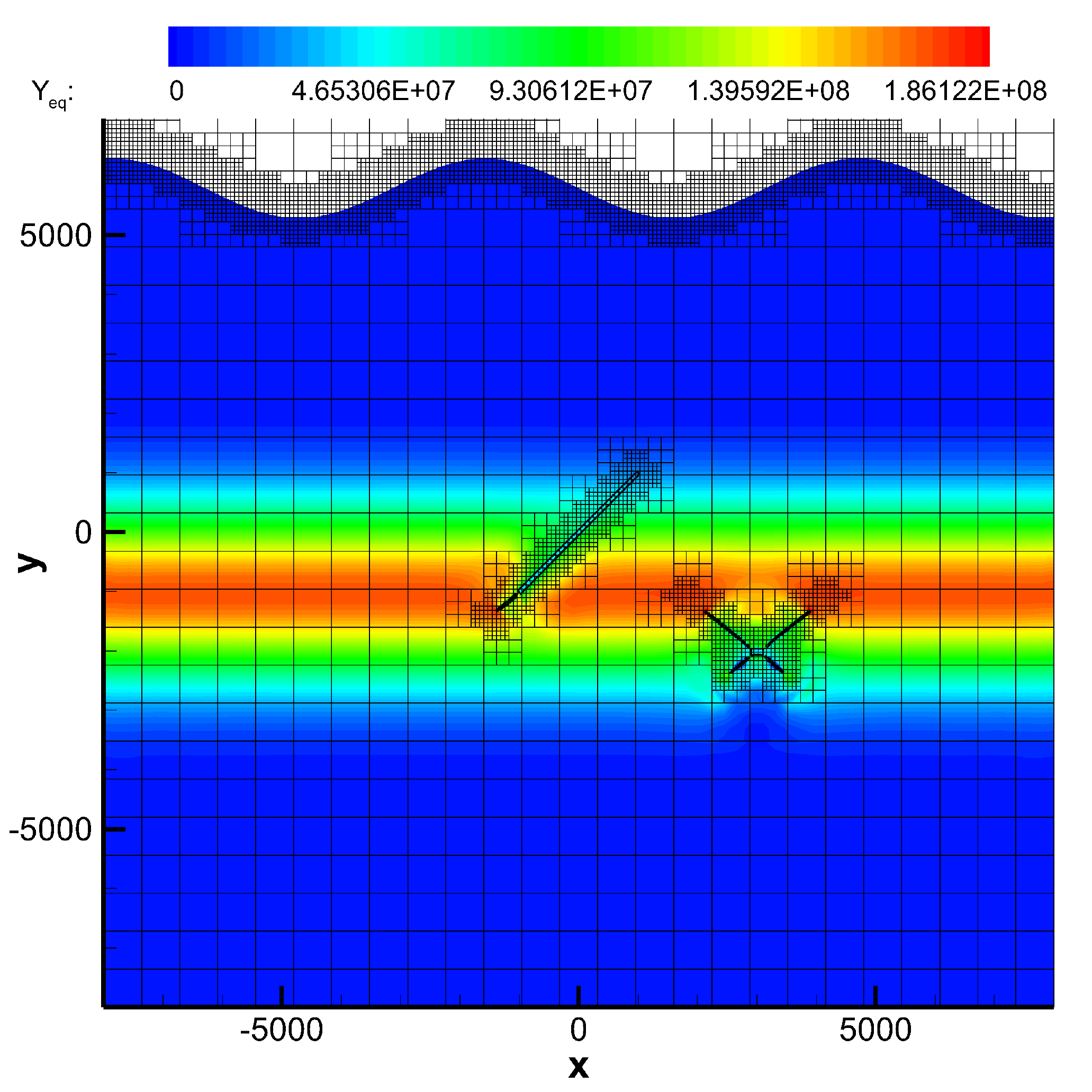} \\
        \includegraphics[width=0.4\textwidth]{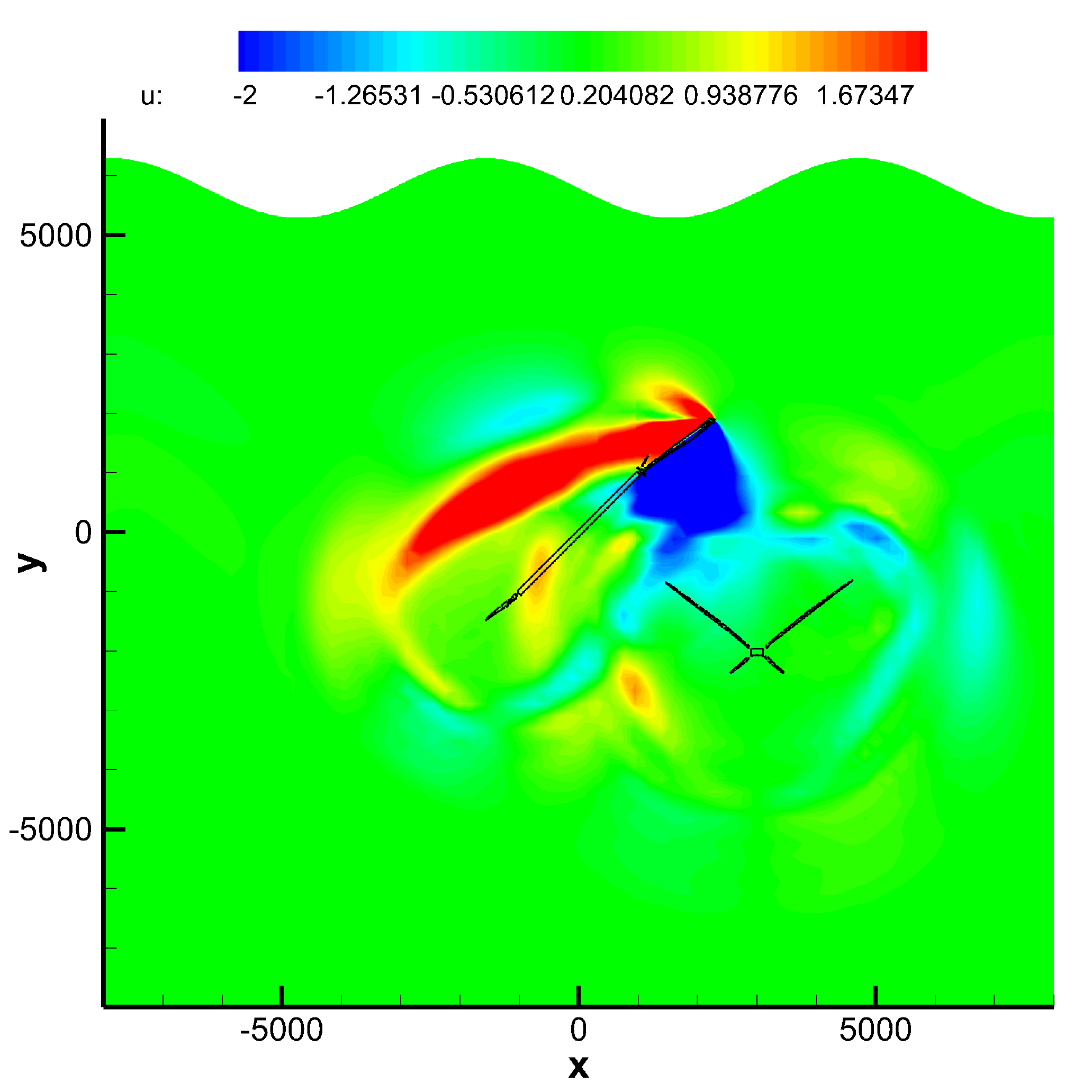}
        \includegraphics[width=0.4\textwidth]{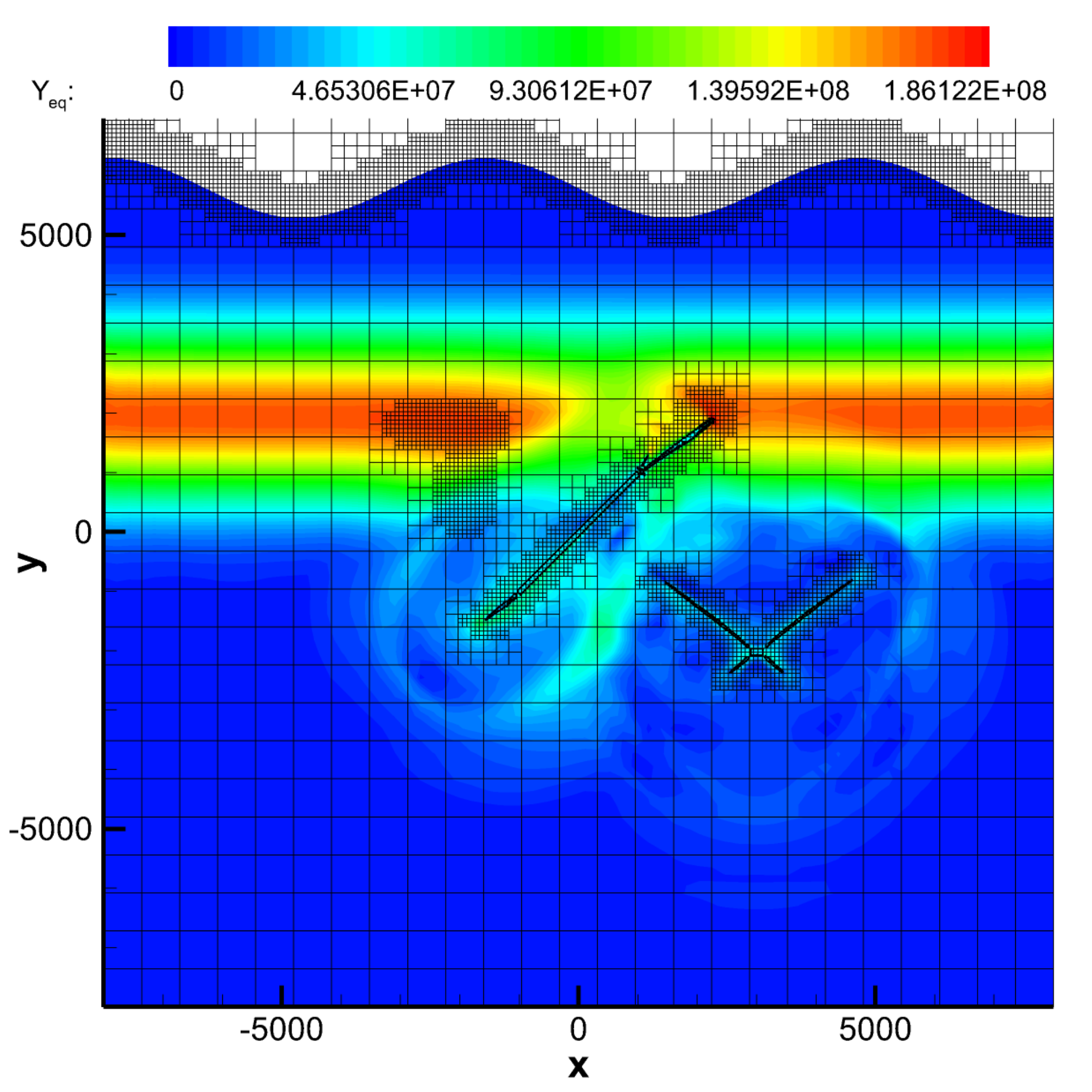} \\
        \includegraphics[width=0.4\textwidth]{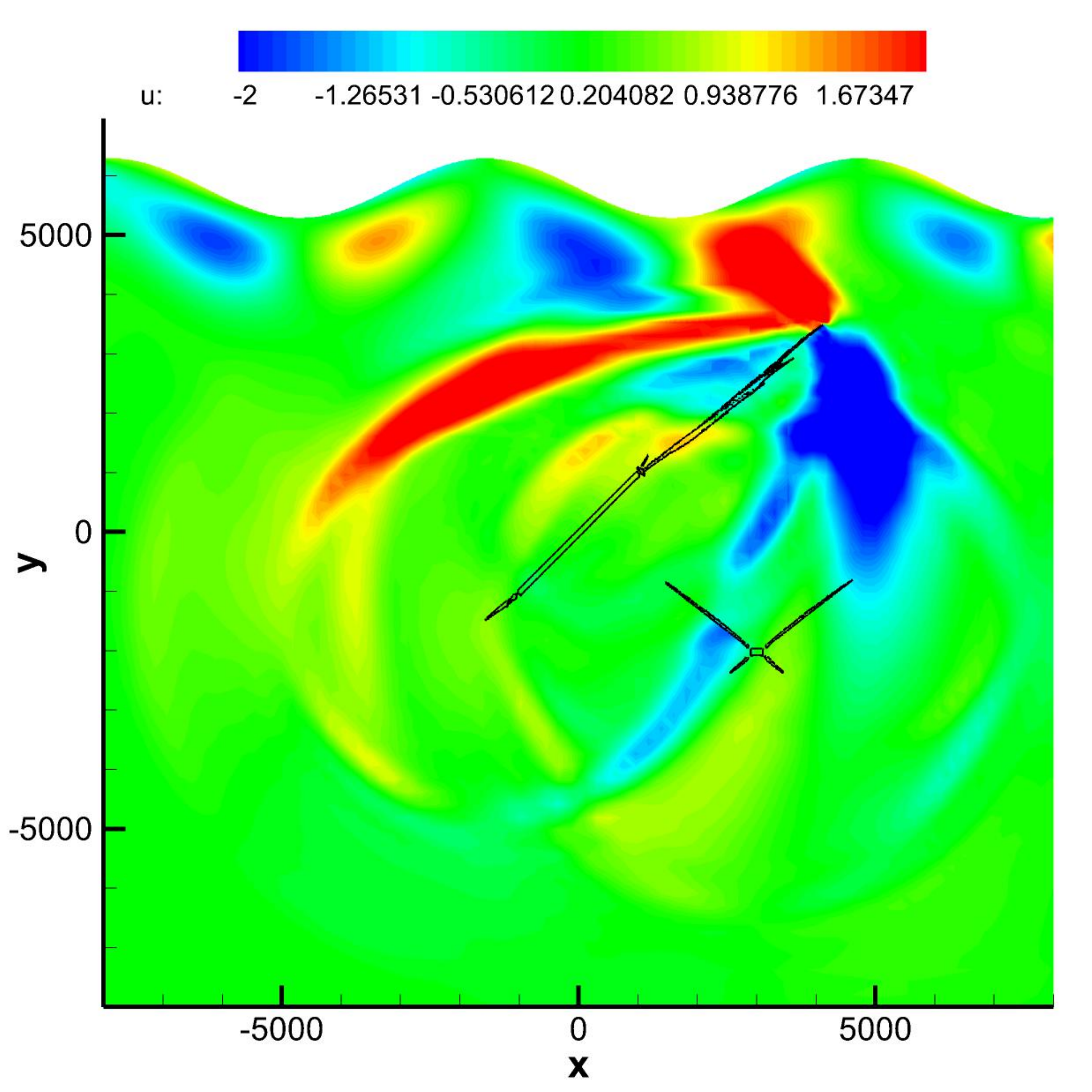}
        \includegraphics[width=0.4\textwidth]{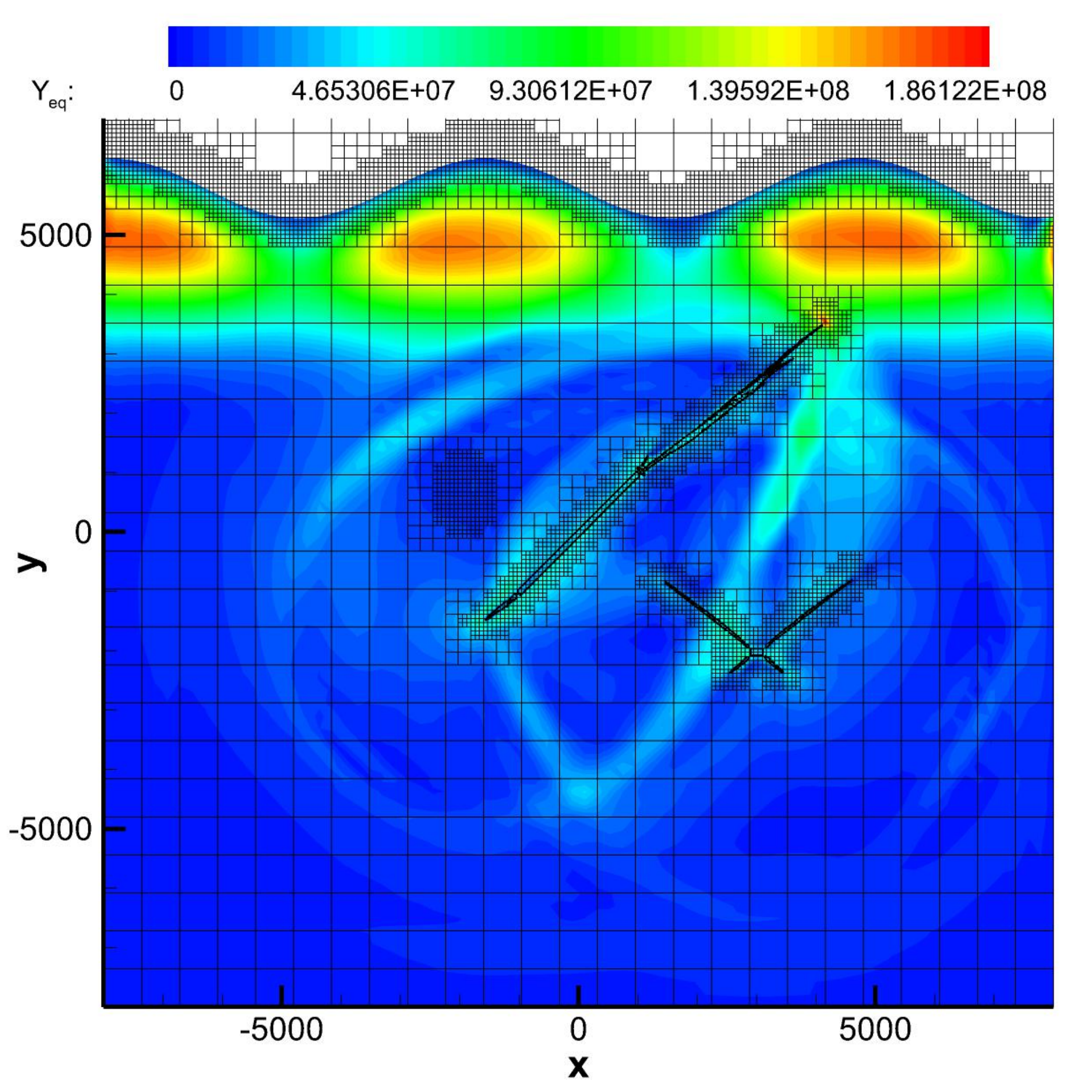} 
    \end{center}
    \caption{P-wave hitting a complex geometry with free crack generation at times, from top to bottom, $t=0.5,1.0$ and $1.5\,\up{s}$. In the left column we report the velocity component $u$ that is only generated due to sliding and complex geometry, in the right column we report the equivalent stress value that depict the traveling p-wave and the interaction with the underground structure. }
    \label{2DCG}
\end{figure}

As a last showcase, in this section we combine dynamic rupture as nonlinear coupling of fracture and wave propagation in complex geometry in 2D and 3D.  
The computational domain is taken as $\Omega=[-8000,8000]^2$ and a level zero mesh with $25 \times 25$ elements is used, in combination with a polynomial approximation degree of $N=M=3$.
Two AMR refinement levels are allowed. As initial condition we set a p-wave traveling in the $y$-direction, namely 
 \begin{eqnarray}
\boldsymbol{\sigma} &=& + \epsilon \, R_\sigma \cdot \mbox{exp}\left({-\frac{1}{2}\frac{(x-x_0)^2}{\delta^2}}\right) \nonumber \\
v&=& -\epsilon \, c_p \cdot \mbox{exp}\left({-\frac{1}{2}\frac{(x-x_0)^2}{\delta^2}}\right) \nonumber \\
u&=& w = 0
\end{eqnarray}
where $R_\sigma=(\lambda,\lambda+2\mu,\lambda,0,0,0)$, $x_0=-4000$, $\delta=1000.0$ and $\epsilon=1.79\cdot 10^{8}/2\mu$. The material is assumed to be ROCK 1, with an equivalent stress simply defined through the von Mises stress. The complex geometry of the free surface of the solid material is described by the following profile:
\begin{equation}
    y_s(x)=6000-500\sin(10^{-3}x).
\end{equation}
The diffuse interface function $\alpha$ is then assumed to be $\alpha=1$ below the surface ($y<y_s$) and $\alpha = 0$ above, i.e. $y>y_s$. The step function in $\alpha$ is then appropriately smoothed with a characteristic smoothing length of $d=300$, see \cite{Tavelli2019} for details. Finally, we set $\xi=1$ in a point placed in  $(x,y)=(3000,-2000)$ and on a line with inclination $\theta=\pi/4$ and horizontal length of $2000\sqrt{2}$. In Figure \ref{2DCG} we report the evolution of the traveling p-wave in terms of equivalent stress (right column) and horizontal velocity (left column). We note how the crack generates and propagates away from the pre-damaged point in a similar way as observed in section \ref{sec.doublecoin}. On the damaged line we see the dislocation through velocity gradients and some bifurcation effects on the main crack line. 

\begin{figure}[!bp]
    \begin{center}
        \includegraphics[width=0.45\textwidth]{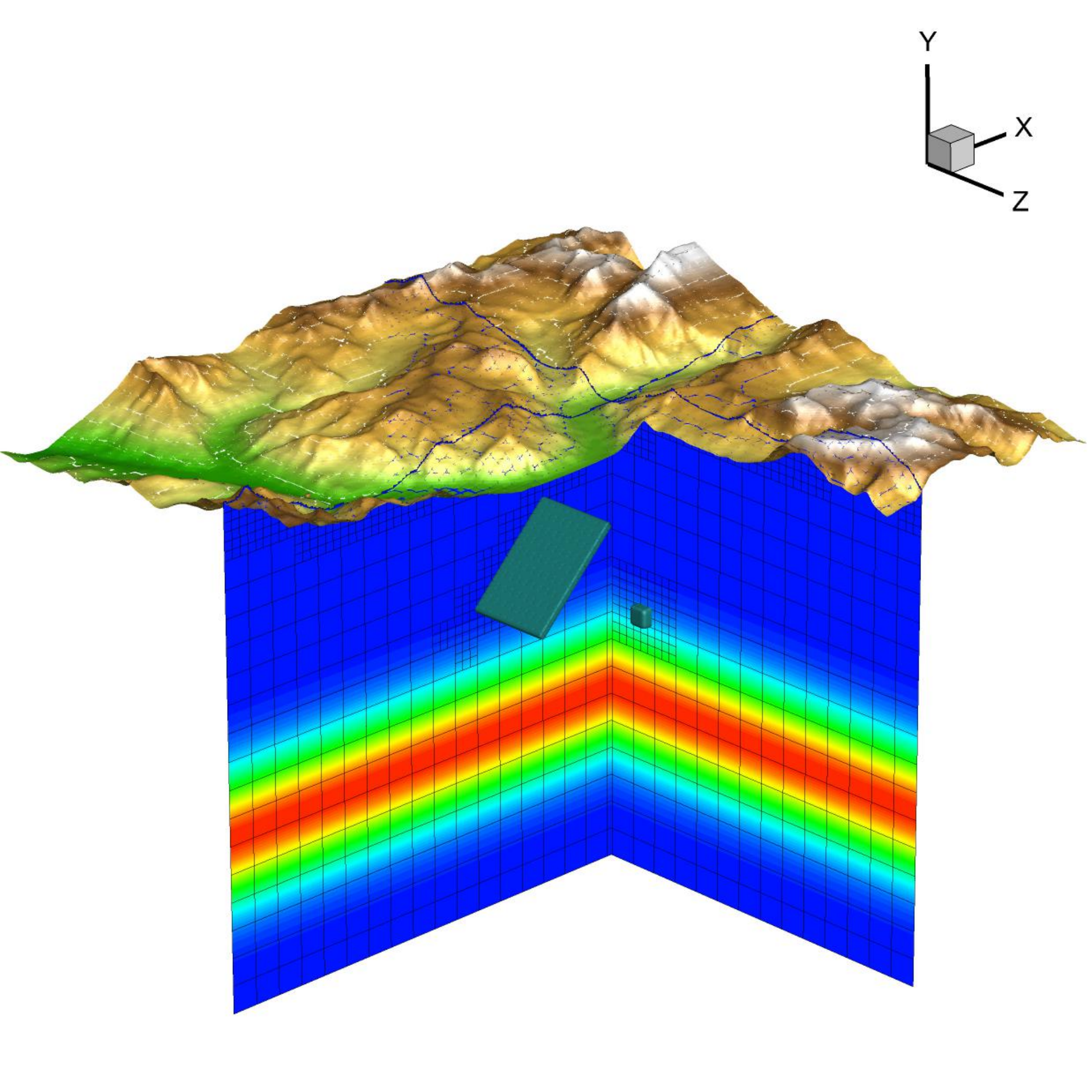}
        \includegraphics[width=0.45\textwidth]{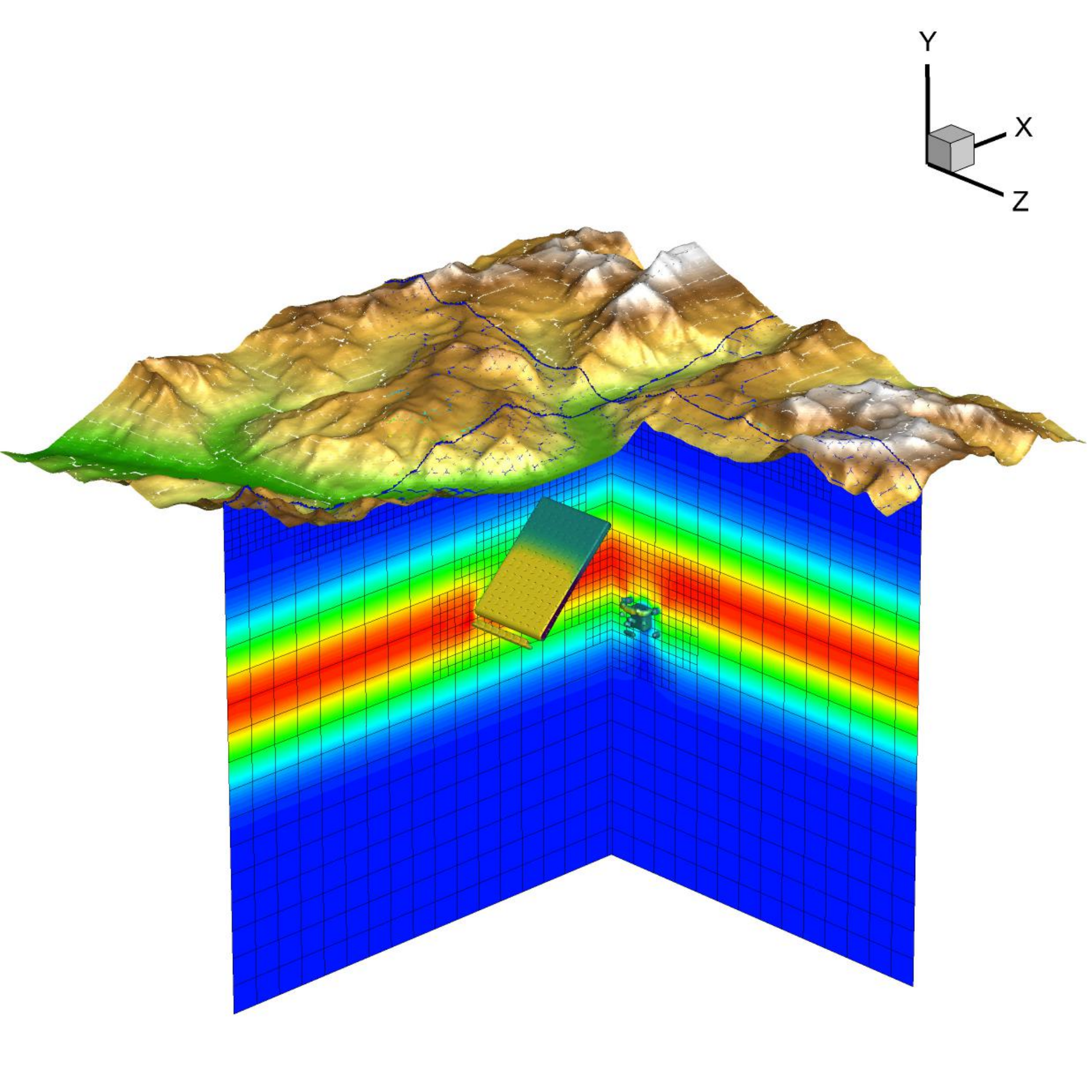} \\
        \includegraphics[width=0.45\textwidth]{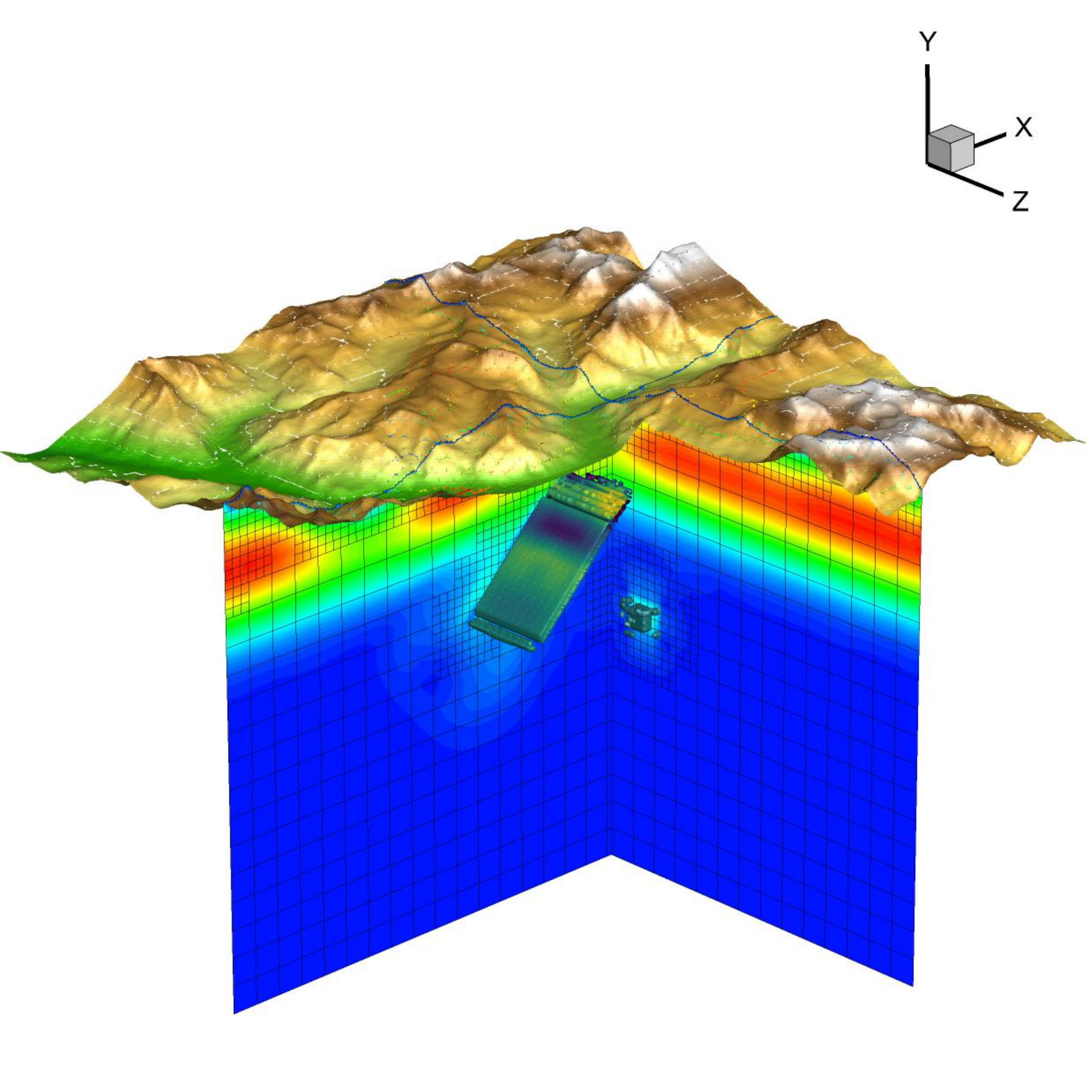}
        \includegraphics[width=0.45\textwidth]{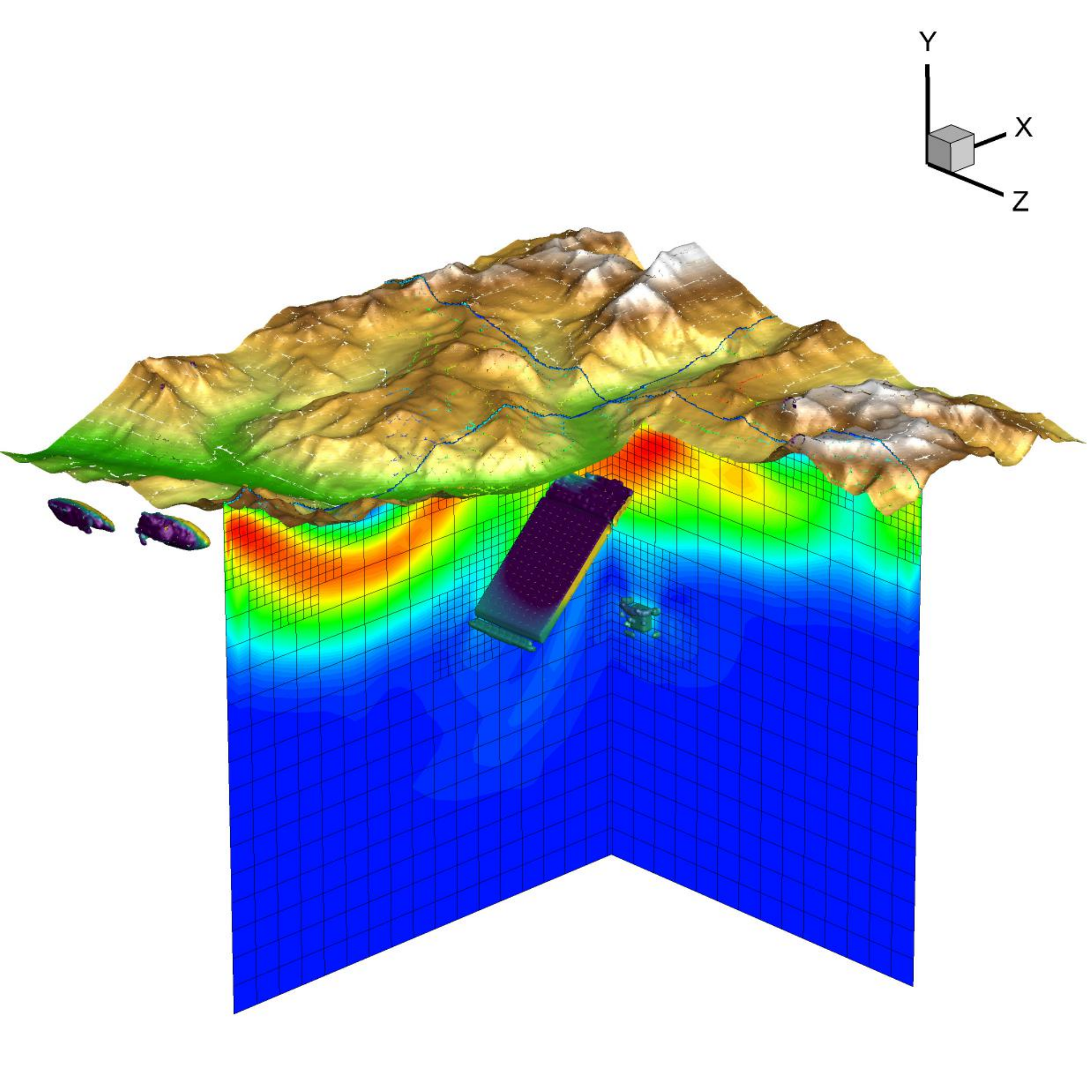}
    \end{center}
    \caption{P-wave propagating in a $3D$ complex geometry, dynamically triggering spontaneous crack generation across a predefined weak fault zone at times $t=[0.0,0.5,1.0,1.4]$ 
        from top left to bottom right. The mountain shape is obtained by extracting the isosurface $\alpha=0.5$, whose values are determined from a real DTM data of extent $16 km$ in 
        each direction, centered in $\vec{x}=[4456.397222711, 2596.544914552]$.}
    \label{3DCG}
\end{figure}

\begin{figure}[!bp]
    \begin{center}
        \includegraphics[width=0.45\textwidth]{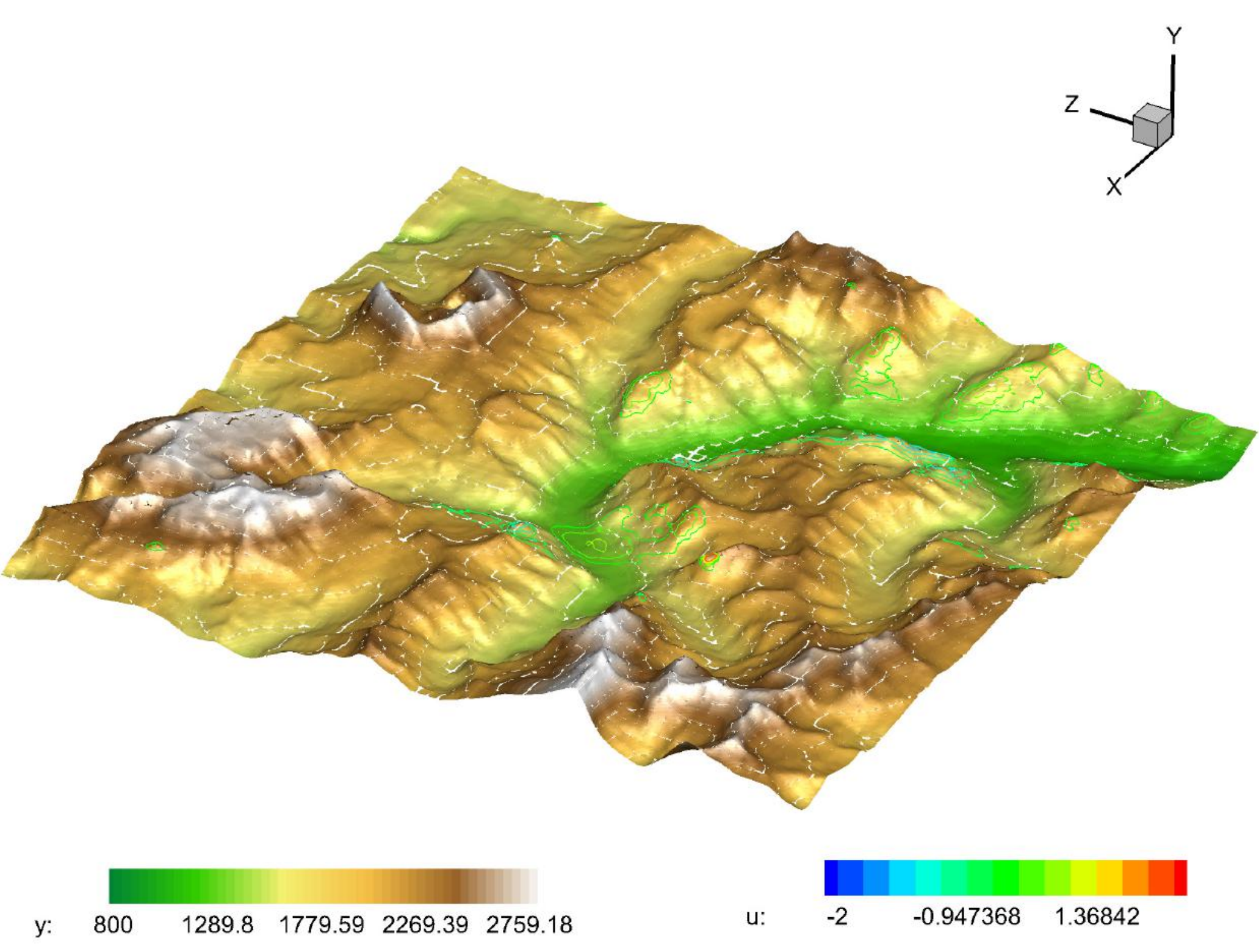}
        \includegraphics[width=0.45\textwidth]{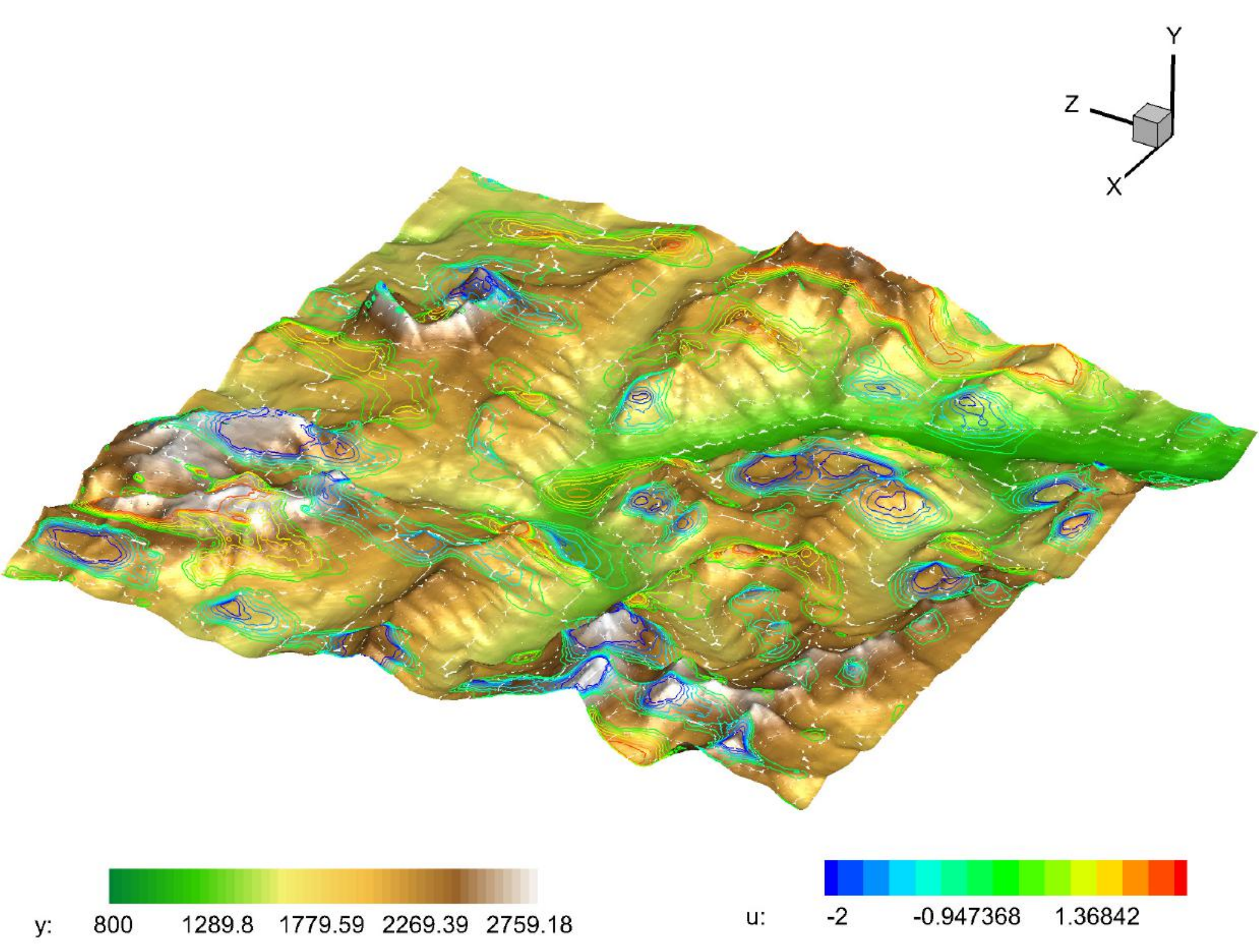}
    \end{center}
    \caption{Velocity component $u$ computed on the surface at the beginning of the impact $t=1.0$ (left) and after the impact at $t=1.4$ (right). We underline that the component $u$ is generated only by sliding and complex geometry shape as for the $2D$ case.}
    \label{3DCG_2}
\end{figure}

\begin{figure}[!bp]
    \begin{center}
        \includegraphics[width=0.9\textwidth]{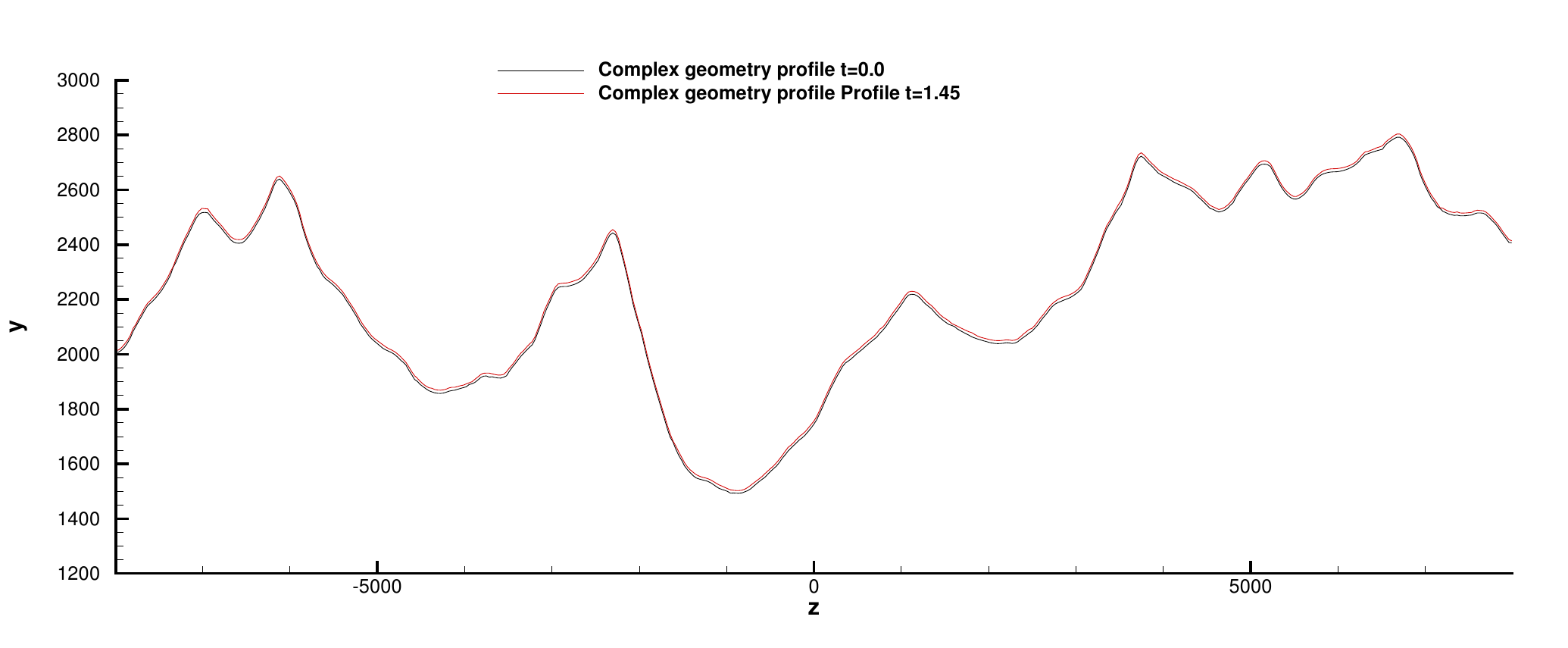}
    \end{center}
    \caption{Elevation profile through the isoline $z \approx 4300$ extracted from the isosurface $\alpha=0.5$ at times $t=0.0$ and $t=1.4$ shows the displacement generated by the impact of the p-wave with the surface due to the transport contribution in the equation for the color function $\alpha$.}
    \label{3DCG_3}
\end{figure}

We now perform a similar test in three-space dimensions. Here the domain $\Omega=[-8000,8000]^3$ is covered with an initially regular Cartesian mesh of  
$25 \times 25 \times 25$ elements and polynomial approximation degree $N=M=2$. 
Then, one level of AMR refinement is allowed. 
For the definition of the complex geometry we extract a zone of 
size $(16\,\up{km})^3$ from real DTM data centered in UTM coordinate  $\vec{x}=[4456.397222711, 2596.544914552]\,\up{km}$, 
which is within the Trentino-S\"udtirol region in Northern Italy. As material we consider a weaker 
version of ROCK 1, where $Y_0=1.75\,\up{GPa}$. The pre-damaged fault line is then extended 
in the $z-$ direction for a total size of $2\,\up{km}$. The time evolution of the p-wave, 
hitting first the pre-damaged zone and then the complex geometry is reported in 
Figure \ref{3DCG}. The dynamic stresses transferred by the propagating waves cause complex crack propagation 
across the predefined weak fault zone. In Figure \ref{3DCG_2} we report the iso-lines of the horizontal  
velocity component $u$ on the complex geometry. We underline here that the $u$ component 
is generated due to the complex topology and underground structure. Finally we report the 
profile $\alpha=0.5$ at position $\{z \approx 4300\}$ before and after the impact in 
Figure \ref{3DCG_3}, where we can observe a small displacement. As already observed 
in \cite{Tavelli2019}, this is possible because according to Eqn. \eqref{eqn.alpha}, the material 
parameters and the solid volume fraction function are free to move according to the local velocity field. 

\section{Conclusion} 
\label{sec:conclusion}

In this paper we have extended the Godunov-Peshkov-Romenski (GPR) model of continuum mechanics in Eulerian coordinates \cite{GodRom1972,Rom1998,PeshRom2014,DumbserPeshkov2016} by an additional dynamic process that describes material failure. More precisely, the extension consists in the addition of a scalar material \textcolor{black}{damage variable} $\xi \in [0,1]$, where $\xi=0$ indicates the fully intact material and $\xi=1$ denotes the fully damaged material. The parameter is governed by an advection-reaction equation, with a highly nonlinear reaction source term that is to a certain extent similar to the reaction source terms that are used to describe reaction kinetics in chemistry. In the present model, the reaction source term depends essentially on the ratio of the equivalent stress (e.g. the von Mises stress) to the yield stress. The underlying ideas go back to previous work of Resnyansky et al. and Romenski \cite{Resnyansky2003,Romenskii2007}, but in the present paper the GPR model with material damage has for the first time solved with high order schemes on space-time adaptive AMR meshes and has furthermore been combined with the diffuse interface approach introduced in \cite{Tavelli2019,FrontierADERGPR}, which allows an easy description of moving elastic bodies on space-time adaptive Cartesian meshes simply via the evolution of an additional scalar solid volume fraction function. The reaction source terms in the governing equation for the \textcolor{black}{damage variable} are extremely stiff, and can in general not be properly integrated in time via a simple implicit Euler method. For this purpose, in this paper we rely on simple operator splitting, combined with a very efficient and accurate exponential time integrator. The proposed time integrator has been carefully validated against standard software packages for stiff ODE. 

Via numerical experiments we have shown that the proposed mathematical model is able to describe \textit{brittle material} as well as \textit{ductile material} if the model parameters are properly chosen. The model is also able to produce \textit{rate-dependent} stress-strain diagrams, and also contains all the necessary ingredients to describe the phenomenon of \textit{material fatigue}. To the best knowledge of the authors, this is the first time that fatigue has been studies in the context of the GPR model of continuum mechanics.      

We have shown numerical results for some test problems with material failure, for which experimental reference data were available. Last but not least, we have also provided a first simplified proof of concept that the computational approach proposed in this paper can at least in principle be also used in the context of describing nonlinear seismic wave propagation and subsequent dynamic rupture processes in complex 3D geometries at a regional scale.    

In the future we plan to use the methodology developed in this paper more extensively for the numerical simulation of nonlinear wave propagation and dynamic rupture processes in computational seismology. \textcolor{black}{We also want to investigate in more detail the rheologies of broken brittle material, which should behave more like a granular fluid rather than an elastic solid. In this case, the correct asymptotic behaviour of the relaxation time $\tau_1$ and of the shear sound speed $c_s$ in the limit $\xi \to 1$ will be of fundamental importance. Further research on this topic will be carried out in the future. }

%=============================================================================
%==========    A C K N O W L E D G M E N T S
\section*{Acknowledgments}
The research presented in this paper has been financed by the European Union's Horizon 2020 Research and  
Innovation Programme under the project \textit{ExaHyPE}, grant agreement number no. 671698 (call 
FETHPC-1-2014). 
Results by E.R. obtained in Sects. 2 and 4 were supported by the Russian Science Foundation grant
(project 19-77-20004).
S.C. acknowledges the financial support received by
the Deutsche Forschungsgemeinschaft (DFG) under the project 
\textit{Droplet Interaction Technologies (DROPIT)}, grant no. GRK 2160/1.
A.-A.G. acknowledges additional funding from the  European Union's Horizon 2020 research and innovation programme under the projects \textit{TEAR}, 
grant agreement No. 852992 (ERC-StG) and \textit{ChEESE}, grant agreement No. 823844.
M.D. also acknowledges the financial support received from the Italian Ministry of Education, University and Research (MIUR) 
in the frame of the Departments of Excellence Initiative 2018--2022 attributed to DICAM of the University of Trento 
(grant L. 232/2016) and in the frame of the PRIN 2017 project \textit{Innovative numerical methods for evolutionary partial differential equations and applications}. M.D. has also received funding from the University of 
Trento via the  \textit{Strategic Initiative Modeling and Simulation}. M.D. is member of the GNCS-INdAM group. 

%=============================================================================
%==========  B I B L I O G R A P H Y
%\input{Bib}
% \section*{References}
\bibliographystyle{plain}
\bibliography{./biblio}
%=============================================================================

\appendix 

\clearpage

\section{Initialization of the distortion field $A$ from a given stress tensor}
\label{sec.init.A}
A peculiar feature of the GPR model is that we evolve the distortion field $\mathbf{A}$ instead 
of the stress tensor $\boldsymbol{\sigma}$, or the strain tensor $\boldsymbol{\epsilon}$. 
The symmetric stress tensor (6 independent components) can be easily computed from 
$\mathbf{A}$ (9 independent components) using the definition $\sigma_{ik} = -\rho A_{ji} E_{A_{jk}}$. For the opposite direction $\boldsymbol{\sigma} \rightarrow  \mathbf{A}$, we first need to remove the rotational degree of freedom contained in $\mathbf{A}$.
In order to do that we write the matrix $\mathbf{A}$ as $\mathbf{A} = \mathbf{L}_d \cdot \mathbf{R}(\theta)$ where
$$\mathbf{R}=\mathbf{R}(\theta)=\left( 
\begin{array}{ccc}
\cos(\theta) & -\sin (\theta) & 0 \\
\sin (\theta) & \cos(\theta) & 0 \\
0& 0& 1
\end{array}
\right), 
\qquad 
\mathbf{L}_d=\left( 
\begin{array}{ccc}
s_1 & 0 & 0 \\
s_2 & s_3 & 0 \\
s_4& s_5& s_6
\end{array}
\right) $$
is a rotation matrix and a lower triangular matrix, respectively; $\theta$ is a chosen as the main direction. 
Now taking $E$ as in Eq. \eqref{E1a}--\eqref{E2a}, $\boldsymbol{\sigma}$ and  $p$ are computed according to their definition and then our target is that the quantity  
$$f(s)=\left(
\begin{array}{c}
\sigma_{11}-p \\
\sigma_{22}-p \\
\sigma_{33}-p \\
\sigma_{12} \\
\sigma_{23} \\
\sigma_{13} \\
\end{array}
\right)$$
approaches our initial stress field $\Sigma^0=(\sigma_{11}^0,\sigma_{22}^0,\sigma_{33}^0,\sigma_{12}^0,\sigma_{23}^0,\sigma_{13}^0)$.
In order to do that, we employ a simple Newton algorithm, that requires the computation of the Jacobian matrix
$$  
J_{ij}=\frac{\partial f_i}{\partial s_j} \qquad \forall\, i,j \in [1,\hdots,6]. 
$$
The explicit expression of the Jacobian for any generic angle $\theta$ can be cumbersome, but can be done at the aid of modern computer algebra systems. The algorithm is then defined as follows: start from an initial guess vector $\vec{s}^0=(s_1^0 \ldots s_6^0)$, then $s^{k+1}=s^k-J^{-1}F(s^k)$ up to convergence of $F(s)=f(s)-\Sigma^0 = 0$.

\end{document}